\documentclass{article}
\usepackage{amsmath}
\usepackage{amsfonts}
\usepackage{amssymb}
\usepackage{amsthm}
\usepackage{bm}
\usepackage{indentfirst}
\usepackage{graphicx}
\usepackage{subfigure}
\usepackage{array}
\usepackage{fullpage}

\DeclareMathAlphabet{\mathsfsl}{OT1}{cmss}{m}{sl}

\newcommand{\PreserveBackslash}[1]{\let\temp=\\#1\let\\=\temp}
\newcolumntype{C}[1]{>{\PreserveBackslash\centering}p{#1}}
\newcolumntype{R}[1]{>{\PreserveBackslash\raggedleft}p{#1}}
\newcolumntype{L}[1]{>{\PreserveBackslash\raggedright}p{#1}}

\numberwithin{equation}{section}
\newtheorem{thm}{Theorem}[section]
\newtheorem{lem}[thm]{Lemma}

\theoremstyle{definition}

\newtheorem{rem}[thm]{Remark}

\usepackage{hyperref}
\usepackage{bm}
\usepackage{amsmath}
\usepackage{amsthm}
\usepackage{amsfonts}
\usepackage{graphicx}
\usepackage[usenames,dvipsnames]{xcolor}

\usepackage[utf8]{inputenc} 
\usepackage[T1]{fontenc}    
\usepackage{hyperref}       
\usepackage{url}            
\usepackage{booktabs}       
\usepackage{nicefrac}       
\usepackage{microtype}      
\usepackage{lipsum}		
\usepackage{fancyvrb}
\usepackage{color}
\usepackage{hyperref}
\usepackage{txfonts}
\usepackage{mathrsfs}
\usepackage{amsmath}
\usepackage{amsfonts}
\usepackage{amsthm}
\usepackage{color}
\usepackage{setspace}
\usepackage{float}
\usepackage{url}
\usepackage{multicol}
\usepackage{subfigure}
\usepackage{float}
\usepackage{soul}
\usepackage[all]{xy}
\usepackage{amssymb,mathtools}
\usepackage{stmaryrd}
\usepackage{hyperref}
\usepackage{pdfsync}
\usepackage[T1]{fontenc} 
\usepackage{lmodern}
\usepackage{bm}
\usepackage{multicol}
\usepackage{multirow}

\newcommand{\fal}{\forall}

\renewcommand{\vartheta}{\Theta}

\definecolor{mygreen}{rgb}{0.1,0.75,0.2}

\title{AN ASYMPTOTICALLY COMPATIBLE APPROACH FOR NEUMANN-TYPE BOUNDARY CONDITION ON NONLOCAL PROBLEMS}

\begin{document}
\date{}
\maketitle


\vspace{-1cm}
\noindent 
\textsf{Huaiqian You,}
\textsf{\textit{\small Department of Mathematics, Lehigh University, Bethlehem, PA, USA}}\\
\textsf{Xin Yang Lu,} \textsf{\textit{\small Department of Mathematical Sciences, Lakehead University, Thunder Bay, ON, Canada}}\\
\textsf{Nathaniel Trask,} \textsf{\textit{\small Center for Computing Research, Sandia National Laboratories, Albuquerque, NM, USA}}\\
\textsf{Yue Yu,} \textsf{\textit{\small Department of Mathematics, Lehigh University, Bethlehem, PA, USA}}, \url{yuy214@lehigh.edu}




\begin{abstract}
In this paper we consider 2D nonlocal diffusion models with a finite nonlocal horizon parameter $\delta$ characterizing the range of nonlocal interactions, and consider the treatment of Neumann-like boundary conditions that have proven challenging for discretizations of nonlocal models. %
While existing 2D nonlocal flux boundary conditions have been shown to exhibit at most first order convergence to the local counter part as $\delta\rightarrow 0$, %
we present a new generalization of classical local Neumann conditions that recovers the local case as $O(\delta^2)$ in the $L^{\infty}(\Omega)$ norm. This convergence %
rate is optimal considering the $O(\delta^2)$ convergence of the nonlocal equation to its local limit away from the boundary. We analyze the application of this new boundary treatment %
to the nonlocal diffusion problem, and present conditions under which the solution of the nonlocal boundary value problem converges to the solution of the corresponding local Neumann %
problem as the horizon is reduced. %
To demonstrate the applicability of this nonlocal flux boundary condition to more %
complicated scenarios, we extend the approach to less regular domains, numerically verifying that we preserve second-order convergence for domains with corners. %
Based on the new formulation for nonlocal boundary condition, %
we develop an asymptotically compatible meshfree discretization, obtaining a solution to the nonlocal diffusion equation with mixed boundary conditions that converges with $O(\delta^2)$ convergence. 
\end{abstract}

\vspace{.5cm}

\vspace{.5cm}


\section{Background}

In recent years, there has been great interest in using nonlocal integro-differential equations (IDEs) as a means to describe physical systems, due to their natural ability to describe %
physical phenomena at small scales and their reduced regularity requirements which lead to greater flexibility %
\cite{silling_2000,bazant2002nonlocal,zimmermann2005continuum,emmrich2007analysis,
emmrich2007well,zhou2010mathematical, du2011mathematical,du2016multiscale,podlubny1998fractional,mainardi2010fractional,magin2006fractional,
burch2011classical,du2014nonlocal,defterli2015fractional,lischke2018fractional,du2014peridynamics,antoine2005approximation,dayal2007real,sachs2013priori,chiarello2018global,erbay2018convergence}. 
In particular, nonlocal problems with Neumann-type boundary constraints %
have received particular attention \cite{cortazar2007boundary,cortazar2008approximate,dipierro2014nonlocal,grubb2014local,montefusco2012fractional,barles2014neumann,
du2013nonlocal,du2015integral,du2017peridynamic,tao2017nonlocal,ren2013compact,aksoylu2010results,mengesha2016characterization,zhou2010mathematical} %
due to their prevalence in describing problems related to: interfaces \cite{alali2015peridynamics}, free boundaries, and multiscale/multiphysics coupling problems \cite{littlewood2015strong,seleson2013interface,Yu2018paper,Astorino_Chouly_Fernandez_2009, Badia_Nobile_Vergara_2008}. %
Unlike classical PDE models, in the nonlocal IDEs the boundary conditions must be defined on
a region with non-zero volume outside the surface \cite{cortazar2008approximate,du2013nonlocal,tao2017nonlocal}, %
in contrast to more traditional engineering scenarios where boundary conditions are typically imposed on a sharp co-dimension one %
surface. Therefore, theoretical and numerical challenges arise from how to mathematically impose inhomogeneous Neumann-type boundary conditions properly in the nonlocal %
model. For instance, in the peridynamic theory of solid mechanics \cite{silling_2000,gerstle2007peridynamic,demmie2007approach,askari2006peridynamic,
xu2008peridynamic,weckner2007damage,foster2009dynamic,madenci2018weak,emmrich2015survey,lipton2014dynamic,madenci2016peridynamic,taylor2015two}, the classical description of material deformation locally via a deformation gradient is replaced by a nonlocal interaction described with integral operators. In these models, it has been shown that the careless imposition of traction conditions on the nonlocal boundary induces an unphysical strain energy concentration, leading in turn to the material being softer near the boundary. Such artificial phenomena are referred to in the literature as a "surface" or "skin" effect\cite{ha2011characteristics,bobaru2011adaptive}.

A key feature in the discretization of nonlocal models has been the concept of \textit{asymptotic compatibility}, originally introduced by Tian and Du \cite{tian2014asymptotically}, which %
describes the ability of a nonlocal discretization to recover a corresponding local model as both $\delta$ and a characteristic discretization lengthscale are reduced at the same rate. %
We advocate the development of both nonlocal boundary treatment and discretization with the objective of preserving this limit. In so doing, we ensure that nonlocal models recover a %
well-understood classical limit, avoiding phenomena such as the surface effect. To this end, we introduce here a non-local boundary treatment that is designed to recover the classical theory. %
After rigorously proving that this nonlocal boundary value problem recovers the desired local Neumann problem as $\delta \rightarrow 0$, we have a firm mathematical foundation upon which to %
demonstrate asymptotic compatibility, where we will develop an asymptotically compatible numerical method and demonstrate its high-order convergence and a lack of artificial surface phenomena.

In this paper we study compactly supported nonlocal integro-differential equations (IDEs) with radial kernels. For concreteness, we focus on the nonlocal diffusion equation
\begin{equation}\label{eqn:nonlocaldelta}
 L_{\delta} u_\delta:=-2\int_{B(\mathbf{x},\delta)} J_\delta(|\mathbf{x}-\mathbf{y}|)(u_\delta(\mathbf{y})-u_\delta(\mathbf{x}))d\mathbf{y}=f(\mathbf{x}),\quad \mathbf{x}\in\Omega\subset \mathbb{R}^N,
\end{equation}
although the proposed technique is applicable to more general problems. %
Here $B(\mathbf{x},\delta)$ is the ball centered at $\mathbf{x}$ with radius $\delta$, $u_\delta(\mathbf{x})$ is the solution, $\Omega$ is a bounded and connected domain %
in $\mathbb{R}^N$ ($N=2$), $f(\mathbf{x})$ is given data, and %
the kernel function $J_\delta:\mathbb{R}\rightarrow\mathbb{R}$ is parameterized by a positive horizon parameter $\delta$ which measures the extent of nonlocal %
interaction. We further take a popular choice of $J_\delta$ as a rescaled kernel given by
\begin{equation}\label{eqn:Jdelta}
 J_{\delta}(|\bm{\xi}|)=\dfrac{c}{\delta^{N+2}}J\left(\dfrac{|\bm{\xi}|}{\delta}\right),
\end{equation}
where $J:[0,\infty)\rightarrow\mathbb{R}$ is a nonnegative and continuous function with $\int_{\mathbb{R}^N} J(|\mathbf{z}|) |\mathbf{z}|^2 d\mathbf{z}=N$, $\int_{\mathbb{R}^N} J(|\mathbf{z}|) |\mathbf{z}|^k d\mathbf{z}=o(1)$ %
for all $k\geq 3$. Similar as in \cite{tao2017nonlocal}, we also assume that $J(r)$ is nonincreasing in $r$, strictly positive in $r\in[0,1]$ and vanishes when $r>1$. %
In this work we aim to design a new formulation of Neumann-type constraint for the nonlocal problem \eqref{eqn:nonlocaldelta} with mixed boundary conditions of Dirichlet, Neumann and mixed type, and present a numerical discretization of the resulting problem.

We pose three requirements for this formulation:
\begin{enumerate}
\item The constraint should be a proper nonlocal analogue to the local Neumann-type boundary conditions, so the formulation provides an approximation of physical boundary conditions on a sharp surface. %
\item A boundary value problem given by the nonlocal Neumann-type constraint with the nonlocal diffusion equation \eqref{eqn:nonlocaldelta} should be well-posed. %
Rigorous mathematical analysis on the existence, uniqueness and continuous dependence on data should be addressed for the associated variational problem.
\item The nonlocal Neumann-type boundary value problem should recover the classical Neumann problem  as $\delta \rightarrow 0$, preferably with an optimal convergence rate of $\textit{O}(\delta^2)$ in the $L^{\infty}$ norm. %
\end{enumerate}


In the first part of the paper, we provide analysis of the boundary value problem (BVP) to establish the consistency and well-posedness of the boundary value problem. We establish here second-order convergence on non-trivial geometry, improving upon the first-order, one-dimensional analysis found in the literature.\cite{cortazar2008approximate,du2015integral,tao2017nonlocal}. In the second part of this paper, we will present a new asymptotically compatible meshfree discretization of the proposed nonlocal BVP \cite{bessa2014meshfree,parks2008pdlammps,Trask2018paper}. We pursue an extension of previous work by Trask et al. \cite{Trask2018paper} utilizing an optimization-based approach to meshfree quadrature. This framework is attractive due to its demonstrated ability to achieve high-order asymptotically compatible solutions on unstructured data, which is complementary to the objective of developing boundary conditions consistent for irregular geometries. By introducing the new boundary treatment we will demonstrate improved second-order convergence over the previously demonstrated first order-convergence shown for Neumann problems \cite{Trask2018paper}.

The paper is organized as follows. We first present in Section \ref{sec:flux} a definition of the nonlocal Neumann-type boundary condition and the corresponding nonlocal variational problem, together with the associated nonlocal operator and natural energy space. In Section \ref{sec:l2}, we study the well-posedness of the nonlocal variational problem for convex and sufficiently regular domains. We provide a consistency result for the nonlocal BVP by showing that the weak solution of the proposed nonlocal Neumann-type constrained value problem (denoted as $u_\delta$) converges to the solution of the corresponding classical diffusion problem (denoted as $u_0$) as the interaction horizon $\delta\rightarrow 0$ in the $L^2(\Omega)$ norm. Although the main proof in this section has assumed homogeneous Neumann boundary conditions, a discussion on the extension to inhomogeneous Neumann-type boundary condition is also provided. Furthermore, in Section \ref{sec:linfty} we prove the $O(\delta^2)$ convergence rate of the continuous nonlocal solution $u_\delta$ to $u_0$ in the $L^{\infty}(\Omega)$ norm without extra regularity assumptions on $u_\delta$. Numerical results utilizing a meshfree quadrature rule are presented in Section \ref{sec:meshfree}. To verify the asymptotic compatibility of the combined boundary treatment and the meshfree numerical scheme, in Section \ref{sec:test} we use manufactured solutions to demonstrate the convergence of the discrete model to the local solution as both the discretization length scale $h$ and the nonlocal interaction length scale $\delta\rightarrow 0$. Furthermore, in Section \ref{sec:corner} we extend the approach to domains with corners which indicate that the conclusions of the model problem and the convergence rates extrapolate to nontrivial problems of interest to the broader engineering community. Section \ref{sec:conclusion} summarizes our findings and discusses future research. The appendices include additional technical details on the theoretical analysis.

\section{A Nonlocal Flux Condition for 2D Diffusion Problem}\label{sec:flux}

\begin{figure}[!htb]\centering
 \subfigure{\includegraphics[width=0.35\textwidth]{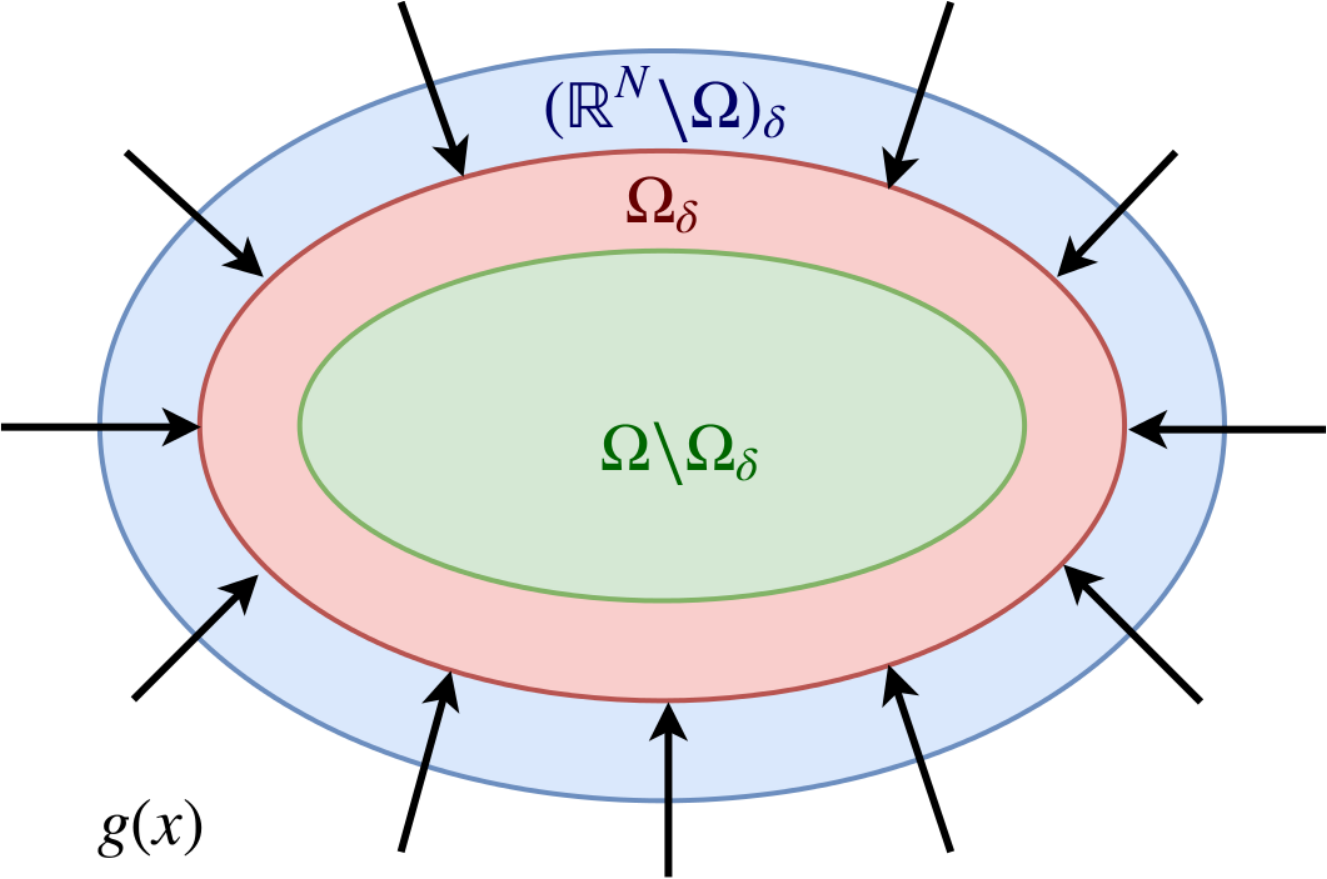}}\hspace*{1in}
 \subfigure{\includegraphics[width=0.2\textwidth]{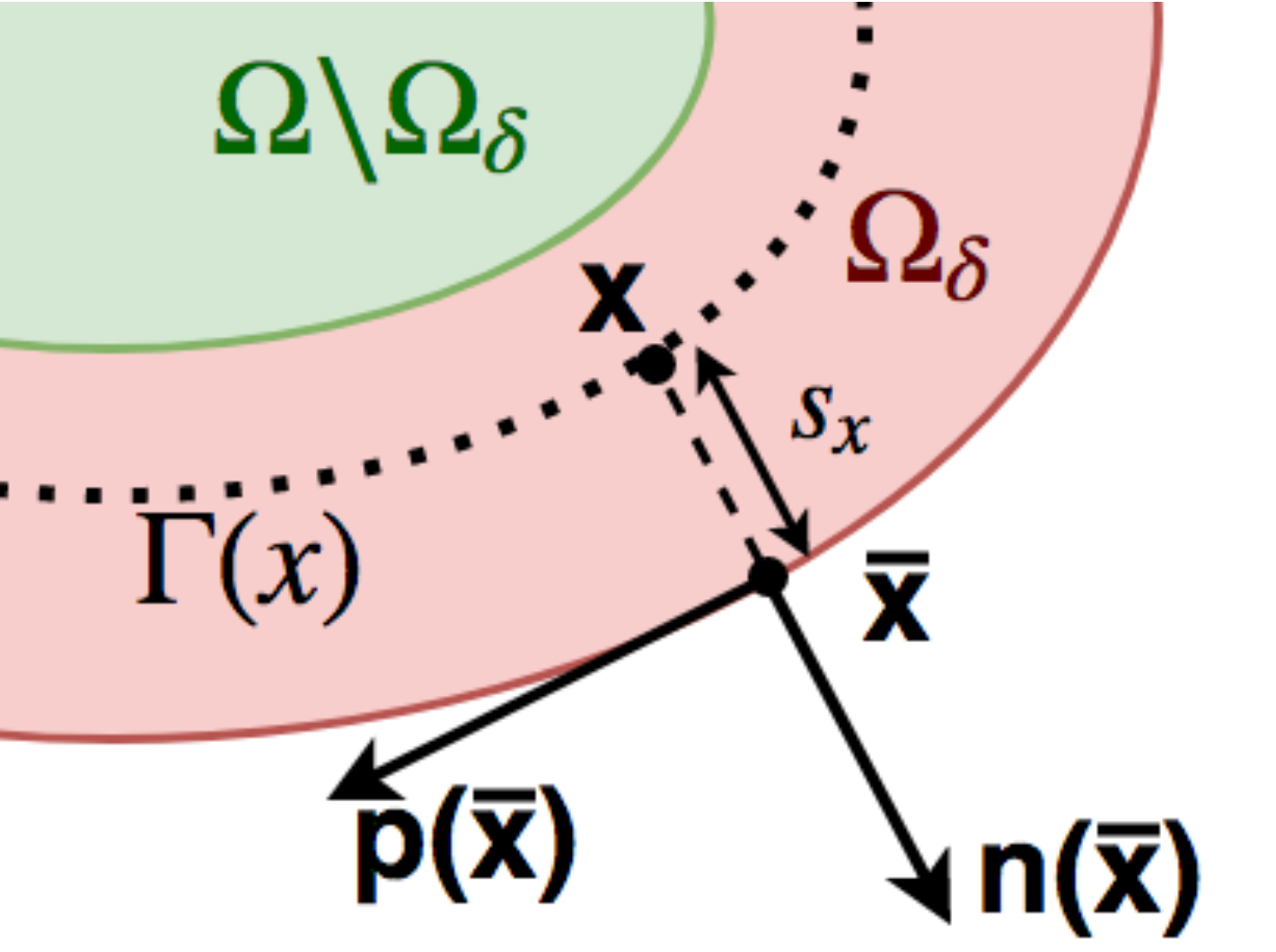}}
 \caption{Left: Notations for the domain, where $\Omega$ is represented by the green and red regions together, and the nonlocal Neumann boundary condition is applied on the red region $\Omega_\delta$. Right: Notations for the projection of point $\mathbf{x}\in \Omega_\delta$, the corresponding unit tangential vector $\mathbf{p}(\overline{\mathbf{x}})$ and the unit normal vector $\mathbf{n}(\overline{\mathbf{x}})$. 
 }
 \label{fig1}
\end{figure}

In this section, we first introduce a nonlocal flux boundary condition, and then provide a corresponding nonlocal variational problem along with the associated energy space for the purpose of analysis. Given that $\Omega\in \mathbb{R}^N$ $(N=2)$ is a bounded, convex, connected and {$C^{3}$} domain, we seek a nonlocal analogue to the local Neumann boundary condition $\dfrac{\partial u}{\partial \mathbf{n}}=g(\mathbf{x})$, $\mathbf{x}\in\partial\Omega$ in the following classical problem
\begin{equation}\label{eqn:localeqn}
 \left\{\begin{array} {cl}
  L_0u_0:=-\triangle u_0=f(\mathbf{x}),&\text{ in }\Omega\\
  \dfrac{\partial u_0}{\partial \mathbf{n}}=g(\mathbf{x}),&\text{ on }\partial \Omega\\\\
  \int_\Omega u_0(\mathbf{x})d\mathbf{x}=0.\\
 \end{array}\right.
\end{equation}
Here $\mathbf{n}({\mathbf{x}})$ is the unit exterior normal to $\Omega$ at ${\mathbf{x}}$. %
Moreover, we will use $\mathbf{p}({\mathbf{x}})$ to represent the unit tangential vector with orientation clockwise to $\mathbf{n}({\mathbf{x}})$. Before introducing our nonlocal formulation, we denote the following notation (see Figure \ref{fig1} for illustration)
\begin{align*}
\Omega_{\delta}:=\{\mathbf{x}\in\Omega|\text{dist}(\mathbf{x},\partial\Omega)<\delta\},\quad
(\mathbb{R}^N\backslash\Omega)_{\delta}:=\{\mathbf{x}\in\mathbb{R}^N\backslash\Omega|\text{dist}(\mathbf{x},\partial\Omega)<\delta\}.
\end{align*}
We further assume sufficient regularity in the boundary that we may take $\delta$ sufficiently small so that for any $\mathbf{x}\in \Omega_{\delta}$, there exists a unique orthogonal projection of $\mathbf{x}$ onto $\partial\Omega$. We denote this projection as $\overline{\mathbf{x}}$. Therefore, one has
$\overline{\mathbf{x}}-\mathbf{x}=s_x\mathbf{n}(\overline{\mathbf{x}})$ for $\mathbf{x}\in \Omega_{\delta}$, where $0<s_x<\delta$. We also assume that for $\mathbf{x}\in \Omega_\delta$, we can find a contour $\Gamma(\mathbf{x})$ which is parallel to $\partial \Omega$. In the following contents, we denote $\mathbf{x}_{l}$ as the point with distance $l$ to $\mathbf{x}$ along $\Gamma(\mathbf{x})$ following the $\mathbf{p}(\overline{\mathbf{x}})$ direction, and $\mathbf{x}_{-l}$ as the point with distance $l$ to $\mathbf{x}$ in the opposite direction. Moreover, we employ the following notations for the directional components of the Hessian matrix of a scalar function $v$:
\begin{align*}
[v({{\mathbf{x}}})]_{pp}:=\mathbf{p}^T(\overline{\mathbf{x}}) \nabla^2v({{\mathbf{x}}}) \mathbf{p}(\overline{\mathbf{x}}),\quad
[v({{\mathbf{x}}})]_{nn}:=\mathbf{n}^T(\overline{\mathbf{x}}) \nabla^2v({{\mathbf{x}}}) \mathbf{n}(\overline{\mathbf{x}}),\quad
[v({{\mathbf{x}}})]_{pn}:=\mathbf{p}^T(\overline{\mathbf{x}}) \nabla^2v({{\mathbf{x}}}) \mathbf{n}(\overline{\mathbf{x}}),
\end{align*}
and the higher order derivative components are similarly defined.

Since $B(\mathbf{x},\delta)\cap(\mathbb{R}^N\backslash\Omega)_\delta \neq \emptyset$ for $\mathbf{x}\in \Omega_{\delta}$, from \eqref{eqn:nonlocaldelta} we have
\begin{align*}
L_{\delta} u_\delta=&-2\int_{B(\mathbf{x},\delta)} J_\delta(|\mathbf{x}-\mathbf{y}|)(u_\delta(\mathbf{y})-u_\delta(\mathbf{x}))d\mathbf{y}\\
=&-2\int_{B(\mathbf{x},\delta)\cap\Omega} J_\delta(|\mathbf{x}-\mathbf{y}|)(u_\delta(\mathbf{y})-u_\delta(\mathbf{x}))d\mathbf{y}
-2\int_{B(\mathbf{x},\delta)\cap (\mathbb{R}^N\backslash\Omega)_\delta} J_\delta(|\mathbf{x}-\mathbf{y}|)(u_\delta(\mathbf{y})-u_\delta(\mathbf{x}))d\mathbf{y},
\end{align*}
hence we need to approximate the integral in $B(\mathbf{x},\delta)\cap (\mathbb{R}^N\backslash\Omega)_\delta$ and obtain a formulation with correction terms. Specifically, we propose the following flux boundary condition for \eqref{eqn:nonlocaldelta}: for $\mathbf{x}\in \Omega_\delta$
\begin{align}
 \nonumber&-2\int_{\Omega} J_{\delta}(|\mathbf{x}-\mathbf{y}|)(u_{\delta}(\mathbf{y})-u_{\delta}(\mathbf{x}))d\mathbf{y}
 -\int_{\mathbb{R}^N\backslash\Omega} J_{\delta}(|\mathbf{x}-\mathbf{y}|)(\mathbf{y}-\mathbf{x})\cdot\mathbf{n}(\overline{\mathbf{x}})%
 (g({\mathbf{x}})+g({\mathbf{y}}))d\mathbf{y}\\
&
 -\int_{\mathbb{R}^N\backslash\Omega} J_{\delta}(|\mathbf{x}-\mathbf{y}|)|(\mathbf{y}-\mathbf{x})\cdot\mathbf{p}(\overline{\mathbf{x}})|^2%
 d\mathbf{y}[ u_{\delta}({{\mathbf{x}}})]_{pp}=f(\mathbf{x}),\label{eqn:formula1}
\end{align}
where the second and third terms aim to provide an approximation for
$$-2\int_{\mathbb{R}^N\backslash\Omega} J_{\delta}(|\mathbf{x}-\mathbf{y}|)(u_{\delta}(\mathbf{y})-u_{\delta}(\mathbf{x}))d\mathbf{y}.$$
Since the boundary condition $g(\mathbf{x})$ is defined only on $\partial\Omega$, the $g(\mathbf{x})$ and $g(\mathbf{y})$ terms in \eqref{eqn:formula1} will be approximated with the following (local) extensions
\begin{align}
 \nonumber&g(\mathbf{x})\approx g(\overline{\mathbf{x}})- (\mathbf{x}-\overline{\mathbf{x}})\cdot\mathbf{n}(\overline{\mathbf{x}})f({{\mathbf{x}}})-%
 (\mathbf{x}-\overline{\mathbf{x}})\cdot\mathbf{n}(\overline{\mathbf{x}}) [ u_{\delta}({{\mathbf{x}}})]_{pp},\\
 \nonumber&g(\mathbf{y})\approx g(\overline{\mathbf{x}})-(\mathbf{y}-\overline{\mathbf{x}})\cdot\mathbf{n}(\overline{\mathbf{x}}) f({{\mathbf{x}}})-%
 (\mathbf{y}-\overline{\mathbf{x}})\cdot\mathbf{n}(\overline{\mathbf{x}}) [ u_{\delta}({{\mathbf{x}}})]_{pp}.
\end{align}
{Furthermore, we replace $[ u_{\delta}({{\mathbf{x}}})]_{pp}$ with its approximation %
$2\int_{-\delta}^{\delta} H_{\delta}(|l|) (u_{\delta}(\mathbf{x}_{l})-u_{\delta}(\mathbf{x})) d\mathbf{x}_{l}-\kappa(\overline{\mathbf{x}})g(\overline{\mathbf{x}})$, where $d\mathbf{x}_{l}$ is the line integral along the contour $\Gamma(\mathbf{x})$, $\kappa(\overline{\mathbf{x}})$ is the curvature of $\partial\Omega$ at $\overline{\mathbf{x}}$, and  $H_{\delta}(|r|)=\dfrac{c}{\delta^3}H\left(\dfrac{|r|}{\delta}\right)$ is the kernel for %
1D nonlocal diffusion model. 
Similar to the requirements for $J$, we assume here $H:[0,\infty)\rightarrow\mathbb{R}$ to be a nonnegative and continuous function with $\int_{\mathbb{R}} H(|z|) |z|^2 dz=1$, $\int_{\mathbb{R}} H(|z|) |z|^k dz=o(1)$ for all $k\geq 3$. %
$H(r)$ is nonincreasing in $r$, strictly positive in $[0,1]$ and vanishes for $|z|>1$. Moreover, we add a further requirement on $H$ that $\int_{\mathbb{R}} H(z) dz:=C_H< \infty$. Here we note that $2\int_{-\delta}^{\delta} H_{\delta}(|l|) (u_{\delta}(\mathbf{x}_{l})-u_{\delta}(\mathbf{x})) d\mathbf{x}_{l}$ is a nonlocal version of the Laplace-Beltrami operator defined on $\Gamma(\mathbf{x})$. 
Substituting the above two approximations into \eqref{eqn:formula1}, we obtain the following model
\begin{align}
 \nonumber&-2\int_{\Omega} J_{\delta}(|\mathbf{x}-\mathbf{y}|)(u_{\delta}(\mathbf{y})-u_{\delta}(\mathbf{x}))d\mathbf{y}-%
 2M_\delta(\mathbf{x}) \int_{-\delta}^{\delta} H_{\delta}(|l|) (u_{\delta}(\mathbf{x}_{l})-u_{\delta}(\mathbf{x})) d\mathbf{x}_{l}\\
 \nonumber=&f({{\mathbf{x}}})-\int_{\mathbb{R}^N\backslash\Omega} J_{\delta}(|\mathbf{x}-\mathbf{y}|)%
 \left[|(\mathbf{y}-\overline{\mathbf{x}})\cdot\mathbf{n}(\overline{\mathbf{x}})|^2-|(\mathbf{x}-\overline{\mathbf{x}})\cdot\mathbf{n}%
 (\overline{\mathbf{x}})|^2\right]d\mathbf{y}f({{\mathbf{x}}})\\
 &+\left(2\int_{\mathbb{R}^N\backslash\Omega} J_{\delta}(|\mathbf{x}-\mathbf{y}|)(\mathbf{y}-\mathbf{x})\cdot\mathbf{n}(\overline{\mathbf{x}})%
 d\mathbf{y}-M_\delta(\mathbf{x})\kappa(\overline{\mathbf{x}})\right)g(\overline{\mathbf{x}}).\label{eqn:formula3}
\end{align}
where
$$M_\delta(\mathbf{x}):=\int_{\mathbb{R}^N\backslash\Omega} J_{\delta}(|\mathbf{x}-\mathbf{y}|)\left[|(\mathbf{y}-\mathbf{x})%
 \cdot\mathbf{p}(\overline{\mathbf{x}})|^2-|(\mathbf{y}-\overline{\mathbf{x}})\cdot\mathbf{n}(\overline{\mathbf{x}})|^2+|(\mathbf{x}-\overline{\mathbf{x}})\cdot\mathbf{n}(\overline{\mathbf{x}})|^2\right]%
 d\mathbf{y}.$$
}
Thus, by defining the nonlocal operator
\begin{align}
 \nonumber L_{N\delta} u:=&-2\int_{\Omega} J_{\delta}(|\mathbf{x}-\mathbf{y}|)(u(\mathbf{y})-u(\mathbf{x}))d\mathbf{y}-2M_\delta(\mathbf{x}){\int_{-\delta}^{\delta} H_{\delta}(|l|) (u(\mathbf{x}_{l})-u(\mathbf{x})) d\mathbf{x}_{l}}
\end{align}
and
\begin{align}
 \nonumber{f}_\delta(\mathbf{x}):=&f({{\mathbf{x}}})-\int_{\mathbb{R}^N\backslash\Omega} J_{\delta}(|\mathbf{x}-\mathbf{y}|)%
 \left[|(\mathbf{y}-\overline{\mathbf{x}})\cdot\mathbf{n}(\overline{\mathbf{x}})|^2-|(\mathbf{x}-\overline{\mathbf{x}})\cdot\mathbf{n}%
 (\overline{\mathbf{x}})|^2\right]d\mathbf{y}f({{\mathbf{x}}})\\
 &+\left(2\int_{\mathbb{R}^N\backslash\Omega} J_{\delta}(|\mathbf{x}-\mathbf{y}|)(\mathbf{y}-\mathbf{x})\cdot\mathbf{n}(\overline{\mathbf{x}})%
 d\mathbf{y}-M_\delta(\mathbf{x})\kappa(\overline{\mathbf{x}})\right)g(\overline{\mathbf{x}}), \label{eqn:2.8}
\end{align}
the proposed algorithm is equivalent to the following nonlocal integral equation
\begin{equation}\label{eqn:nonlocaleqn}
 \left\{\begin{array} {cl}
  L_{N\delta} u_\delta=f_\delta,&\text{ in }\Omega\\
  \int_\Omega u_\delta d\mathbf{x}=0.&\\
 \end{array}\right.
\end{equation}
The corresponding nonlocal weak formulation can then be introduced
\begin{equation}\label{eqn:weakformula}
 B_\delta(u_\delta,v)=({f}_\delta,v)_{L^2(\Omega)},
\end{equation}
where $B_\delta(u,v)$ denotes a nonsymmetric bilinear form $B_\delta(u,v):=(L_{\delta} u,v)$. We note that
\begin{align*}
 &-2\int_{\Omega}\int_{\Omega} J_{\delta}(|\mathbf{x}-\mathbf{y}|)(u(\mathbf{y})-u(\mathbf{x}))d\mathbf{y}v(\mathbf{x})d\mathbf{x}\\
 =&-\int_{\Omega}\int_{\Omega} J_{\delta}(|\mathbf{x}-\mathbf{y}|)(u(\mathbf{y})-u(\mathbf{x}))v(\mathbf{x})d\mathbf{y}d\mathbf{x}
 -\int_{\Omega}\int_{\Omega} J_{\delta}(|\mathbf{y}-\mathbf{x}|)(u(\mathbf{x})-u(\mathbf{y}))v(\mathbf{y})d\mathbf{y}d\mathbf{x}\\
 \nonumber =&\int_{\Omega}\int_{\Omega} J_{\delta}(|\mathbf{x}-\mathbf{y}|)[u(\mathbf{y})-u(\mathbf{x})][v(\mathbf{y})-v(\mathbf{x})]d\mathbf{y}d\mathbf{x}
\end{align*}
and
\begin{align*}
 &-2\int_{\Omega}M_\delta(\mathbf{x}){\int_{-\delta}^{\delta} H_{\delta}(|l|) (u(\mathbf{x}_{l})-u(\mathbf{x})) d\mathbf{x}_{l}}v(\mathbf{x})d\mathbf{x}\\
\nonumber=&-2\int_{\Omega}\int_{\Omega}M_\delta(\mathbf{x}) H_{\delta}(|l|)D(\mathbf{x},\mathbf{y}) (u(\mathbf{y})-u(\mathbf{x})) v(\mathbf{x})d\mathbf{y}d\mathbf{x}\\
\nonumber=&\int_{\Omega}\int_{\Omega}M_\delta(\mathbf{x}) H_{\delta}(|l|)D(\mathbf{x},\mathbf{y}) (u(\mathbf{y})-u(\mathbf{x})) [v(\mathbf{y})-v(\mathbf{x})]d\mathbf{y}d\mathbf{x}\\
&-\int_{\Omega}\int_{\Omega}M_\delta(\mathbf{y}) H_{\delta}(|l|)D(\mathbf{y},\mathbf{x}) (u(\mathbf{x})-u(\mathbf{y})) v(\mathbf{x})d\mathbf{x}d\mathbf{y}\\
&-\int_{\Omega}\int_{\Omega}M_\delta(\mathbf{x}) H_{\delta}(|l|)D(\mathbf{x},\mathbf{y}) (u(\mathbf{y})-u(\mathbf{x})) v(\mathbf{x})d\mathbf{y}d\mathbf{x}\\
\nonumber=&\int_{\Omega}\int_{\Omega}M_\delta(\mathbf{x}) H_{\delta}(|l|)D(\mathbf{x},\mathbf{y}) [u(\mathbf{y})-u(\mathbf{x})] [v(\mathbf{y})-v(\mathbf{x})]|\mathbf{r}'(\mathbf{y})|d\mathbf{y}d\mathbf{x}\\
&+\int_{\Omega}\int_{\Omega} \left[M_\delta(\mathbf{y})\dfrac{|\mathbf{r}'(\mathbf{x})|}{|\mathbf{r}'(\mathbf{y})|}-M_\delta(\mathbf{x})\right]H_{\delta}(|l|)D(\mathbf{x},\mathbf{y}) (u(\mathbf{y})-u(\mathbf{x})) v(\mathbf{x})|\mathbf{r}'(\mathbf{y})|d\mathbf{y}d\mathbf{x}\\
 \nonumber =&\int_{\Omega_\delta} M_\delta(\mathbf{x}) \int_{-\delta}^{\delta} H_{\delta}(|l|) [u(\mathbf{x}_{l})-u(\mathbf{x})]%
 [v(\mathbf{x}_{l})-v(\mathbf{x})] d\mathbf{x}_{l}d\mathbf{x}\\
&+\int_{\Omega_\delta} \int_{-\delta}^{\delta} \left[M_\delta(\mathbf{x}_{l})\dfrac{|\mathbf{r}'(\mathbf{x})|}{|\mathbf{r}'(\mathbf{x}_l)|}-M_\delta(\mathbf{x})\right] H_{\delta}(|l|) [u(\mathbf{x}_{l})-u(\mathbf{x})]%
 d\mathbf{x}_l v(\mathbf{x}) d\mathbf{x},
\end{align*}
where $\mathbf{r}$ is the bijective parametrization of $\Gamma(\mathbf{x})$, $|\mathbf{r}'(\mathbf{x})|$ is the Jacobian of $\mathbf{r}$, and $D(\mathbf{x},\mathbf{y})$ denotes a Dirac-Delta function:
\begin{displaymath}
D(\mathbf{x},\mathbf{y}):=\lim_{\epsilon\rightarrow 0}\epsilon^{-1}\psi(\text{dist}(\mathbf{y},\Gamma(\mathbf{x}))/\epsilon),\;\text{where }\psi \text{ is a mollifier function on }\mathbb{R}.
\end{displaymath}
Therefore
{\begin{align}
 \nonumber B_\delta(u,v)=&\int_{\Omega}\int_{\Omega} J_{\delta}(|\mathbf{x}-\mathbf{y}|)[u(\mathbf{y})-u(\mathbf{x})][v(\mathbf{y})-v(\mathbf{x})]d\mathbf{y}d\mathbf{x}\\
 \nonumber&+\int_{\Omega_\delta} M_\delta(\mathbf{x}) \int_{-\delta}^{\delta} H_{\delta}(|l|) [u(\mathbf{x}_{l})-u(\mathbf{x})]%
 [v(\mathbf{x}_{l})-v(\mathbf{x})] d\mathbf{x}_{l}d\mathbf{x}\\
&+\int_{\Omega_\delta} \int_{-\delta}^{\delta} \left[M_\delta(\mathbf{x}_{l})\dfrac{|\mathbf{r}'(\mathbf{x})|}{|\mathbf{r}'(\mathbf{x}_l)|}-M_\delta(\mathbf{x})\right] H_{\delta}(|l|) [u(\mathbf{x}_{l})-u(\mathbf{x})]%
 d\mathbf{x}_l v(\mathbf{x}) d\mathbf{x}.\label{eqn:B}
\end{align}}
We then consider the nonlocal energy seminorm $||\cdot||_{S_\delta}$ as 
{\begin{align*}
||u||^2_{S_\delta}=&\int_{\Omega}\int_{\Omega} J_{\delta}(|\mathbf{x}-\mathbf{y}|)[u(\mathbf{y})-u(\mathbf{x})]^2d\mathbf{y}d\mathbf{x}
+\int_{\Omega_\delta} M_\delta(\mathbf{x}) \int_{-\delta}^{\delta} H_{\delta}(|l|) [u(\mathbf{x}_{l})-u(\mathbf{x})]^2d\mathbf{x}_{l}d\mathbf{x},
\end{align*} }
with corresponding constrained energy space given by
\begin{displaymath}
 S_{\delta}(\Omega)=\left\{u\in L^2(\Omega): ||u||_{S_\delta}<\infty,\, \int_\Omega u d\mathbf{x}=0\right\}.
\end{displaymath}
Given the nonlocal Poincare inequality which will be addressed in the next section, we will see that $||\cdot||_{S_\delta}$ is actually a full norm. Similar to \cite{mengesha2013analysis}, one can show that the constrained energy space $S_{\delta}(\Omega)$ is a Hilbert space under the given assumptions for the kernels $J$ and $H$. 


\begin{rem}
A similar form of the flux condition \eqref{eqn:formula1} has been proposed in the previous literature, e.g., \cite{cortazar2008approximate,tao2017nonlocal}. By comparing the second term of \eqref{eqn:formula1} with the first case in \cite{cortazar2008approximate}, one can see that the second term of \eqref{eqn:formula1} can be obtained by taking $G=G_1$ in  \cite{cortazar2008approximate} and modifying the correction term $\int_{\mathbb{R}^N\backslash\Omega} G_{\delta}(\mathbf{x},\mathbf{x}-\mathbf{y})%
 g({\mathbf{y}})d\mathbf{y}$ as $\int_{\mathbb{R}^N\backslash\Omega} G_{\delta}(\mathbf{x},\mathbf{x}-\mathbf{y})%
 \dfrac{1}{2}(g({\mathbf{x}})+g({\mathbf{y}}))d\mathbf{y}$. Actually, this modification is sufficient to provide a nonlocal Neumann-type condition with second order accuracy in the 1D case, as shown in \cite{tao2017nonlocal}. However, in higher dimensional cases we need to add the third term of \eqref{eqn:formula1} to achieve second order accuracy.
\end{rem}

\begin{rem}
Note that in the current paper we focus on the $2D$ nonlocal diffusion problem, while the idea can be further extended to the $3D$ cases and to more general nonlocal IDEs, which will be addressed in future work.
\end{rem}


\section{Well-Posedness and Asymptotic Property}\label{sec:l2}

In this section, we first address the well-posedness of the proposed nonlocal Neumann volume-constrained problem by providing a nonlocal Poincar\'{e}-type inequality based on the estimates for boundary curvature $\kappa(\mathbf{x})$ and its derivative $\kappa'(\mathbf{x})$. The coercivity and boundedness of the nonsymmetric bilinear operator $B_\delta(\cdot,\cdot)$ defined in \eqref{eqn:B} follow, which yield the well-posedness of the variational problem. Furthermore, we study the consistency of the nonlocal problem with the classical local model. Specifically, following the framework introduced in \cite{tian2014asymptotically} we prove the uniform embedding property and the precompact property of the proposed norm $S_\delta$, and then show the asymptotic property of the solution of \eqref{eqn:nonlocaleqn} as $\delta\rightarrow 0$, i.e., the solution $u_\delta$ converges to the solution $u_0$ from the limiting local model \eqref{eqn:localeqn}. Here for simplicity we consider the case when $g(\mathbf{x})=0$, and defer discussion of inhomogeneous boundary conditions until Remark \ref{remark:inhomo}. For the limiting local model one can define the corresponding inner product $||u||_{S_0}=||\nabla u||_{L^2(\Omega)}$, the bilinear form %
$B_0(u,v)=(\nabla u,\nabla v)$ and the constrained energy space $S_0=\left\{u\in H^1(\Omega):\int_\Omega u d\mathbf{x}=0\right\}$. %
Throughout this section, we consider the symbol ``$C$'' to indicate a generic constant that is independent of $\delta$, but may have different numerical values in different situations.
Moreover, we introduce the following notation for simplicity:
\begin{align*}
 b_\delta(u,v)& :=\int_{\Omega}\int_{\Omega} J_\delta(|\mathbf{x}-\mathbf{y}|)(u(\mathbf{y})-u(\mathbf{x}))(v(\mathbf{y})-v(\mathbf{x}))
  d\mathbf{y} d\mathbf{x},\\
  h_\delta(u,v)& :={\int_{\Omega_\delta} M_\delta(\mathbf{x}) \int_{-\delta}^{\delta} H_{\delta}(|l|) [u(\mathbf{x}_{l})-u(\mathbf{x})]%
 [v(\mathbf{x}_{l})-v(\mathbf{x})] d\mathbf{x}_{l}d\mathbf{x}},\\
I_\delta(\mathbf{x},\mathbf{y})&:= |(\mathbf{y}-\mathbf{x})\cdot \mathbf{p}(\overline{\mathbf{x}})|^2 -
|(\mathbf{y}-\overline{\mathbf{x}})\cdot \mathbf{n}(\overline{\mathbf{x}})|^2   
+|(\mathbf{x}-\overline{\mathbf{x}})\cdot \mathbf{n}(\overline{\mathbf{x}})|^2.
\end{align*}
{We first have the following estimates of the function $M_\delta(\mathbf{x})$ for each $\mathbf{x}\in \Omega_\delta$: 
 \begin{lem}\label{Mdelta}
For $l\in[-\delta,\delta]$, and assuming that there exist constants $d,D>0$ such that $|\kappa'(\mathbf{z})|\leq D$, $|\kappa(\mathbf{z})|\leq D$ and %
$\sup_{|\xi|\leq d}\left| \dfrac{\kappa'(\mathbf{z}_\xi)}{\kappa(\mathbf{z}) }\right|\leq D$ for almost every $\mathbf{z}\in \partial\Omega$, there exists a $0<\overline{\delta}\leq d$ such that for $\delta\leq \overline{\delta}$ for almost every $\mathbf{x}\in \partial \Omega$ we have $0\leq M_\delta(\mathbf{x})\leq C$ and
\begin{equation}\label{eqn:estM-}
 \left|M_\delta(\mathbf{x})-M_\delta(\mathbf{x}_l)\dfrac{|\mathbf{r}'(\mathbf{x})|}{|\mathbf{r}'(\mathbf{x}_l)|}\right|\leq C_M \delta^2,
\end{equation}
\begin{equation}\label{eqn:estM-1} 
\bigg|\frac{M_\delta(\mathbf{x})|\mathbf{r}'(\mathbf{x}_l)|-M_\delta(\mathbf{x}_l)|\mathbf{r}'(\mathbf{x})|}{M_\delta(\mathbf{x})|\mathbf{r}'(\mathbf{x}_l)|}\bigg| 
 \le C_N \delta,
\end{equation}
where $C_M$, $C_N$ are constants independent of $\delta$.
\end{lem}

\begin{proof}
\begin{figure}[!htb]
\centering
\subfigure{\includegraphics[scale=.4]{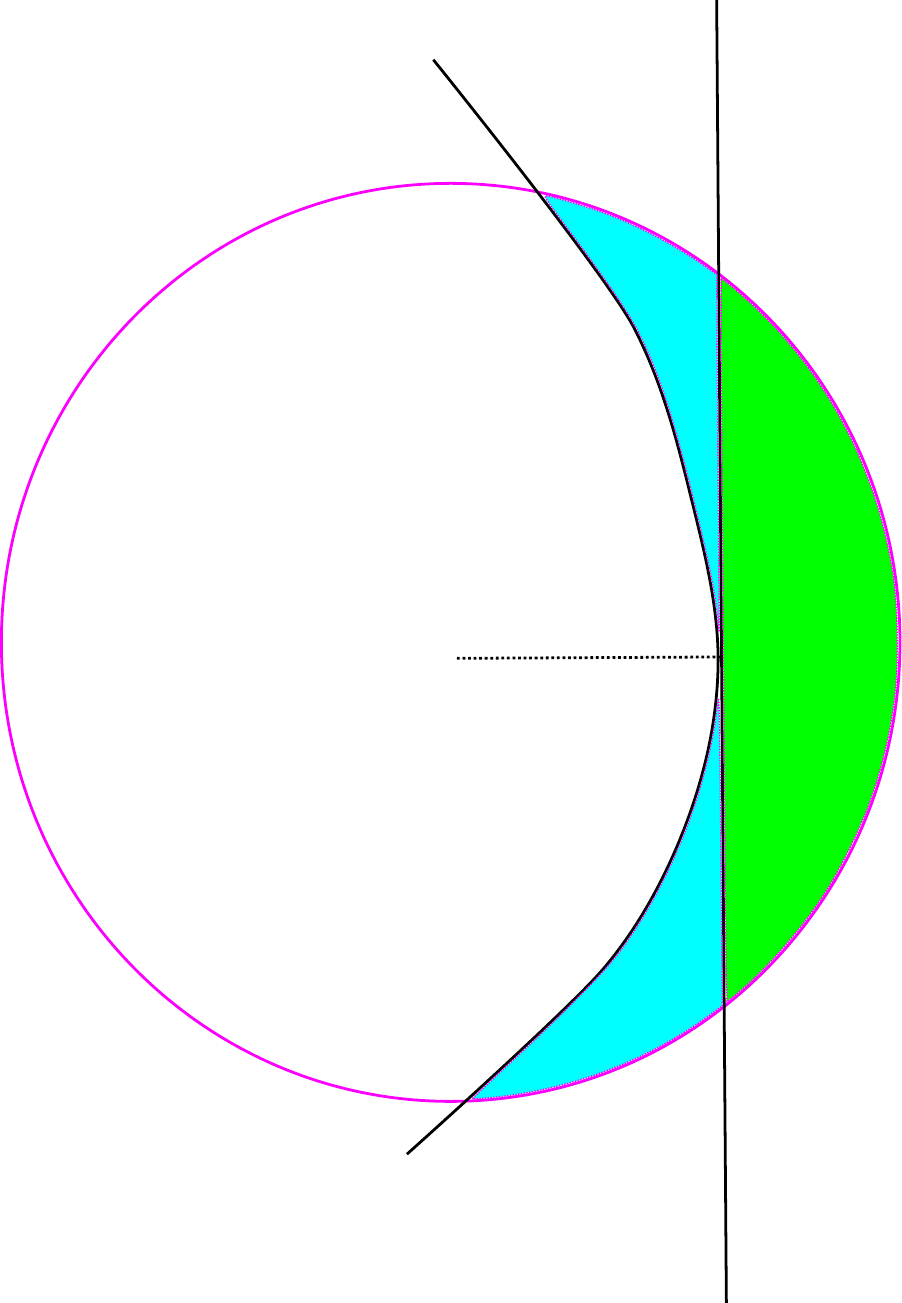}
 \put(-55,70){$\mathbf{x}$}\put(-20,70){$\overline{\mathbf{x}}$}
 \put(-20,10){$\tau(\overline{\mathbf{x}})$}
 \put(-15,90){$D_\delta$}
 \put(-33,115){$A_\delta$}\put(-35,35){$A_\delta$}
 \put(-50,135){$\partial\Omega$}
 \put(-100,90){$B(\mathbf{x},\delta)$}}\qquad\qquad
\subfigure{\includegraphics[scale=0.5]{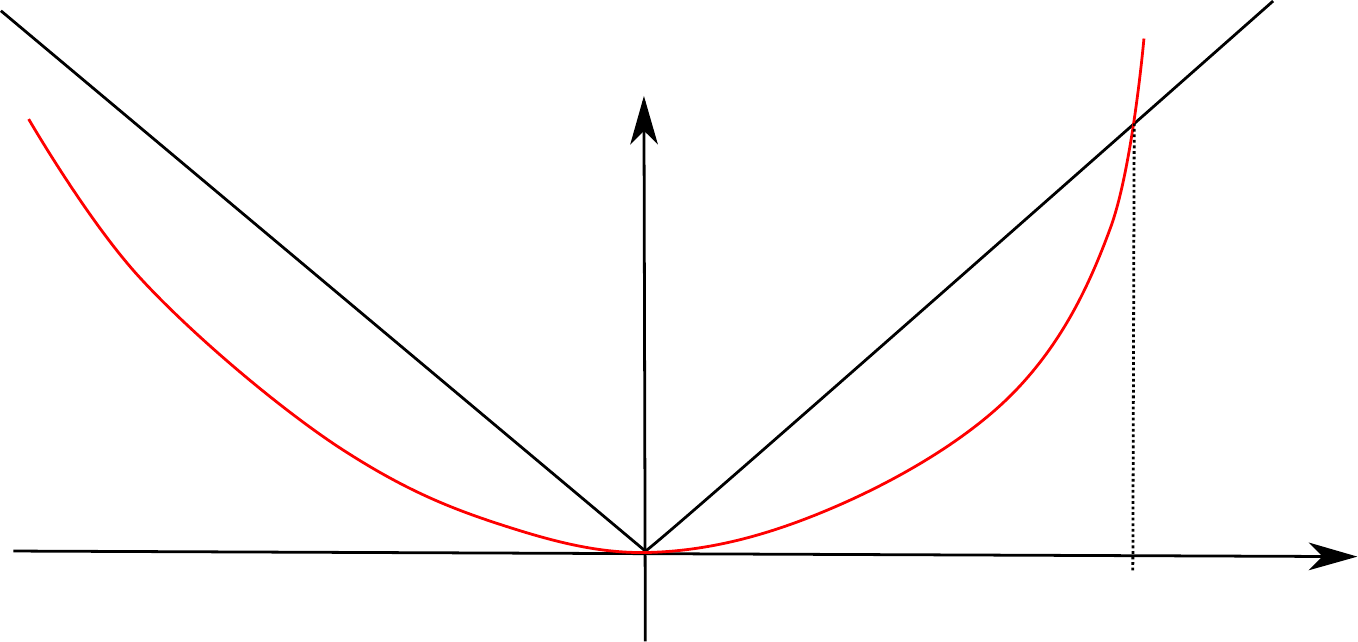}
 \put(-100,2){$0$}
 \put(-120,85){$y\parallel\mathbf{n}(\overline{\mathbf{x}}) $}
 \put(-10,2){$x\parallel\mathbf{p}(\overline{\mathbf{x}}) $}
 \put(-40,2){$x_0$}
 \put(-180,30){$\substack{y=f(x)\\ (=\partial\Omega)}$}
 \put(-180,90){$y=|x|$}} 
 \caption{Notation for the geometric estimates in Lemma \ref{Mdelta}. 
 Left: illustration of regions $D_\delta$ and $A_\delta$. Green represents $D_\delta$, the region in $B(\mathbf{x},\delta)$ which lies on the other side of the tangential line at $\overline{\mathbf{x}}$ with respect to $\Omega$. Cyan represents $A_\delta$, the region in $B(\mathbf{x},\delta)$ which lies between $\partial\Omega$ and the tangential line. Right: Representation of the Cartesian coordinate system locally near $\overline{\mathbf{x}}$. Here the region $A_\delta$ lies below the red curve $y=f(x)$.}
 \label{graph}
\end{figure}
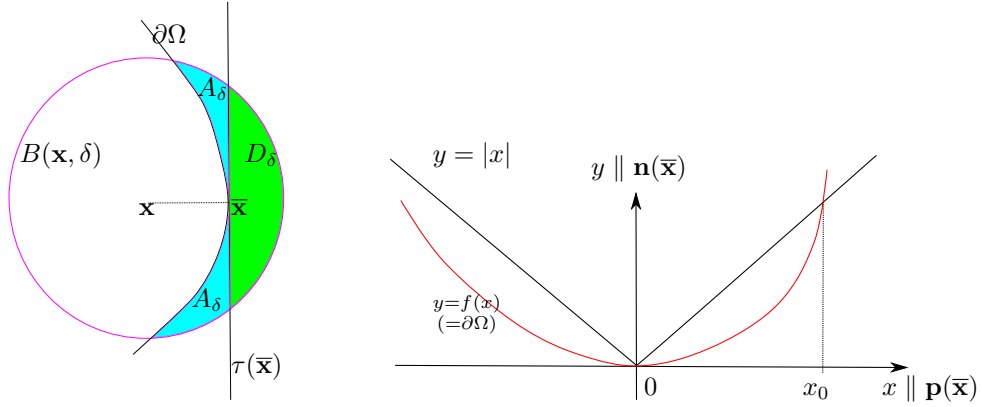

We show now that $0\leq M_\delta(\mathbf{x})\leq C$. Note that
\[M_\delta(\mathbf{x})=\int_{(\mathbb{R}^N\backslash \Omega)_\delta}  J_\delta( |\mathbf{x}-\mathbf{y}| ) I_\delta d \mathbf{y} = 
\int_{D_\delta}  J_\delta( |\mathbf{x}-\mathbf{y}| ) I_\delta d \mathbf{y}
+\int_{A_\delta}  J_\delta( |\mathbf{x}-\mathbf{y}| ) I_\delta d \mathbf{y}.\]
With $\tau(\overline{\mathbf{x}})$ representing the tangent line to $\partial\Omega$
 at  $\overline{\mathbf{x}}$, here $D_\delta$ is the region of $B(\mathbf{x},\delta)$ on the side of $\tau(\overline{\mathbf{x}})$ {\em not} containing $\Omega$ (as shown in the green region in the left plot of Figure \ref{graph}), and $A_\delta:=B(\mathbf{x},\delta) \backslash (D_\delta\cup \Omega)$ (as shown in the cyan region in the left plot of Figure \ref{graph}). %
We consider first the $D_\delta$ part. One can rewrite %
$\mathbf{y}\in D_\delta
$ as $\mathbf{x}+(r\cos(\theta),r\sin(\theta))$ with $s_x<r<\delta$ and %
$-\pi/2<-\arccos(s_x/r)\leq\theta\leq \arccos(s_x/r)\leq \pi/2$, %
which yields
\begin{align}
 \nonumber\int_{D_\delta}  J_\delta( |\mathbf{x}-\mathbf{y}| ) I_\delta d \mathbf{y}&=%
 \int_{s_x}^\delta J_\delta(r)\int_{-\arccos(s_x/r)}^{\arccos(s_x/r)} I_\delta (\mathbf{x},\mathbf{y})rd\theta dr\\
 &=2\int_{s_x}^\delta J_\delta(r) r^2 \sqrt{1-(s_x/r)^2} s_xdr.\label{eqn:intT}
\end{align}
From \eqref{eqn:intT} we can see that 
$\int_{D_\delta} J_\delta(|\mathbf{x}-\mathbf{y}|)I_\delta d\mathbf{y}\ge0$ and
\begin{align}
M_\delta(\mathbf{x}) &\ge \int_{A_\delta} J_\delta(|\mathbf{x}-\mathbf{y}|)[|(\mathbf{y}-\overline{\mathbf{x}})\cdot \mathbf{p}(\overline{\mathbf{x}})|^2 -
|(\mathbf{y}-\overline{\mathbf{x}})\cdot \mathbf{n}(\overline{\mathbf{x}})|^2  ] d\mathbf{y} . \label{convex}
 \end{align}
Therefore it suffices to show now that 
\begin{equation}
 |(\mathbf{y}-\overline{\mathbf{x}})\cdot \mathbf{p}(\overline{\mathbf{x}})|\ge
|(\mathbf{y}-\overline{\mathbf{x}})\cdot \mathbf{n}(\overline{\mathbf{x}})| \qquad \fal \mathbf{y}\in A_\delta.
 \label{geometry}
\end{equation}
We adopt a Cartesian coordinate system as shown in the right plot of Figure \ref{graph}, assuming that $\overline{\mathbf{x}}$ coincides with the origin, $\mathbf{p}(\overline{\mathbf{x}})$ is oriented along the positive direction of the $x$-axis while $\mathbf{n}(\overline{\mathbf{x}})$ coincides with the negative direction of the $y$-axis. We then have $\overline{\mathbf{x}}=(0,0)$, $\tau(\overline{\mathbf{x}})=\{y=0\}$, $\Omega\subseteq\{y\ge 0\}$, and let $y=f(x)$ be the curve describing $\partial\Omega$. We note that any point $\mathbf{y}$ lying below $y=|x|$ satisfies \eqref{geometry}. Assuming that there exists a point $\mathbf{z}\in \partial\Omega$ lying above $y=|x|$, there exists $x_0\neq 0$ such that $f(x_0)=|x_0|$ and $(x_0,f(x_0))\in \partial\Omega$. For simplicity we consider the case where $x_0>0$ since the case where $x_0<0$ is analogous. Since $f'(0)=0$, by continuity there exists at least one point $x_1\in (0,x_0)$ such that $f'(x_1)\geq 1$. Let $x_2:=\inf\{t>0: f'(t) \ge 1 \}\le x_1\le x_0$, then by the regularity of $f$ we have $x_2>0$. Thus $f'(x_2)-f'(0) =1 =\int_0^{x_2}f''(s)ds$. Moreover, the unsigned curvature of the graph of $f$ can be given by $|f''(x)|(1+f'(x)^2)^{-3/2}$. Due to the finiteness of the curvature of $\partial \Omega$, and the fact that $f'(x)^2\le 1$ for all $x\in [0,x_2]$, we obtain $D\geq |f''(x)|(1+f'(x)^2)^{-3/2}$ and therefore
\begin{align*}
|f''(x)| \le D(1+f'(x)^2)^{3/2}
\le 2\sqrt{2}D, \qquad \fal x\in [0,x_2].
\end{align*}
Hence
\begin{align*}
1 &=\int_0^{x_2}f''(s)ds \le 
\int_0^{x_2}|f''(s)|ds \le 2\sqrt{2}x_2D\Rightarrow x_2\ge \frac{1}{2\sqrt{2}D}.
\end{align*}
But since $x_2\le x_0$, this means that the first intersection point between
$y=f(x)$ and $y=|x|$ (which we denote as $\mathbf{w}=(x_0,|x_0|)$) has distance at least $x_2$ from $\overline{\mathbf{x}}=(0,0)$. 
Thus, for sufficiently small $\delta\leq\frac{1}{4\sqrt{2}D}=:\delta_1$, we get
\begin{align*}
|\mathbf{w}-\overline{\mathbf{x}}| &\ge x_2\ge \frac{1}{2\sqrt{2}D} \geq 2\delta=\sup_{p,q\in B(\mathbf{x},\delta)}|p-q|.
\end{align*}
Therefore, $\mathbf{w}\notin B(\mathbf{x},\delta)$, and the entire region $A_\delta$ lies below $y=|x|$. Consequently, any $\mathbf{y}\in A_\delta$ satisfies
$ |(\mathbf{y}-\overline{\mathbf{x}})\cdot \mathbf{p}(\overline{\mathbf{x}})|\ge
|(\mathbf{y}-\overline{\mathbf{x}})\cdot \mathbf{n}(\overline{\mathbf{x}})|$ and in turn $M_\delta\ge 0$. On the other hand, with the $C^3$ regularity of $\Omega$ and by Taylor expansion $B(\mathbf{x},\delta)\cap \partial\Omega$ is the graph of a function of the form $y=f(x)=\frac{\kappa(\overline{\mathbf{x}})}{2}x^2+O(x^3)$. Therefore, the area $|A_\delta|\leq C |\kappa(\overline{\mathbf{x}})|(\delta^2-s_x^2)^{3/2}\leq CD\delta^3$. Hence
\begin{align*}
M_\delta(\mathbf{x}) \leq & 
\left|\int_{D_\delta}  J_\delta( |\mathbf{x}-\mathbf{y}| ) I_\delta d \mathbf{y}\right|
+\left|\int_{A_\delta}  J_\delta( |\mathbf{x}-\mathbf{y}| ) I_\delta d \mathbf{y}\right|\\
\leq& C \sup_{r}J(|r|)(\delta^{-4}s_x(\delta^2-s_x^2)^{3/2}+ D \delta)\leq C.
\end{align*}

To show \eqref{eqn:estM-}, denoting by $D_{\delta,\mathbf{x}_{l}},A_{\delta,\mathbf{x}_{l}}$ the analogous sets of $D_\delta,A_\delta$ at $\mathbf{x}_{l}$ instead of $\mathbf{x}$, we then have
\begin{align*}
M_\delta(\mathbf{x})-M_\delta(\mathbf{x}_l)\dfrac{|\mathbf{r}'(\mathbf{x})|}{|\mathbf{r}'(\mathbf{x}_l)|}
=&\int_{D_{\delta}}  J_\delta(|\mathbf{x}-\mathbf{y}|) I_\delta(\mathbf{x},\mathbf{y}) d \mathbf{y}-\dfrac{|\mathbf{r}'(\mathbf{x})|}{|\mathbf{r}'(\mathbf{x}_l)|}\int_{D_{\delta,\mathbf{x}_{l}}}  J_\delta(|\mathbf{x}_{l}-\mathbf{y}|) I_{\delta}(\mathbf{x}_{l},\mathbf{y})d \mathbf{y}\\
&+\int_{A_\delta}  J_\delta( |\mathbf{x}-\mathbf{y}| ) I_\delta(\mathbf{x},\mathbf{y}) d \mathbf{y}-\dfrac{|\mathbf{r}'(\mathbf{x})|}{|\mathbf{r}'(\mathbf{x}_l)|}\int_{A_{\delta,\mathbf{x}_{l}}}  J_\delta(|\mathbf{x}_{l}-\mathbf{y}|) I_{\delta}(\mathbf{x}_{l},\mathbf{y})d \mathbf{y}.
\end{align*}
{With the definition of $\mathbf{x}_{l}$ and the regularity assumptions on $\Omega$, it holds $s_{\mathbf{x}_{l}}:= \text{dist}(\mathbf{x}_{l},\partial\Omega)
=\text{dist}(\mathbf{x},\partial\Omega)$. We obtain
 \[\int_{D_{\delta,\mathbf{x}_{l}}}  J_\delta(|\mathbf{x}_{l}-\mathbf{y}|) I_{\delta}(\mathbf{x}_{l},\mathbf{y})d \mathbf{y}=\int_{D_{\delta}}  J_\delta(|\mathbf{x}-\mathbf{y}|) I_\delta(\mathbf{x},\mathbf{y}) d \mathbf{y}. \]
 Moreover, with the coordinate system as shown in the right plot of Figure \ref{graph}, we have $|\mathbf{r}'(\mathbf{x})|=1$ and $|\mathbf{r}'(\mathbf{x}_l)|=\sqrt{1+(f'(x))^2}$. Since for any point $\mathbf{x}_l=(x_l,f(x_l))$ in $B(\mathbf{x},\delta)$, $|f'(x_l)|=|x_lf''(\xi)|\leq C\delta$ for some $\xi\in[0,x_l]$, therefore
\begin{equation}\label{eqn:rratio}
\left|\dfrac{|\mathbf{r}'(\mathbf{x})|}{|\mathbf{r}'(\mathbf{x}_l)|}-1\right|=\left|\dfrac{1}{\sqrt{1+(f'(x_l))^2}}-1\right|=\dfrac{1}{2}(f'(x_l))^2+O(\delta^4)\leq C \delta^2
\end{equation}
 and hence together with \eqref{eqn:intT} we obtain
 \[\int_{D_{\delta}}  J_\delta(|\mathbf{x}-\mathbf{y}|) I_\delta(\mathbf{x},\mathbf{y}) d \mathbf{y}-\dfrac{|\mathbf{r}'(\mathbf{x})|}{|\mathbf{r}'(\mathbf{x}_l)|}\int_{D_{\delta,\mathbf{x}_{l}}}  J_\delta(|\mathbf{x}_{l}-\mathbf{y}|) I_{\delta}(\mathbf{x}_{l},\mathbf{y})d \mathbf{y}\leq C\delta^2.\]
 To estimate  $\int_{A_{\delta}}  J_\delta( |\mathbf{x}-\mathbf{y}| ) I_{\delta}(\mathbf{x},\mathbf{y}) d \mathbf{y}-\int_{A_{\delta,\mathbf{x}_{l}}}  J_\delta( |\mathbf{x}_{l}-\mathbf{y}| ) I_{\delta}(\mathbf{x}_{l},\mathbf{y}) d \mathbf{y}$, let $R$ be the rototranslation mapping such that
 \[R(\overline{\mathbf{x}}_{l})=\overline{\mathbf{x}},\quad R(\mathbf{p}(\overline{\mathbf{x}}_{l}))
 =\mathbf{p}(\overline{\mathbf{x}}),\quad R(\mathbf{n}(\overline{\mathbf{x}}_{l}))=\mathbf{n}(\overline{\mathbf{x}}).\]
With such construction we note that the curves $R(B({\mathbf{x}}_{l}, \delta) \cap \partial \Omega )$ and $B({\mathbf{x}}, \delta) \cap \partial \Omega$ share the same tangential lines at $\overline{\mathbf{x}}$. Meanwhile, $B({\mathbf{x}}, \delta) \cap \partial \Omega$ and $R(B({\mathbf{x}}_{l}, \delta) \cap \partial \Omega)$ have different curvatures $\kappa(\overline{\mathbf{x}})$ and $\kappa(\overline{\mathbf{x}}_{l})$, respectively. When $\delta\ll 1/D$, we have the arc lengths of $R(B({\mathbf{x}}_{l}, \delta) \cap \partial \Omega )$ and $B({\mathbf{x}}, \delta) \cap \partial \Omega$ satisfying $|R(B({\mathbf{x}}_{l}, \delta) \cap \partial \Omega )|\leq 2\sqrt{\delta^2-s_x^2}+C \kappa(\overline{\mathbf{x}}_{l})(\delta^2-s_x^2)$ and $|B({\mathbf{x}}, \delta) \cap \partial \Omega|\leq 2\sqrt{\delta^2-s_x^2}+C \kappa(\overline{\mathbf{x}})(\delta^2-s_x^2)$. Moreover, the spread $d_{\mathcal{H}}(B({\mathbf{x}},\delta)\cap \partial\Omega,  R(B({\mathbf{x}}_{l}, \delta) \cap \partial \Omega) )$ is bounded by 
 \begin{align*}
  &d_{\mathcal{H}}(B({\mathbf{x}},\delta)\cap \partial\Omega,  R(B({\mathbf{x}}_{l}, \delta) \cap \partial \Omega) )\notag\\
 =&\max\Big\{\sup_{z\in B({\mathbf{x}},\delta)\cap\partial \Omega}\text{dist}(z, 
 R(B({\mathbf{x}}_{l}, \delta) \cap \partial \Omega)),
 \sup_{z\in R(B({\mathbf{x}}_{l}, \delta) \cap \partial \Omega)}\text{dist}(z,B({\mathbf{x}},\delta)\cap \partial \Omega)\Big\}\\
 \le& C|\kappa(\overline{\mathbf{x}}_{l})-\kappa(\overline{\mathbf{x}})|(\delta^2-s_x^2).
 \end{align*}
Therefore, noting that the quantities $J_\delta$, $I_\delta$ and $d \mathbf{y}$ are invariant under $R$, and $|I_{\delta}(\mathbf{x}_{l},\mathbf{y})|$, $|I_{\delta}(\mathbf{x},\mathbf{y})|\le 3\delta^2$, one has
 \begin{align}
 \nonumber &\bigg|\int_{A_{\delta,\mathbf{x}_{l}}}  J_\delta( |\mathbf{x}_{l}-\mathbf{y}| ) I_{\delta}(\mathbf{x}_{l},\mathbf{y}) d \mathbf{y}
  -\int_{A_{\delta}}  J_\delta( |\mathbf{x}-\mathbf{y}| ) I_\delta(\mathbf{x},\mathbf{y}) d \mathbf{y}\bigg|\\ 
 \nonumber =& \bigg|\int_{R(A_{\delta,\mathbf{x}_{l}})}  J_\delta( |\mathbf{x}-\mathbf{y}| ) I_{\delta}(\mathbf{x},\mathbf{y}) %
  d \mathbf{y} -\int_{A_{\delta}}  J_\delta( |\mathbf{x}-\mathbf{y}| ) I_\delta(\mathbf{x},\mathbf{y}) d \mathbf{y}\bigg| \\
  \nonumber\le& C\int_{A_{\delta}\Delta R(A_{\delta,\mathbf{x}_{l}})}  J_\delta( |\mathbf{x}-\mathbf{y}| ) \delta^2  d \mathbf{y}
  \le C|\kappa(\overline{\mathbf{x}}_{l})-\kappa(\overline{\mathbf{x}})|\delta^{-2} (\delta^2-s_x^2)^{3/2} \\
  \le& C \sup_{|\xi|\leq |l|}|\kappa'(\overline{\mathbf{x}}_\xi)|\delta^{-1} (\delta^2-s_x^2)^{3/2}\leq C \delta^{-1} (\delta^2-s_x^2)^{3/2},\label{M up}
 \end{align}
where the constant $C$ depends on $\sup_{r} J(|r|)$ and is independent of $\delta$. Moreover, with \eqref{eqn:rratio} and $\left|\int_{A_\delta}  J_\delta( |\mathbf{x}-\mathbf{y}| ) I_\delta d \mathbf{y}\right|\leq C \sup_{r}J(|r|\delta$ we have
\[\left|\left(1-\dfrac{|\mathbf{r}'(\mathbf{x})|}{|\mathbf{r}'(\mathbf{x}_l)|}\right)\int_{A_{\delta,\mathbf{x}_{l}}}  J_\delta(|\mathbf{x}_{l}-\mathbf{y}|) I_{\delta}(\mathbf{x}_{l},\mathbf{y})d \mathbf{y}\right|\leq C \delta^3.\]
Thus, we obtain the bound in \eqref{eqn:estM-}.

We now work on the proof of \eqref{eqn:estM-1} by combining \eqref{M up} and establishing a lower bound for $M_\delta$. 
We firstly prove that 
\begin{equation}\label{eqn:estM-1step1}
\bigg|\frac{M_\delta(\mathbf{x})-M_\delta(\mathbf{x}_l)}{M_\delta(\mathbf{x})}\bigg|\leq C\delta.
\end{equation}
With the previous calculation, we have 
\begin{align*} 
&\int_{D_\delta} J_\delta( \mathbf{x}-\mathbf{y} ) I_\delta d \mathbf{y} =2\int_{s_x}^\delta J_\delta(r) r^2 \sqrt{1-(s_x/r)^2} s_xdr
\geq C \delta^{-4}s_x(\delta^2-s_x^2)^{3/2}\\
\geq &\dfrac{C}{D}|\kappa(\overline{\mathbf{x}})| \delta^{-4}s_x(\delta^2-s_x^2)^{3/2}=C|\kappa(\overline{\mathbf{x}})| \delta^{-4}s_x(\delta^2-s_x^2)^{3/2}
\end{align*}
and $\int_{A_\delta} J_\delta( \mathbf{x}-\mathbf{y} ) I_\delta d \mathbf{y}\geq 0$. When $s^2_x\geq \delta^2/2$ 
one has
\begin{align*}
\bigg|\frac{M_\delta(\mathbf{x})-M_\delta(\mathbf{x}_l)}{M_\delta(\mathbf{x})}\bigg|\leq &C\dfrac{|\kappa(\overline{\mathbf{x}}_{l})-\kappa(\overline{\mathbf{x}})|\delta^{-2} (\delta^2-s_x^2)^{3/2}}{|\kappa(\overline{\mathbf{x}})| \delta^{-4}s_x(\delta^2-s_x^2)^{3/2}}
= C\dfrac{\delta^3 }{ s_x}\sup_{|\xi|\leq d}
\left|\dfrac{\kappa'(\overline{\mathbf{x}}_\xi)}{\kappa(\overline{\mathbf{x}})}\right|\leq C\delta\sup_{|\xi|\leq d}
\left|\dfrac{\kappa'(\overline{\mathbf{x}}_\xi)}{\kappa(\overline{\mathbf{x}})}\right|
\end{align*}
and therefore \eqref{eqn:estM-1step1} holds true. For $s^2_x< \delta^2/2$, we just need to bound
$\int_{A_\delta} J_\delta( \mathbf{x}-\mathbf{y} ) I_\delta d \mathbf{y}$
 from below. For notational simplicity, we assume here the Cartesian coordinate system shown in the right plot of Figure \ref{graph}.  
The following properties hold:
 \begin{align}
  (\mathbf{y}-\overline{\mathbf{x}})\cdot \mathbf{p}(\overline{\mathbf{x}}) = x \text{ coordinate of } \mathbf{y},\\
 \quad (\mathbf{y}-\overline{\mathbf{x}})\cdot \mathbf{n}(\overline{\mathbf{x}}) = - (y \text{ coordinate of } \mathbf{y}).
 \label{xy}
 \end{align}
We first assume that $|\kappa(\overline{\mathbf{x}})|>0$. By Taylor approximation,
 $B(\mathbf{x}, \delta)\cap \partial \Omega$ is the graph of a function of the form
$y=f(x)=\frac{\kappa(\overline{\mathbf{x}})}{2}x^2+O(x^3)$. Integrating it yields that the area $|A_\delta| = C|\kappa(\overline{\mathbf{x}})|(\delta^2-s_x^2)^{3/2}= C|\kappa(\overline{\mathbf{x}})|\delta^3$. Let $h\in (0,1)$ be a point where the area of $A_\delta\cap \{x\ge h\delta\}$ 
 is  $C|\kappa(\overline{\mathbf{x}})|\delta^3/2$. With the convexity assumption of $\partial\Omega$, one has $h>1/2$. When $\delta\leq \delta_1<\frac{1}{2D}$, the slope of $f$ (i.e. the slope of the tangent derivative of $B(\mathbf{x}, \delta)\cap \partial \Omega$) can reach at most $\delta D<1/2$. Thus the graph of $f$ lies below the line $y=x/2$ and \eqref{xy} gives
\begin{align*}
I_\delta\ge &|(\mathbf{y}-\mathbf{x})\cdot \mathbf{p}(\overline{\mathbf{x}})|^2 -|(\mathbf{y}-\overline{\mathbf{x}})\cdot \mathbf{n}(\overline{\mathbf{x}})|^2\geq|(\mathbf{y}-\overline{\mathbf{x}})\cdot \mathbf{p}(\overline{\mathbf{x}})|^2-
 |(\mathbf{y}-\overline{\mathbf{x}})\cdot \mathbf{n}(\overline{\mathbf{x}})|^2\\
 \geq &\frac34 | (\mathbf{y}-\overline{\mathbf{x}})\cdot \mathbf{p}(\overline{\mathbf{x}})|^2
\ge \frac34 h^2 \delta^2\ge \dfrac{3}{16} \delta^2
\end{align*}
 for all $\mathbf{y} \in A_\delta\cap \{x\ge h\delta\}$. 
 {Recalling that $J(r)$ is strictly positive for $0\leq r\leq 1$ and therefore $\min_{r\leq1}J(|r|)>0$}, we infer
 \begin{align}
\nonumber M_\delta\ge& \int_{A_\delta\cap \{x\ge h\}} \min_{r\leq1}J(|r|)\delta^{-4}(|(\mathbf{y}-\overline{\mathbf{x}})\cdot \mathbf{p}(\overline{\mathbf{x}})|^2-
 |(\mathbf{y}-\overline{\mathbf{x}})\cdot \mathbf{n}(\overline{\mathbf{x}})|^2) d \mathbf{y}\\
 \ge&  \frac{3C}{32}|\kappa(\overline{\mathbf{x}})|  \delta\ge C|\kappa(\overline{\mathbf{x}})|  \delta. 
\label{M low}
\end{align}
Combining with \eqref{M up}, we thus obtain the bound \eqref{eqn:estM-1step1}.

For $|\kappa(\overline{\mathbf{x}})|=0$, with domain $C^3$ regularity assumption and $\sup_{|\xi|\leq d}
\left|\dfrac{\kappa'(\overline{\mathbf{x}}_\xi)}{\kappa(\overline{\mathbf{x}})}\right|\leq D$ a.e., we have $\kappa(\overline{\mathbf{x}}_\xi)\equiv 0$ for $|\xi|\leq d$ almost everywhere and therefore $M_\delta(\overline{\mathbf{x}}_l)=M_\delta(\overline{\mathbf{x}})$ for $|l|\leq \delta$ and $\delta\leq d/2$. \eqref{eqn:estM-1step1} can then be trivially proved.  

We can now prove \eqref{eqn:estM-1}:
\begin{align*}
\bigg|\frac{M_\delta(\mathbf{x})|\mathbf{r}'(\mathbf{x}_l)|-M_\delta(\mathbf{x}_l)|\mathbf{r}'(\mathbf{x})|}{M_\delta(\mathbf{x})|\mathbf{r}'(\mathbf{x}_l)|}\bigg|\leq&\dfrac{|\mathbf{r}'(\mathbf{x})|}{|\mathbf{r}'(\mathbf{x}_l)|}\bigg|\frac{M_\delta(\mathbf{x})-M_\delta(\mathbf{x}_l)}{M_\delta(\mathbf{x})}\bigg|+\left|\dfrac{|\mathbf{r}'(\mathbf{x})|}{|\mathbf{r}'(\mathbf{x}_l)|}-1\right|
\leq C \delta.
\end{align*}

}\end{proof}}

\begin{rem}
Note that in the previous proof we assumed $J(r)$ is strictly positive in $[0,1]$ such that $J(|r|)\geq C_1>0$. However, the proof can be extended for a more general positive $J$ whose support is the entire ball $B(0,1)$. It suffices to note:
\begin{itemize}
\item It easily follows from the previous proof that the set 
\[\~A_\delta:=\{\mathbf{z}\in A_\delta: |\mathbf{z}-\mathbf{x}| \in [\delta/3,\delta/2]\}\]
 has area $C|A_\delta|$ for some constant $C\in (0,1)$, and on $\~A_\delta$ it holds
$ I_\delta \ge C_1\delta^2$, again for some constant $C_1\in (0,1)$.
\item Since $J(r)$ is nonincreasing on $r$ and its support is the entire ball $B(0,1)$ there exists another constant $C_2>0$ such that $J(r)\ge C_2$ for $r\in[1/3,1/2]$.
\end{itemize}
Combining the above two facts, we obtain
\begin{align*}
M_\delta &\ge\int_{A_\delta} J_\delta( |\mathbf{x}-\mathbf{y}| ) I_\delta d \mathbf{y}
 \ge \int_{\~A_\delta} \min_{\delta/3\le r\leq \delta/2}J_\delta(r)I_\delta d \mathbf{y}
 \ge  C|\kappa(\overline{\mathbf{x}})|  \delta.
\end{align*}

\end{rem}

\begin{rem}
When $u\in C^{\infty}(\Omega)$, the above bounds of $M_\delta(\mathbf{x})$ yield
\begin{align*}
0\leq h_\delta(u,u)=&{\int_{\Omega_\delta} M_\delta(\mathbf{x}) \int_{-\delta}^{\delta} H_{\delta}(|l|) [u(\mathbf{x}_{l})-u(\mathbf{x})]^2 d\mathbf{x}_{l}d\mathbf{x}}\\
\leq&{\int_{\Omega_\delta} M_\delta(\mathbf{x}) \int_{-\delta}^{\delta} H_{\delta}(|l|) \left[\left|\dfrac{\partial u(\mathbf{x})}{\partial \mathbf{p}}\right|^2|l|^2+C|l|^3\right] d\mathbf{x}_{l}d\mathbf{x}}\\
\leq&\int_{\Omega_\delta} M_\delta(\mathbf{x})\left|\dfrac{\partial u(\mathbf{x})}{\partial \mathbf{p}}\right|^2d\mathbf{x}+C|\Omega_\delta|\delta
\leq C|\Omega_\delta|\left(\left|\dfrac{\partial u(\mathbf{x})}{\partial \mathbf{p}}\right|^2+\delta\right).
\end{align*}
Combining with the results in \cite{tian2014asymptotically}, we have $\lim_{\delta\rightarrow 0}||u||_{S_\delta}=||u||_{S_0}$.
\end{rem}

We will now show a nonlocal Poincare-type inequality: 
\begin{lem}\label{lemma:poincare}
{ 
There exists a $0<\tilde{\delta}\leq 1$ such that
\begin{equation}
 ||u||_{L^2(\Omega)}^2\leq C B_\delta(u,u)
\end{equation}
for all $u\in S_{\delta}$ and $\delta\leq \tilde{\delta}$. Note that here $\tilde{\delta}$ depends on both $u$ and $\Omega$.
}
\end{lem}

\begin{proof} With \cite[Proposition~2]{mengesha2014nonlocal} we have the bound for the first term in \eqref{eqn:B}: there exist $\delta_0$ such that for all $\delta<\delta_0$,
\begin{equation}\label{eqn:C*}
 ||u||_{L^2(\Omega)}^2\leq C^* b_\delta(u,u),
\end{equation}
and here we assume $C^*> 0$ without loss of generality. %
To estimate the remaining two terms, we first work on the case where $\delta\Omega$ is a straight line. For $\mathbf{x}\in\Omega_\delta$ we have {$M_\delta (\mathbf{x}_{l})\dfrac{|\mathbf{r}'(\mathbf{x})|}{|\mathbf{r}'(\mathbf{x}_l)|}=M_\delta (\mathbf{x})$}, and therefore the last term of $B_\delta(u,u)$ vanishes. For the second term of $B_\delta(u,u)$, with Lemma \ref{Mdelta} we have $M_\delta(\mathbf{x})\geq 0$, and therefore $h_\delta(u,u)=\int_{\Omega_\delta} M_\delta(\mathbf{x}) \int_{-\delta}^{\delta} H_{\delta}(|l|) [u(\mathbf{x}_{l})-u(\mathbf{x})]^2 d\mathbf{x}_{l}d\mathbf{x}\geq0$. %
We then have the Poincare-type inequality: there exists constants $C$ and $\delta_0$ such that for all $u\in S_\delta$ and $\delta\leq\delta_0$:
$$||u||_{L^2(\Omega)}^2\leq C^* b_\delta(u,u)\leq C^*(b_\delta(u,u)+h_\delta(u,u))=C^*B_\delta(u,u).$$
We now proceed to finish the proof. 
Here we assume that $\|u\|_{L^2(\Omega)}>0$, {otherwise the result is trivial}. For simplicity, we now denote $\delta_1$ as $\min(\delta_0,\overline{\delta})$ where $\overline{\delta}$ is defined in Lemma \ref{Mdelta} and $\delta_0$ as in \eqref{eqn:C*}. With \eqref{eqn:C*} and Lemma \ref{Mdelta} we  still have $||u||_{L^2(\Omega)}^2\leq C^* b_\delta(u,u)$ and $h_\delta(u,u)\geq 0$. {We now proceed to estimate the last term in $B_\delta(u,u)$:
\begin{align}
\nonumber&\int_{\Omega_\delta} \int_{-\delta}^{\delta} \left[M_\delta(\mathbf{x}_{l})\dfrac{|\mathbf{r}'(\mathbf{x})|}{|\mathbf{r}'(\mathbf{x}_l)|}-M_\delta(\mathbf{x})\right] H_{\delta}(|l|) [u(\mathbf{x}_{l})-u(\mathbf{x})]%
 d\mathbf{x}_l u(\mathbf{x}) d\mathbf{x}\\
\nonumber\geq& - \dfrac{1}{2}\int_{\Omega_\delta} M_\delta(\mathbf{x})  \int_{-\delta}^{\delta} H_{\delta}(|l|) [u(\mathbf{x}_{l})-u(\mathbf{x})]^2 d\mathbf{x}_l
d\mathbf{x} \\
\nonumber&-\dfrac{1}{2}\int_{\Omega_\delta} \int_{-\delta}^{\delta} H_{\delta}(|l|)\dfrac{|M_\delta(\mathbf{x})-M_\delta(\mathbf{x}_l)|\mathbf{r}'(\mathbf{x})|/|\mathbf{r}'(\mathbf{x}_l)||^2}{|M_\delta(\mathbf{x})|}d\mathbf{x}_l %
 |u(\mathbf{x})|^2d\mathbf{x}\\
 \nonumber\geq& - \dfrac{1}{2}h_\delta(u,u) -\dfrac{C_NC_M }{2}
 \int_{\Omega_\delta} \int_{-\delta}^{\delta} H_{\delta}(|l|)\delta^3 d\mathbf{x}_l |u(\mathbf{x})|^2d\mathbf{x}\\
 \geq& - \dfrac{1}{2}h_\delta(u,u) -\dfrac{C_NC_M C_H \delta}{2}
 \|u\|_{L^2(\Omega)}^2.\label{eqn:lastterm}
\end{align}}
Hence, when
\begin{equation}
\delta<\min\left\{\delta_1,\frac{1}{C^*C_NC_MC_H}\right\}=:\tilde{\delta}  \label{dt}
\end{equation}
we have
{\begin{align*}
B_\delta(u,u)
&\geq \left(\dfrac{1}{C^*}-\dfrac{C_NC_M C_H \delta}{2}\right)\|u\|_{L^2(\Omega)}^2+\dfrac{1}{2}h_\delta(u,u)\geq \dfrac{1}{2C^*}\|u\|_{L^2(\Omega)}^2.
\end{align*}}
\end{proof}

The uniform boundedness of $L_\delta^{-1}$ then follows
\begin{lem}\label{thm:boundL}
{ Assuming that $\Omega$ and $\tilde{\delta}$ satisfy the conditions in Lemma \ref{lemma:poincare}, there exists a constant $C$ such that
\begin{equation}
 ||L_\delta^{-1}||_{L^2(\Omega)}\leq C.
\end{equation}
}
\end{lem}

Moreover, with the definition of $||\cdot||_{S_\delta}$, we can show the boundedness and coercivity of the nonsymmetric bilinear operator $B_\delta(\cdot,\cdot)$: 
\begin{lem}\label{thm:boundcoer}
{ 
There exists a $0<\tilde{\delta}\leq 1$ such that for all $\delta<\tilde{\delta}$ the following inequalities hold
\begin{equation}\label{eqn:boundedness}
 \forall u,v\in S_{\delta}, \; B_\delta(u,v)\leq C_1 ||u||_{S_\delta}||v||_{S_\delta},
\end{equation}
\begin{equation}\label{eqn:coercivity}
 \forall u\in S_{\delta}, \; B_\delta(u,u)\geq C_2 ||u||^2_{S_\delta},
\end{equation}
for two constants $C_1,C_2>0$.
}
\end{lem}

\begin{proof} 
{
We first show \eqref{eqn:boundedness}. For the first two terms in $B_\delta(u,v)$, with the Cauchy-Schwarz inequality one may obtain
$b_\delta(u,v)\leq C\sqrt{ b_\delta(u,u)b_\delta(v,v)}$ and $h_\delta(u,v)\leq C\sqrt{h_\delta(u,u)h_\delta(v,v)}$. 
Moreover, with Lemma \ref{Mdelta}, similar as in \eqref{eqn:lastterm} we can show that  
\begin{align}
 \nonumber&\int_{\Omega_\delta} \int_{-\delta}^{\delta} \left[M_\delta(\mathbf{x}_{l})\frac{|\mathbf{r}'(\mathbf{x})|}{|\mathbf{r}'(\mathbf{x}_l)|}-M_\delta(\mathbf{x})\right] H_{\delta}(|l|) [u(\mathbf{x}_{l})-u(\mathbf{x})]%
 d\mathbf{x}_l v(\mathbf{x}) d\mathbf{x}\\
\leq &C\sqrt{h_\delta(u,u) \left(C_NC_MC_H \delta \|v\|_{L^2(\Omega)}^2\right)} \leq C\sqrt{h_\delta(u,u) b_\delta (v,v)}.\label{eqn:term3}
\end{align}
Therefore
 \begin{align*}
 B_\delta(u,v)^2 
 &\le C(b_\delta(u,u)b_\delta(v,v)+h_\delta(u,u) b_\delta (v,v)+h_\delta(u,u)h_\delta(v,v))\leq C ||u||^2_{S_\delta}||v||^2_{S_\delta}.
  \end{align*}
On the other hand, \eqref{eqn:coercivity} can be obtained when $\tilde{\delta}$ is taken as in \eqref{dt} and follow a similar proof as in Lemma \eqref{lemma:poincare}. 
}
\end{proof}

With the above properties, we can see that there exists a unique solution $u_\delta\in S_\delta$ solving \eqref{eqn:weakformula} (cf, \cite[Theorem~2.5.6]{brenner2007mathematical}). The well-posedness of the proposed variational problem is therefore obtained. To further show the asymptotic property of solution when $\delta\rightarrow 0$, we need the following embedding property:
\begin{lem}\label{lemma:poincare2}
For all $u\in S_{0}$ there exists a constant $C$ such that
\begin{equation}
 B_\delta(u,u)\leq C||\nabla u||_{L^2(\Omega)}^2
\end{equation}
for any $\delta$ satisfying the condition in Lemma \ref{thm:boundcoer}.
\end{lem}

\begin{proof} 
Given $u\in S_0$, from \cite[Theorem~1]{bourgain2001another} we have that
 \[b_\delta(u,u)\le C \|u\|_{H^1(\Omega)}^2\le C \|\nabla u\|_{L^2(\Omega)}^2.\]
To bound the second and the third terms of $B_\delta(u,u)$, we start with the case of boundary curvature$\equiv 0$, where we only need to %
show that $h_\delta(u,u)=\int_{\Omega_\delta} M_\delta(\mathbf{x}) \int_{-\delta}^{\delta} H_{\delta}(|l|) [u(\mathbf{x}_{l})-u(\mathbf{x})]^2 %
 d\mathbf{x}_{l}d\mathbf{x}\leq C||\nabla u||_{L^2(\Omega)}^2$. Since %
$M_\delta(\mathbf{x})\leq C$, 
it suffices to estimate $\int_{\Omega_\delta}\int_{-\delta}^{\delta} H_{\delta}(|l|) [u(\mathbf{x}_{l})-u(\mathbf{x})]^2 d\mathbf{x}_{l} d\mathbf{x}$. {With the H\"older inequality and the fact that $\int_{\Omega_\delta} |\nabla u(\mathbf{x}_{t})|^2 d\mathbf{x}=\int_{\Omega_\delta} |\nabla u(\mathbf{x})|^2 d\mathbf{x}$ for %
all $|t|\leq \delta$, we have
\begin{align*}
 &\int_{\Omega_\delta}\int_{-\delta}^{\delta} H_{\delta}(|l|) [u(\mathbf{x}_{l})-u(\mathbf{x})]^2 d\mathbf{x}_{l} d\mathbf{x}\\
 \leq &\sup_{r\leq1} H(|r|)\dfrac{1}{\delta^3}\int_{\Omega_\delta}\int_{-\delta}^{\delta} [u(\mathbf{x}_{l})-u(\mathbf{x})]^2 d\mathbf{x}_{l} d\mathbf{x}
 \leq \dfrac{C}{\delta^2}\int_{\Omega_\delta}\int_{-\delta}^{\delta} \int_{0}^l |\nabla u(\mathbf{x}_{t})|^2 dt d\mathbf{x}_{l} d\mathbf{x}\\
 \leq &\dfrac{C}{\delta^2}\int_{\Omega_\delta}\int_{-\delta}^{\delta} \int_{0}^{\delta} |\nabla u(\mathbf{x}_{t})|^2 dt d\mathbf{x}_{l} d\mathbf{x}%
 =\dfrac{C}{\delta}\int_{\Omega_\delta}\int_{0}^{\delta} |\nabla u(\mathbf{x}_{t})|^2 dt d\mathbf{x}\\
 =&\dfrac{C}{\delta}\int_{0}^{\delta}\int_{\Omega_\delta} |\nabla u(\mathbf{x}_{t})|^2 d\mathbf{x}dt %
=C\int_{\Omega_\delta} |\nabla u(\mathbf{x})|^2 d\mathbf{x}\leq C||\nabla u||_{L^2(\Omega)}^2.
\end{align*}
Therefore, the Lemma holds true when the boundary curvature $\kappa(\mathbf{x})\equiv 0$, a.e..} We now work on the case of nonzero curvature. 
Similar as in the curvature$\equiv 0$ case we can obtain $h_\delta(u,u)\leq C\|\nabla u\|_{L^2(\Omega)}^2$. For the last term of $B_\delta(u,u)$, with \eqref{eqn:term3} we have
\begin{align*}
&\int_{\Omega_\delta} \int_{-\delta}^{\delta} \left[M_\delta(\mathbf{x}_{l})\dfrac{|\mathbf{r}'(\mathbf{x})|}{|\mathbf{r}'(\mathbf{x}_l)|}-M_\delta(\mathbf{x})\right] H_{\delta}(|l|) [u(\mathbf{x}_{l})-u(\mathbf{x})]%
 d\mathbf{x}_l u(\mathbf{x}) d\mathbf{x}\\
 \leq &C\sqrt{h_\delta(u,u) b_\delta (u,u)}\leq C\|\nabla u\|_{L^2(\Omega)}^2.
\end{align*}
\end{proof}

Before studying the limiting behavior of the nonlocal operator, we need a compactness property:
\begin{lem}\label{thm:precompact}
Suppose $u_n\in S_{\delta_n}$ and $\delta_n\rightarrow 0$, then given $\sup_n B_{\delta_n}(u_n,u_n)\leq \infty$, %
$u_n$ is precompact in $L^2(\Omega)$. Moreover, any limit point $u\in S_0$.
\end{lem}

\begin{proof}
Since $S_{\delta_n} \subseteq L^2(\Omega)$ and $h_\delta(u_n,u_n)\geq 0$, {similar to \eqref{eqn:lastterm} we have,
\begin{align*}
B_\delta(u,u)
\geq &\dfrac{1}{2}\int_{\Omega}\int_{\Omega} J_{\delta_n}(\mathbf{x}-\mathbf{y})(u_n(\mathbf{y})-u_n%
(\mathbf{x}))^2d\mathbf{y}d\mathbf{x}
+\left(\dfrac{1}{2C^*}-\dfrac{C_NC_M C_H \delta}{2}\right)\|u\|_{L^2(\Omega)}^2
\end{align*}
where $C^*$ denotes the constant in \eqref{eqn:C*}. %
Therefore, when $\tilde{\delta}$ is taken as in \eqref{dt}, then for all $\delta<\tilde{\delta}$
\begin{displaymath}
B_{\delta_n}(u_n,u_n)\geq \dfrac{1}{2}\int_{\Omega}\int_{\Omega} J_{\delta_n}(\mathbf{x}-\mathbf{y})(u_n(\mathbf{y})-u_n%
(\mathbf{x}))^2d\mathbf{y}d\mathbf{x}.
\end{displaymath}}
We have $u_n\in L^2(\Omega)$ and
\begin{displaymath}
 \int_{\Omega}\int_{\Omega} J_{\delta_n}(\mathbf{x}-\mathbf{y})(u_n(\mathbf{y})-u_n%
(\mathbf{x}))^2d\mathbf{y}d\mathbf{x} \leq \infty.
\end{displaymath}
From \cite[Theorem~1.2]{ponce2004estimate}, any limit of $\{u_n\}$ is in ${L^2(\Omega)}$, or equivalently, %
$u_n$ is precompact, and any limit point $u\in S_0$.
\end{proof}

With the above lemmas, we obtain the following $L^2$ convergence result for an intermediate solution as $\delta\rightarrow 0$:
\begin{lem}\label{thm:l2conv}
Suppose $\tilde{u}_\delta$ is the weak solution of
 \begin{equation}\label{eqn:nonlocaleqn_m}
 \left\{\begin{array} {cl}
  L_{N\delta} \tilde{u}_\delta=f,&\text{ in }\Omega\\
  \int_\Omega \tilde{u}_\delta d\mathbf{x}=0,&\\
 \end{array}\right.
\end{equation}
and $u_0$ is the weak solution of \eqref{eqn:localeqn}, then
\begin{equation}\label{eqn:L2conv}
 \lim_{\delta\rightarrow 0}||\tilde{u}_\delta-u_0||_{L^2(\Omega)}=0.
\end{equation}
\end{lem}

\begin{proof}
The proof follows a similar strategy as in \cite{tian2014asymptotically,tao2017nonlocal}. A detailed derivation is provided in Appendix \ref{App:l2conv}.
\end{proof}

We now have the main theorem of this section for $f\in C(\overline{\Omega})$:
\begin{thm}\label{thm:l2main}
 Suppose $u_\delta$ is the weak solution of \eqref{eqn:nonlocaleqn} and $u_0$ is the weak solution of \eqref{eqn:localeqn}, then
\begin{equation}\label{eqn:L2conv_final}
 \lim_{\delta\rightarrow 0}||{u}_\delta-u_0||_{L^2(\Omega)}=0.
\end{equation}
\end{thm}

\begin{proof}
With the results in Lemma \ref{thm:l2conv}, we only need to show that $\lim_{\delta\rightarrow 0}||{u}_\delta-\tilde{u}_\delta||_{L^2(\Omega)}=0$. %
Since $L_\delta({u}_\delta-\tilde{u}_\delta)=f_\delta-f$, with Lemma \ref{thm:boundL} we can see that it suffices to show
\begin{equation}
 \lim_{\delta\rightarrow 0}||f_\delta-f||_{L^2(\Omega)}=0,
\end{equation}
or equivalently
\begin{equation}
 \nonumber\lim_{\delta\rightarrow 0}\int_{\Omega_\delta}\left(\int_{(\mathbb{R}^N\backslash\Omega)_\delta} J_{\delta}(|\mathbf{x}-\mathbf{y}|)%
 (|(\mathbf{y}-\overline{\mathbf{x}})\cdot\mathbf{n}(\overline{\mathbf{x}})|^2-|(\mathbf{x}-\overline{\mathbf{x}})\cdot\mathbf{n}%
 (\overline{\mathbf{x}})|^2)d\mathbf{y}f(\mathbf{x})\right)^{2}d\mathbf{x}=0.
\end{equation}
Since
\begin{align}
 \nonumber&\left|\int_{(\mathbb{R}^N\backslash\Omega)_\delta} J_{\delta}(|\mathbf{x}-\mathbf{y}|)%
 (|(\mathbf{y}-\overline{\mathbf{x}})\cdot\mathbf{n}(\overline{\mathbf{x}})|^2-|(\mathbf{x}-\overline{\mathbf{x}})\cdot\mathbf{n}%
 (\overline{\mathbf{x}})|^2)d\mathbf{y}\right|\\
 \nonumber\leq &C\int_{(\mathbb{R}^N\backslash\Omega)_\delta} J_{\delta}(|\mathbf{x}-\mathbf{y}|)|(\mathbf{y}-\mathbf{x})\cdot\mathbf{n}%
 (\overline{\mathbf{x}})|^2d\mathbf{y}\leq C,
\end{align}
we have
\begin{align*}
 \nonumber&\int_{\Omega_\delta}\left(\int_{(\mathbb{R}^N\backslash\Omega)_\delta} J_{\delta}(|\mathbf{x}-\mathbf{y}|)%
 (((\mathbf{y}-\overline{\mathbf{x}})\cdot\mathbf{n}(\overline{\mathbf{x}}))^2-((\mathbf{x}-\overline{\mathbf{x}})\cdot\mathbf{n}%
 (\overline{\mathbf{x}}))^2)d\mathbf{y}f(\mathbf{x})\right)^{2}d\mathbf{x}
 \leq C\int_{\Omega_\delta}|f(\mathbf{x})|^{2}d\mathbf{x}
\end{align*}
which vanishes as $\delta\rightarrow 0$. 
\end{proof}

\begin{rem}\label{remark:inhomo}
For the analysis in this paper we focus on the homogeneous Neumann-type boundary condition $g(\mathbf{x})=0$, while we note that the proposed nonlocal variational formulation can be applied to inhomogeneous boundary conditions. Here we take $J_\delta(r)=\dfrac{4}{\pi \delta^4}$ for simplicity. When $f(\mathbf{x})=0$ and $g(\mathbf{x})\neq 0$, applying a test function $v(\mathbf{x})\in C^{\infty}(\Omega)$ to \ref{eqn:2.8} yields
\begin{align*}
 (f_\delta,v)_{L^2(\Omega)}=&\int_\Omega\left(2\int_{\mathbb{R}^N\backslash\Omega} J_{\delta}(|\mathbf{x}-\mathbf{y}|)(\mathbf{y}-\mathbf{x})\cdot\mathbf{n}(\overline{\mathbf{x}})%
 d\mathbf{y}-M_\delta(\mathbf{x})\kappa(\overline{\mathbf{x}})\right)g(\overline{\mathbf{x}})v(\mathbf{x})d\mathbf{x}\\
 =&2\int_{\Omega_\delta}\int_{\mathbb{R}^N\backslash\Omega} J_{\delta}(|\mathbf{x}-\mathbf{y}|)(\mathbf{y}-\mathbf{x})\cdot\mathbf{n}(\overline{\mathbf{x}})%
 d\mathbf{y}g(\overline{\mathbf{x}})v(\mathbf{x})d\mathbf{x}
 -\int_{\Omega_\delta} M_\delta(\mathbf{x})\kappa(\overline{\mathbf{x}})g(\overline{\mathbf{x}})v(\mathbf{x})d\mathbf{x}.
\end{align*}
For the second part, with the H\"older inequality we have
\begin{align*}
\left|\int_{\Omega_\delta} M_\delta(\mathbf{x})\kappa(\overline{\mathbf{x}})g(\overline{\mathbf{x}})v(\mathbf{x})d\mathbf{x}\right|
\leq C\sqrt{\int_{\Omega_\delta}g^2(\overline{\mathbf{x}})d\mathbf{x}\int_{\Omega_\delta}v^2(\mathbf{x})d\mathbf{x}} \leq  C\sqrt{\delta ||g||_{L^2(\partial\Omega)}||v||_{L^2(\Omega)}}.
\end{align*}
Therefore, $\left|\int_{\Omega_\delta} M_\delta(\mathbf{x})\kappa(\overline{\mathbf{x}})g(\overline{\mathbf{x}})v(\mathbf{x})d\mathbf{x}\right|\rightarrow 0$ as $\delta\rightarrow 0$. To show the asymptotic limit for the first part as $\delta\rightarrow 0$, for each $\mathbf{x}\in\Omega_\delta$ we have
\begin{align*}
\int_{\mathbb{R}^N\backslash\Omega} J_{\delta}(|\mathbf{x}-\mathbf{y}|)(\mathbf{y}-\mathbf{x})\cdot\mathbf{n}(\overline{\mathbf{x}})%
 d\mathbf{y}
 =&\int_{D_\delta} J_{\delta}(|\mathbf{x}-\mathbf{y}|)(\mathbf{y}-\mathbf{x})\cdot\mathbf{n}(\overline{\mathbf{x}})%
 d\mathbf{y}+\int_{A_\delta} J_{\delta}(|\mathbf{x}-\mathbf{y}|)(\mathbf{y}-\mathbf{x})\cdot\mathbf{n}(\overline{\mathbf{x}})%
 d\mathbf{y}.
\end{align*}
For the first part
\[\int_{D_\delta} J_{\delta}(|\mathbf{x}-\mathbf{y}|)(\mathbf{y}-\mathbf{x})\cdot\mathbf{n}(\overline{\mathbf{x}})%
 d\mathbf{y}=\dfrac{4}{\pi \delta^4}\int_{s_x}^\delta 2r\sqrt{\delta^2-r^2}dr=\dfrac{8}{3\pi \delta^4}(\delta^2-s_x^2)^{3/2}\]
and
\begin{align}
\nonumber&2\int_{\Omega_\delta}\int_{D_\delta} J_{\delta}(|\mathbf{x}-\mathbf{y}|)(\mathbf{y}-\mathbf{x})\cdot\mathbf{n}(\overline{\mathbf{x}})%
 d\mathbf{y}g(\overline{\mathbf{x}})v(\mathbf{x})d\mathbf{x}\\
\nonumber=&\dfrac{16}{3\pi \delta^4}\int_{\Omega_\delta}(\delta^2-s_x^2)^{3/2}g(\overline{\mathbf{x}})v(\mathbf{x})d\mathbf{x}=\dfrac{16}{3\pi \delta^4}\int_{\partial\Omega}\int_0^\delta(\delta^2-r^2)^{3/2}dr g(\overline{\mathbf{x}})v(\overline{\mathbf{x}})d\overline{\mathbf{x}}+\text{O}(\delta)\\
=&\int_{\partial\Omega} g(\overline{\mathbf{x}})v(\overline{\mathbf{x}})d\overline{\mathbf{x}}+\text{O}(\delta).\label{eqn:nonhomo_1}
\end{align}
For the second part, since the area of $A_\delta$ is bounded by $C\delta^3$, we have
\begin{equation}\label{eqn:nonhomo_2}
\left|\int_{A_\delta} J_{\delta}(|\mathbf{x}-\mathbf{y}|)(\mathbf{y}-\mathbf{x})\cdot\mathbf{n}(\overline{\mathbf{x}})%
 d\mathbf{y}\right|\leq C\delta.
\end{equation}
Combining \eqref{eqn:nonhomo_1} and \eqref{eqn:nonhomo_2} yields
\begin{align*}
\lim_{\delta\rightarrow 0}(f_\delta,v)_{L^2(\Omega)}=&\lim_{\delta\rightarrow 0}2\int_{\Omega_\delta}\int_{\mathbb{R}^N\backslash\Omega} J_{\delta}(|\mathbf{x}-\mathbf{y}|)(\mathbf{y}-\mathbf{x})\cdot\mathbf{n}(\overline{\mathbf{x}})%
 d\mathbf{y}g(\overline{\mathbf{x}})v(\mathbf{x})d\mathbf{x}
 =\int_{\partial\Omega} g(\overline{\mathbf{x}})v(\overline{\mathbf{x}})d\overline{\mathbf{x}}.
\end{align*}
Therefore, the right hand side converges to the inhomogeneous flux condition as $\delta\rightarrow 0$ in the variational formulation. In fact, the asymptotic convergence property in Theorem \ref{thm:l2main} can be shown for the nonlocal diffusion problem with inhomogeneous flux conditions given the corresponding nonlocal trace theorem, which will be addressed in the future work.
\end{rem}

\section{Convergence rate in the $L^{\infty}(\Omega)$ norm}\label{sec:linfty}

In this section we will estimate the order of convergence rate by considering a problem with the more general setting: $\partial\Omega=\partial\Omega_D\bigcup \partial\Omega_N$ %
and $(\partial\Omega_D)^o\bigcap (\partial\Omega_N)^o=\emptyset$. %
Here $\partial \Omega_D$ and $\partial \Omega_N$ are both 1D curves. %
To define a Dirichlet-type constraint on $\partial \Omega_D$, we denote $\Omega_{D\delta}=\{\mathbf{x}\in\Omega_\delta:%
\overline{\mathbf{x}}\in\partial\Omega_D\}$ where $\overline{\mathbf{x}}$ is the orthogonal projection of $\mathbf{x}$ on $\partial\Omega$. %
Moreover, we denote $\partial\Omega_{D\delta}=\{\mathbf{x}\in{\mathbb{R}^N\backslash\Omega}:%
\text{dist}(\mathbf{x},\Omega_{D\delta})\leq\delta\}$ and assume that the value of $u$ is given on it. To be specific, %
here we assume $u(\mathbf{x})=0$ on $\partial\Omega_{D\delta}$ without loss of generality. %
Similarly, to apply the Neumann-type constraint on $\partial \Omega_N$, we denote $\Omega_{N\delta}=\{\mathbf{x}\in\Omega_\delta:%
\overline{\mathbf{x}}\in\partial\Omega_N\}$ and $\partial\Omega_{N\delta}=\{\mathbf{x}\in{\mathbb{R}^N\backslash\Omega}:%
\text{dist}(\mathbf{x},\Omega_{N\delta})\leq\delta\}$. %
We consider a Neumann-type constraint as an extension of $\dfrac{\partial u}{\partial \mathbf{n}}=g(\mathbf{x})$ on $\partial \Omega_N$, %
by modifying the nonlocal problem discussed in the last section as follows: for $\mathbf{x}\in\Omega\backslash\Omega_{N\delta}$:
\begin{align}
 \nonumber&-2\int_{\Omega\cup\partial\Omega_{D\delta}} J_{\delta}(|\mathbf{x}-\mathbf{y}|)(u_{\delta}(\mathbf{y})-u_{\delta}(\mathbf{x}))d\mathbf{y}=f(\mathbf{x}),
\end{align}
and for $\mathbf{x}\in\Omega_{N\delta}$:
\begin{align}
 \nonumber&-2\int_{\Omega} J_{\delta}(|\mathbf{x}-\mathbf{y}|)(u_{\delta}(\mathbf{y})-u_{\delta}(\mathbf{x}))d\mathbf{y}-%
 2M_\delta(\mathbf{x})\int_{-\delta}^{\delta} H_{\delta}(|l|) (u_{\delta}(\mathbf{x}_{l})-u_{\delta}(\mathbf{x})) d\mathbf{x}_{l}\\
 \nonumber=&f(\mathbf{x})-\int_{\partial\Omega_{N\delta}} J_{\delta}(|\mathbf{x}-\mathbf{y}|)%
 \left[|(\mathbf{y}-\overline{\mathbf{x}})\cdot\mathbf{n}(\overline{\mathbf{x}})|^2-|(\mathbf{x}-\overline{\mathbf{x}})\cdot\mathbf{n}%
 (\overline{\mathbf{x}})|^2\right]d\mathbf{y}f({{\mathbf{x}}})\\
 &+\left(2\int_{\partial\Omega_{N\delta}} J_{\delta}(|\mathbf{x}-\mathbf{y}|)(\mathbf{y}-\mathbf{x})\cdot\mathbf{n}(\overline{\mathbf{x}})%
 d\mathbf{y}- M_\delta(\mathbf{x})\kappa(\overline{\mathbf{x}})\right)g(\overline{\mathbf{x}}),\label{eqn:formula2mix}
\end{align}
where
\[M_\delta(\mathbf{x}):=\int_{\partial\Omega_{N\delta}} J_{\delta}(|\mathbf{x}-\mathbf{y}|)\left[|(\mathbf{y}-\mathbf{x})\cdot\mathbf{p}(\overline{\mathbf{x}})|^2-|(\mathbf{y}-\overline{\mathbf{x}})\cdot\mathbf{n}(\overline{\mathbf{x}})|^2+|(\mathbf{x}-\overline{\mathbf{x}})\cdot\mathbf{n}(\overline{\mathbf{x}})|^2\right]%
 d\mathbf{y} \]
Here we note that it is possible that $\Omega_{D\delta}\bigcap\Omega_{N\delta}\neq\emptyset$. We can then rewrite the nonlocal equation to be solved as
\begin{equation}\label{eqn:nonlocalmix}
 \left\{\begin{array}{cl}
L_\delta u =f,&\text{ on }\Omega\backslash\Omega_{N\delta}\\
L_{N\delta} u =f_\delta,&\text{ on }\Omega_{N\delta}\\
u=0,&\text{ on }\partial\Omega_{D\delta}.
        \end{array}
 \right.
\end{equation}
The corresponding limiting local model is given by
\begin{equation}\label{eqn:localmix}
 \left\{\begin{array}{cl}
-\triangle u =f,&\text{ on }\Omega\\
\dfrac{\partial u}{\partial \mathbf{n}} =g,&\text{ on }\partial\Omega_{N}\\
u=0,&\text{ on }\partial\Omega_{D}.
        \end{array}
 \right.
\end{equation}
In this section we focus on the case with homogeneous Neumann-type constraints, i.e., $g(\mathbf{x})=0$.

For the above problem with mixed constraints, we have the nonlocal maximum principle stated below
\begin{lem}\label{eqn:maxprinciple}
{ For $u\in C(\overline{\Omega})\cap C(\partial\Omega_{D\delta}\backslash\partial\Omega_D)$ and $u$ bounded on $\partial\Omega_{D\delta}$, 
assuming that $u$ satisfies $L_\delta u\leq 0$ for all $x\in \Omega\backslash\Omega_{N\delta}$ and $L_{N\delta} u\leq 0$ for all $x\in \Omega_{N\delta}$, we have
\begin{equation}
 \sup_{\mathbf{x}\in\overline{\Omega}\cup \partial\Omega_{D\delta}} u(\mathbf{x})\leq \sup_{\mathbf{x}\in \partial\Omega_{D\delta}} u(\mathbf{x}). 
\end{equation}
}
\end{lem}

\begin{proof}
Assuming that $\sup_{\mathbf{x}\in\overline{\Omega}\cup \partial\Omega_{D\delta}}u(\mathbf{x})>\sup_{\mathbf{x}\in\partial\Omega_{D\delta}}u(\mathbf{x})$, since $u\in C(\overline{\Omega})$ there exists $\mathbf{x}^*\in(\Omega\cup\partial\Omega_N)$ such that $u(\mathbf{x}^*)=\sup_{\mathbf{x}\in(\overline{\Omega}\cup \partial\Omega_{D\delta})}u(\mathbf{x})$. 

{\em Case 1: $\mathbf{x}^*\in \Omega \setminus \Omega_{N\delta}$.} Then
$L_\delta u(\mathbf{x}^*) =-2\int_{\Omega\cup \partial \Omega_{D\delta}} J_\delta(|\mathbf{x}^*
-\mathbf{y}|) (u(\mathbf{y})-u(\mathbf{x}^*)) d\mathbf{y}\ge 0$. Therefore $L_\delta u(\mathbf{x}^*) =0$ and
\begin{equation}
 u(\mathbf{y})=u(\mathbf{x}^*)=\sup_{\mathbf{x}\in\overline{\Omega}\cup \partial \Omega_{D\delta}} u(\mathbf{x}), \qquad \forall \mathbf{y}\in 
(\overline{\Omega}\cup \partial \Omega_{D\delta})\cap B(\mathbf{x}^*,\delta).
\label{x0}
\end{equation}

{\em Case 2: $\mathbf{x}^*\in \Omega_{N\delta}$.} Then
\begin{align*}
L_{N\delta} u(\mathbf{x}^*) =&-2\int_{\Omega} J_\delta(|\mathbf{x}^*
-\mathbf{y}|) (u(\mathbf{y})-u(\mathbf{x}^*)) d\mathbf{y}
-2M_\delta(\mathbf{x}^*){\int_{-\delta}^{\delta} H(|l|) [u(\mathbf{x}^*_l)-u(\mathbf{x}^*)] d \mathbf{x}^*_l
}
\ge 0.
\end{align*}
 Note that in Lemma \ref{Mdelta} we have proven $M_\delta(\mathbf{x}^*) \ge 0$. 
{
Again, this is possible only when
\begin{equation}
 u(\mathbf{y})=u(\mathbf{x}^*)=\sup_{\mathbf{x}\in\overline{\Omega}\cup \partial \Omega_{D\delta}} u(\mathbf{x}), \qquad \forall \mathbf{y}\in 
\overline{\Omega}\cap B(\mathbf{x}^*,\delta).
\label{xN0}
\end{equation}
Summing up the two cases, in view of \eqref{x0} and \eqref{xN0}, we have
\begin{align}
\mathbf{x}^*\in \Omega \setminus \Omega_{N\delta}&\Rightarrow u(\mathbf{y})=u(\mathbf{x}^*)=\sup_{\overline{\Omega}\cup \partial\Omega_{D\delta}} u,\; \forall \mathbf{y}\in 
(\overline{\Omega}\cup \partial\Omega_{D\delta})\cap B(\mathbf{x}^*,\delta),\label{max1}\\
\mathbf{x}^*\in \Omega_{N\delta}&\Rightarrow u(\mathbf{y})=u(\mathbf{x}^*)=\sup_{\overline{\Omega}\cup \partial\Omega_{D\delta}} u,\; \forall \mathbf{y}\in 
\overline{\Omega}\cap B(\mathbf{x}^*,\delta)\label{max2}.
\end{align}
Now fixing $\mathbf{y}^*\in ((\Omega\cup\partial\Omega_N)\cap B(\mathbf{x}^*,\delta))$, 
we can apply the same arguments with $\mathbf{y}^*$ in place of $\mathbf{x}^*$,
and get \eqref{max1} and \eqref{max2} with $\mathbf{y}^*$ in the role of $\mathbf{x}^*$.
This process can be repeated for all points $\mathbf{y}^*\in ((\Omega\cup\partial\Omega_N)\cap B(\mathbf{x}^*,\delta))$, and together with the continuity assumption of $u$ we obtain:
\begin{align*}
 &u(\mathbf{y})=u(\mathbf{x}^*)=\sup_{\overline{\Omega}\cup \partial \Omega_{D\delta}} u,\quad \forall \mathbf{y}\in 
(\overline{\Omega}\cup \partial \Omega_{D\delta})\cap \left[B(\mathbf{x}^*,\delta) \cup \left(\bigcup_{ \mathbf{y}^*\in (\Omega\cup\partial\Omega_N)\cap B(\mathbf{x}^*,\delta)} B( \mathbf{y}^*,\delta)\right)\right].
\end{align*}
Geometrically, note that
\begin{align*}
&(\overline{\Omega}\cup \partial \Omega_{D\delta})\cap \left[B(\mathbf{x}^*,\delta) \cup \left(\bigcup_{ \mathbf{y}^*\in (\Omega\cup\partial\Omega_N)
\cap B(\mathbf{x}^*,\delta)} B( \mathbf{y}^*,\delta)\right)\right]\\
=&\{z\in \overline{\Omega}\cup \partial \Omega_{D\delta}: \text{dist}(z,
(\Omega\cup\partial\Omega_N)\cap B(\mathbf{x}^*,\delta)  ) \le \delta\} .
\end{align*}
In other words, with this argument we expanded the region where 
$u(\mathbf{z})=\sup_{\overline{\Omega}\cup \partial \Omega_{D\delta}} u$ from $\mathbf{z}\in(\Omega\cup\partial\Omega_N)\cap B(\mathbf{x}^*,\delta)$ to its entire $\delta$-neighborhood lying in $\overline{\Omega}\cup \partial \Omega_{D\delta}$. 
We then apply this argument recursively, so that the region where $u(\mathbf{z})=\sup_{\overline{\Omega}\cup \partial \Omega_{D\delta}} u$  will get expanded to the entire domain of $\overline{\Omega}\cup \partial \Omega_{D\delta}$. In other words, to have a global maximum inside $\Omega$, the only possibility is for $u$ to be constant on $\overline{\Omega}\cup \partial \Omega_{D\delta}$, which contradicts with the assumption that $\sup_{\mathbf{x}\in\overline{\Omega}\cup \partial \Omega_{D\delta}}u(\mathbf{x})>\sup_{\mathbf{x}\in\partial\Omega_{D\delta}}u(\mathbf{x})$. 

} 
\end{proof}


We now assume that $u_\delta$ is the solution of \eqref{eqn:nonlocalmix} and $u_0$ is the solution of \eqref{eqn:localmix}. Denote $e_\delta(\mathbf{x}):={u}_\delta(\mathbf{x})-u_0(\mathbf{x})$, $T_\delta(\mathbf{x}):=(L_0u_0(\mathbf{x})-L_\delta u_0(\mathbf{x}))+(f_\delta(\mathbf{x})-f(\mathbf{x}))$ for $\mathbf{x}\in\Omega\backslash\Omega_{N\delta}$ and $T_\delta(\mathbf{x}):=(L_0u_0(\mathbf{x})-L_{N\delta} u_0(\mathbf{x}))+(f_\delta(\mathbf{x})-f(\mathbf{x}))$ for $\mathbf{x}\in\Omega_{N\delta}$, then for $\mathbf{x}\in \Omega\backslash\Omega_{N\delta}$,
\begin{equation}
 L_\delta e_\delta=L_\delta {u}_\delta-L_\delta u_0=L_0u_0-L_\delta u_0=T_\delta,
\end{equation}
and similarly for $\mathbf{x}\in \Omega_{N\delta}$,
\begin{equation}
 L_\delta e_\delta=L_{N\delta} {u}_\delta-L_{N\delta} u_0=f_\delta-f+L_0u_0-L_{N\delta} u_0=T_\delta.
\end{equation}
We then obtain the following truncation estimate for $T_\delta$:
\begin{lem}\label{lem:T}
{Suppose $u_0$ is the solution to local problem \eqref{eqn:localmix}, then
\begin{equation}
 T_\delta(\mathbf{x})=\textit{O}(\delta^2)
\end{equation}
for $\mathbf{x}\in \Omega\backslash\Omega_{N\delta}$, and
\begin{align}
 \nonumber T_\delta(\mathbf{x})=&2\int_{E_\delta} J_{\delta}(|\mathbf{x}-\mathbf{y}|)\dfrac{\partial u_0({\mathbf{x}})}{\partial \mathbf{p}}%
 ((\mathbf{x}-\mathbf{y})\cdot \mathbf{p}(\overline{\mathbf{x}}))d\mathbf{y}\\
 \nonumber&+\int_{\Omega} J_{\delta}(|\mathbf{x}-\mathbf{y}|)[u_0(\mathbf{x})]_{nnn}%
 ((\mathbf{x}-\mathbf{y})\cdot \mathbf{n}(\overline{\mathbf{x}})) (-|\overline{\mathbf{x}}-\mathbf{x}|^2+\dfrac{1}{3}|(\mathbf{x}-\mathbf{y})\cdot \mathbf{n}(\overline{\mathbf{x}})|^2)d\mathbf{y}\\
 \nonumber&+\int_{\Omega} J_{\delta}(|\mathbf{x}-\mathbf{y}|)[u_0(\mathbf{x})]_{npp}%
 ((\mathbf{x}-\mathbf{y})\cdot \mathbf{n}(\overline{\mathbf{x}}))|(\mathbf{x}-\mathbf{y})\cdot \mathbf{p}(\overline{\mathbf{x}})|^2d\mathbf{y}\\
 &+{\kappa(\overline{\mathbf{x}}) M_\delta({\mathbf{x}}) [u_0(\mathbf{x})]_{nn} %
((\mathbf{x} - \overline{\mathbf{x}}) \cdot \mathbf{n}(\mathbf{\overline{x}}))}
+O(\delta^2)
\end{align}
for $\mathbf{x}\in \Omega_{N\delta}$. Here $E_\delta$ denotes the region in $A_\delta$ which is asymmetric with respect to the $y$ axis (see the right plot of Figure \ref{graph}).
}
\end{lem}

\begin{proof}
The proof is based on the Taylor expansion of $u_0$ and an estimate for the asymmetric part in $A_\delta$. The detailed derivations can be found in Appendix \ref{App:T}.
\end{proof}

Furthermore, with the maximum principle, when $f\in C(\overline{\Omega})$ and $u_\delta=u_0$ continuous in $\partial\Omega_{D\delta}$, we have the following lemma
\begin{lem}\label{thm:phi}
{
Suppose that a nonnegative continuous function $\phi(\mathbf{x})$ is defined on $\overline{\Omega}\cup\partial\Omega_{D\delta}$, and $-L_\delta\phi\geq G(\mathbf{x})>0$ %
for $\mathbf{x}\in \Omega\backslash\Omega_{N\delta}$, $-L_{N\delta}\phi\geq G(\mathbf{x})>0$ %
for $\mathbf{x}\in \Omega_{N\delta}$. Then
\begin{equation}
\sup_{\mathbf{x}\in\Omega\cup\partial\Omega_N}|e_\delta(\mathbf{x})|\leq \sup_{\mathbf{x}\in\partial\Omega_{D\delta}}\phi(\mathbf{x}) \sup_{\mathbf{x}\in\Omega\cup\partial\Omega_N}\dfrac{|T_\delta(\mathbf{x})|}{G(\mathbf{x})}.
\end{equation}
}
\end{lem}

\begin{proof}
Let $K_\delta=\sup_{\mathbf{x}\in\Omega\cup\partial\Omega_N}\dfrac{|T_\delta(\mathbf{x})|}{G(\mathbf{x})}$, then for $K_\delta\phi(\mathbf{x})+e_\delta(\mathbf{x})$ %
we have: For $\mathbf{x}\in \Omega\backslash\Omega_{N\delta}$
\begin{align*}
 \nonumber L_\delta (K_\delta\phi(\mathbf{x})+e_\delta(\mathbf{x}))=&\sup_{\mathbf{x}\in\Omega\cup\partial\Omega_N}\dfrac{|T_\delta(\mathbf{x})|}{G(\mathbf{x})}L_\delta\phi(\mathbf{x})+L_\delta e_\delta(\mathbf{x})
 =\sup_{\mathbf{x}\in\Omega\cup\partial\Omega_N}\dfrac{|T_\delta(\mathbf{x})|}{G(\mathbf{x})}L_\delta\phi(\mathbf{x})+T_\delta \leq 0,
\end{align*}
and a similar argument holds for $\mathbf{x}\in \Omega_{N\delta}$. 
With the maximum principle in Lemma \ref{eqn:maxprinciple} we have
\begin{align*}
 \nonumber\sup_{\mathbf{x}\in\Omega\cup\partial\Omega_N} e_\delta(\mathbf{x})\leq& \sup_{\mathbf{x}\in\Omega\cup\partial\Omega_N}(K_\delta\phi(\mathbf{x})+e_\delta(\mathbf{x}))
 \leq \sup_{\mathbf{x}\in\partial\Omega_{D\delta}}(K_\delta\phi(\mathbf{x})+e_\delta(\mathbf{x}))%
 =K_\delta \sup_{\mathbf{x}\in\partial\Omega_{D\delta}} \phi(\mathbf{x}).
\end{align*}
Similarly, we have $L_\delta(K_\delta\phi(\mathbf{x})-e_\delta(\mathbf{x})) \leq0$ for $\mathbf{x}\in \Omega\backslash\Omega_{N\delta}$ and %
$L_{N\delta}(K_\delta\phi(\mathbf{x})-e_\delta(\mathbf{x})) \leq0$ for $\mathbf{x}\in \Omega_{N\delta}$, hence
\begin{align*}
 \nonumber\sup_{\mathbf{x}\in\Omega\cup\partial\Omega_N} (-e_\delta(\mathbf{x}))\leq& \sup_{\mathbf{x}\in\Omega\cup\partial\Omega_N}(K_\delta\phi(\mathbf{x})-e_\delta(\mathbf{x}))
 \leq  \sup_{\mathbf{x}\in\partial\Omega_{D\delta}}(K_\delta\phi(\mathbf{x})-e_\delta(\mathbf{x}))%
 =K_\delta \sup_{\mathbf{x}\in\partial\Omega_{D\delta}} \phi(\mathbf{x}).
\end{align*}

\end{proof}

 \begin{figure}[ht]\label{fig:4}
 \centering
\includegraphics[scale=0.4]{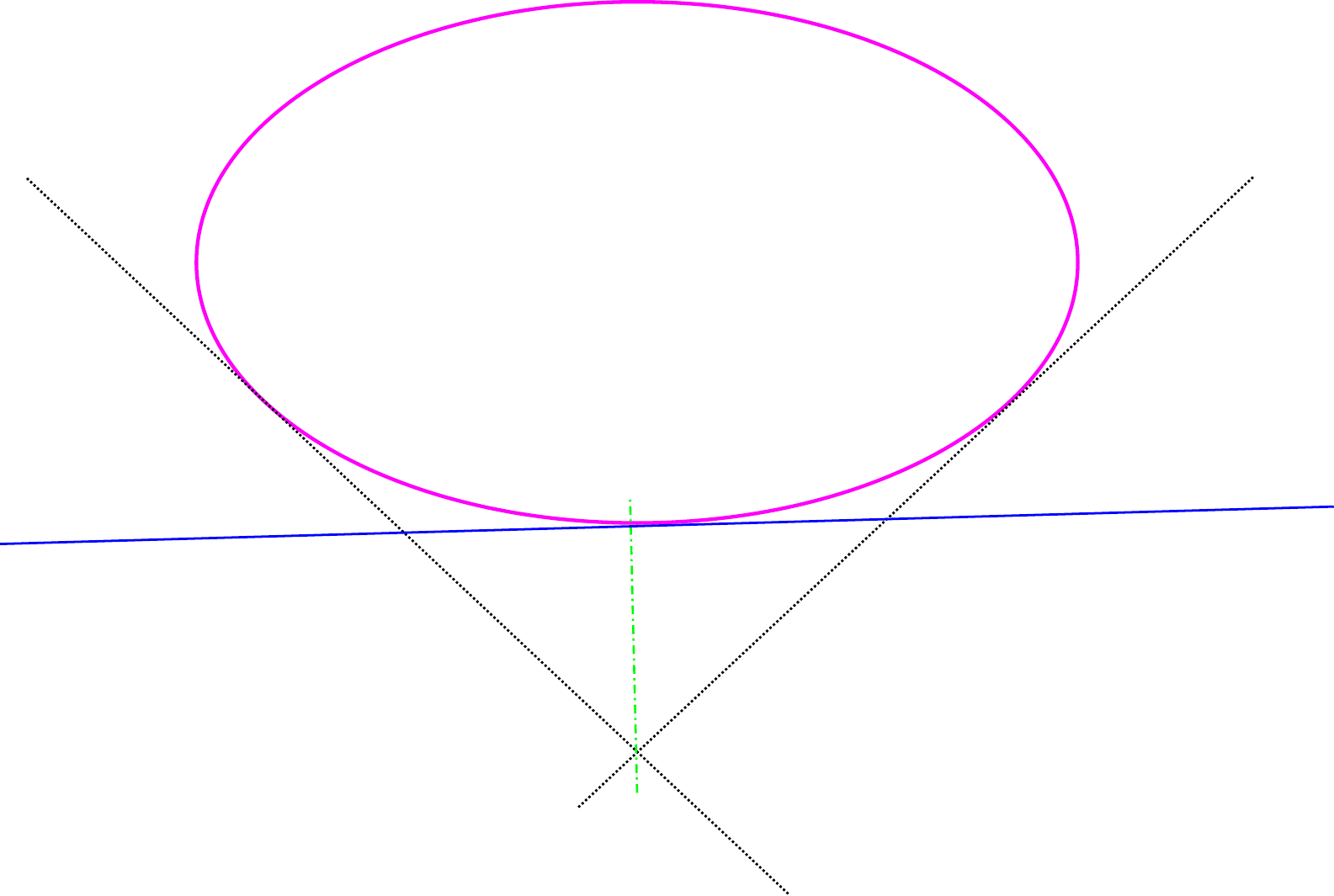}  
\put(-100,115){$\partial\Omega_D$}\put(-80,60){$\partial\Omega_N$}
\put(-155,60){$\mathbf{z}_1$}\put(-45,60){$\mathbf{z}_2$}
\put(-100,45){$\mathbf{z}_3$}
\put(-200,70){$\Pi$}
\put(-180,100){$\tau(\mathbf{z}_1)$}\put(-35,100){$\tau(\mathbf{z}_2)$}
\put(-25,45){$\tau(\mathbf{z}_3)$}
\put(-100,10){$\tilde{\mathbf{z}}$}
\caption{Illustration of the geometric assumption and notation for the barrier function $\phi(\mathbf{x})$ definition.}
 \end{figure}
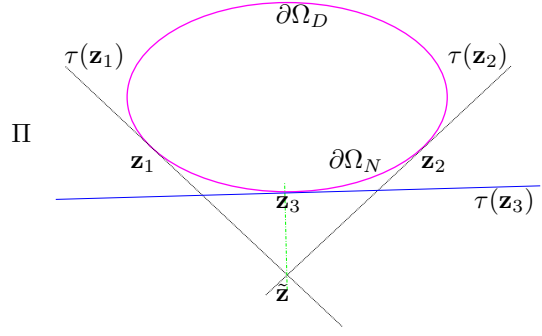
 
We now define a nonnegative continuous function $\phi$ satisfying the conditions given in Lemma \ref{thm:phi}.  In the following we take a specific kernel $J_{\delta}(\mathbf{s})=\dfrac{4}{\pi\delta^4}$ for $|\mathbf{s}|\leq\delta$ for simplicity.  
As shown in Figure \ref{fig:4}, let $\{\mathbf{z}_1,\mathbf{z}_2\}:= \partial \Omega_D\cap \partial \Omega_N$ and $\pi_{\partial\Omega}$ be the projection operator onto $\partial \Omega$.
  Due to the convexity of $\Omega$, the map $\pi_{\partial\Omega}(\mathbf{x})$ is always well defined and single-valued for any point $\mathbf{x}\notin\Omega$. For $\mathbf{x}\in\Omega$, the set where $\pi_{\partial\Omega}(\mathbf{x})$ is not single-valued (i.e. the ``ridge'' of $\partial\Omega$) is $L^2$-negligible \cite{mantegazza2003hamilton}. We then make the following crucial geometric assumption: Let $\tau(\mathbf{z}_1)$ (resp. $\tau(\mathbf{z}_2)$) be the tangent line to $\partial\Omega$
 at $\mathbf{z}_1$ (resp. $\mathbf{z}_2$), then the intersecting point $\tilde{\mathbf{z}}:=\tau(\mathbf{z}_1)\cap \tau(\mathbf{z}_2)$ satisfies
\begin{equation}
 \pi_{\partial\Omega}(\tilde{\mathbf{z}})\in \partial\Omega_N.
\end{equation}
Let $\mathbf{z}_3\in \partial\Omega$ be a point such that $\tau(\mathbf{z}_3)$ is orthogonal to the bisector of the angle $\angle \mathbf{z}_2 \tilde{\mathbf{z}}\mathbf{z}_1$. Set the barrier function as
\begin{equation}\label{eqn:phi} 
\phi(\mathbf{x}):=|\text{dist}(\mathbf{x},\tau(\mathbf{z}_3))+1|^2.
\end{equation}
For any point $\overline{\mathbf{x}}\in\partial\Omega_N$, in the following we denote the angle between $\mathbf{p}(\overline{\mathbf{x}})$ and $\mathbf{n}(\mathbf{z}_3)$ as $\alpha(\overline{\mathbf{x}})$. Note that with the crucial geometric assumption and the fact the $\Omega$ is convex, there exists $0<\tilde{\alpha}<\pi/2$ such that $\tilde{\alpha}\leq\alpha(\overline{\mathbf{x}})\leq \pi-\tilde{\alpha}$, $\forall\overline{\mathbf{x}}\in\partial\Omega_N$. Let $\Pi$ be the half-plane delimited by $\tau(\mathbf{z}_3)$ and containing $\Omega$, we now check the conditions in Lemma \ref{thm:phi} with the following 3 steps: 
 
 {\em Step 1. Convexity of $\phi$}. To check that $\phi$ is convex on $\Pi$, consider arbitrary points
 $\mathbf{x},\mathbf{y}\in \Pi$, and $t\in (0,1)$. 
 We need to check
 \begin{equation}
\phi(  (1-t)\mathbf{x}  +t \mathbf{y} ) \le (1-t)\phi(  \mathbf{x} ) + t\phi(  \mathbf{y} ).  
  \label{convex1}
 \end{equation}
By construction, $\phi$ is invariant in the direction 
 of $\tau(\mathbf{z}_3)$. Letting $\Sigma$ be an arbitrary line orthogonal to $\tau(\mathbf{z}_3)$ and $\mathbf{x}^*$ (resp. $\mathbf{y}^*$) be the projections of $\mathbf{x}$ (resp. $\mathbf{y}$) on $\Sigma$ for the projection of
 $(1-t)\mathbf{x}  +t \mathbf{y}$ on $\Sigma$, we get 
 \[[(1-t)\mathbf{x}  +t \mathbf{y}]^* = (1-t)\mathbf{x}^*  +t \mathbf{y}^*.\] 
Since $\phi$ is invariant in the direction 
 of $\tau(\mathbf{z}_3)$, we get
 \[ \phi( (1-t)\mathbf{x}  +t \mathbf{y} ) = \phi( (1-t)\mathbf{x}^*  +t \mathbf{y}^*) ,\quad
 \phi( \mathbf{x}   )=\phi( \mathbf{x}^*   ),
 \quad \phi( \mathbf{y}   )=\phi( \mathbf{y}^*   ),\]
 and \eqref{convex1} is equivalent to
 \begin{equation}
\phi(  (1-t)\mathbf{x}^*  +t \mathbf{y}^* ) \le (1-t)\phi(  \mathbf{x}^* ) + t\phi(  \mathbf{y}^* ).  \label{convex2}
   \end{equation}
Note that \eqref{convex2} holds true due to the convexity of $\phi$ along the direction 
 $\mathbf{n}(\mathbf{z}_3)\parallel \Sigma$. 
The convexity of $\phi$ gives $[\phi]_{vv}\ge 0$ for any (nonzero) vector $\mathbf{v}$. 
Combining with the facts $0\leq M_\delta (\mathbf{x})\leq C$ as shown in Lemma \ref{Mdelta} and $\int_{-\delta}^{\delta} H(|l|) [\phi(\mathbf{x}_l)-\phi(\mathbf{x})] d \mathbf{x}_l =\dfrac{1}{2}[\phi]_{pp}+\dfrac{\kappa}{2}[\phi]_{nn}(\mathbf{x}-\overline{\mathbf{x}})\cdot \mathbf{x}$ as shown in \eqref{eqn:H}, we infer directly that
\begin{equation}\label{positive1}
M_\delta(\mathbf{x})\int_{-\delta}^{\delta} H(|l|) [\phi(\mathbf{x}_l)-\phi(\mathbf{x})] d \mathbf{x}_l\ge 0.
\end{equation}
It remains to show the bounds for
\begin{align}
   \int_{(\Omega \cup \partial\Omega_{D\delta})\cap B(\mathbf{x},\delta)}(\phi(\mathbf{y})-\phi(\mathbf{x})) d\mathbf{y}&  
 ,\qquad \forall \mathbf{x}\in\Omega\setminus \Omega_{N\delta},\label{positive3}\\
 \int_{\Omega \cap B(\mathbf{x},\delta) }(\phi(\mathbf{y})-\phi(\mathbf{x})) d\mathbf{y} &
 ,\qquad \forall \mathbf{x}\in \Omega_{N\delta}.\label{positive2}
\end{align}

{\em Step 2: bound for \eqref{positive3}.} Note that in this case $B(\mathbf{x},\delta)\subset\Omega \cup \partial\Omega_{D\delta}$. 
Let $\ell(\mathbf{x})$ be the line through $\mathbf{x}$ and parallel to $\ell:=\tau(\mathbf{z}_3)$. 
Noting that $B(\mathbf{x},\delta)$ is symmetric with respect to $\ell(\mathbf{x})$,
for any $\mathbf{y}\in B(\mathbf{x},\delta)$ we denote by $\mathbf{y}^* $ the reflection
of $\mathbf{y}$ across $\ell(\mathbf{x})$. Let $B^+(\mathbf{x},\delta)$ (resp. $B^-(\mathbf{x},\delta)$) be the ``upper'' (resp. ``lower'') half ball, 
 then 
 \begin{align*}
 &\int_{(\Omega \cup \partial\Omega_{D\delta})\cap B(\mathbf{x},\delta)} (\phi(\mathbf{y})-\phi(\mathbf{x})) d\mathbf{y}
 =\int_{ B^+(\mathbf{x},\delta)}(\phi(\mathbf{y})-\phi(\mathbf{x})) +(\phi(\mathbf{y}^*)-\phi(\mathbf{x})) d\mathbf{y}\\
 =&\int_0^\delta  [(\sqrt{\phi(\mathbf{x})} +\rho )^2  +
 (\sqrt{\phi(\mathbf{x})} -\rho )^2-2\phi(\mathbf{x}))] 2\sqrt{\delta^2-\rho^2}d\rho=4\int_0^\delta \rho^2\sqrt{\delta^2-\rho^2}d\rho \\
 \ge&4\int_0^{\delta/2} \rho^2\sqrt{\delta^2-\delta^2/4}d\rho = \frac{\delta^4}{4\sqrt{3}}.
 \end{align*}
Recalling $J_\delta=\dfrac{4}{\pi\delta^4}$ on its support, we obtained 
$-L_\delta\phi(\mathbf{x})\ge \dfrac{2}{\pi\sqrt{3}}$, $\forall \mathbf{x}\in\Omega\setminus \Omega_{N\delta}$.
 
 \begin{figure}
 \centering
 \includegraphics[scale=0.4]{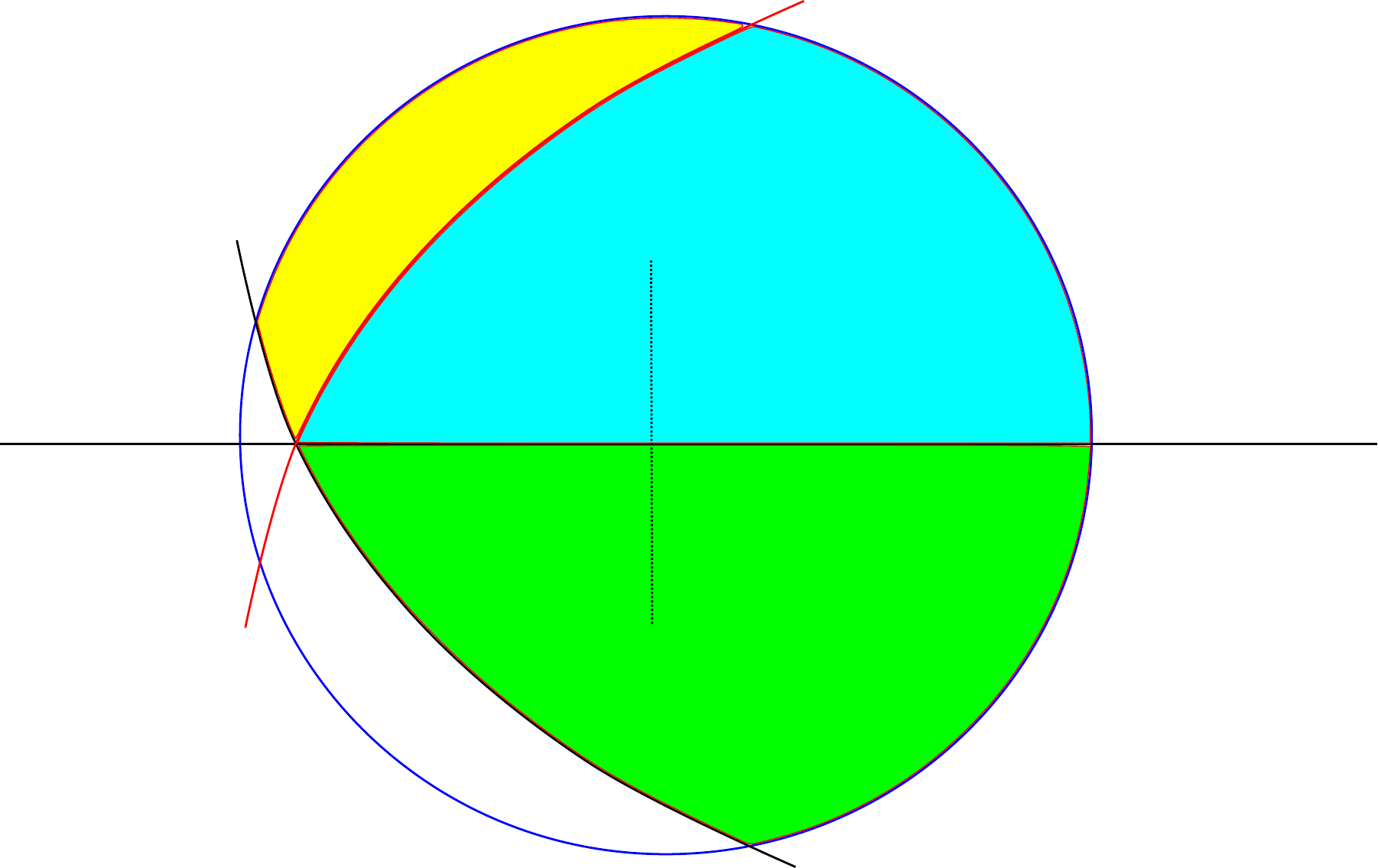}
 \put(-110,67){$\mathbf{x}$}\put(-105,92){$\mathbf{y}$}\put(-105,32){$\mathbf{y}^*$}
 \put(-10,67){$\ell(\mathbf{x})$}
 \put(-118,12){$\partial\Omega $}\put(-118,113){$[\partial\Omega]^*$}
 \caption{Notation for estimating the bound of $\int_{\Omega \cap B(\mathbf{x},\delta)}(\phi(\mathbf{y})-\phi(\mathbf{x})) d\mathbf{y}$ when $\mathbf{x}\in \Omega_{N\delta}$. Here green denotes the region of $B^-(\mathbf{x},\delta)\cap \Omega$, cyan denotes $[B^-(\mathbf{x},\delta)\cap \Omega]^*$, and yellow denotes $[\Omega \cap B^+(\mathbf{x},\delta) ]\setminus[B^-(\mathbf{x},\delta)\cap \Omega]^*$. $[\partial\Omega]^*$
 is the reflection of $\partial\Omega$ across $\ell(\mathbf{x})$.}
\label{fig:5}
\end{figure}

 \begin{figure}
 \subfigure{\includegraphics[scale=0.34]{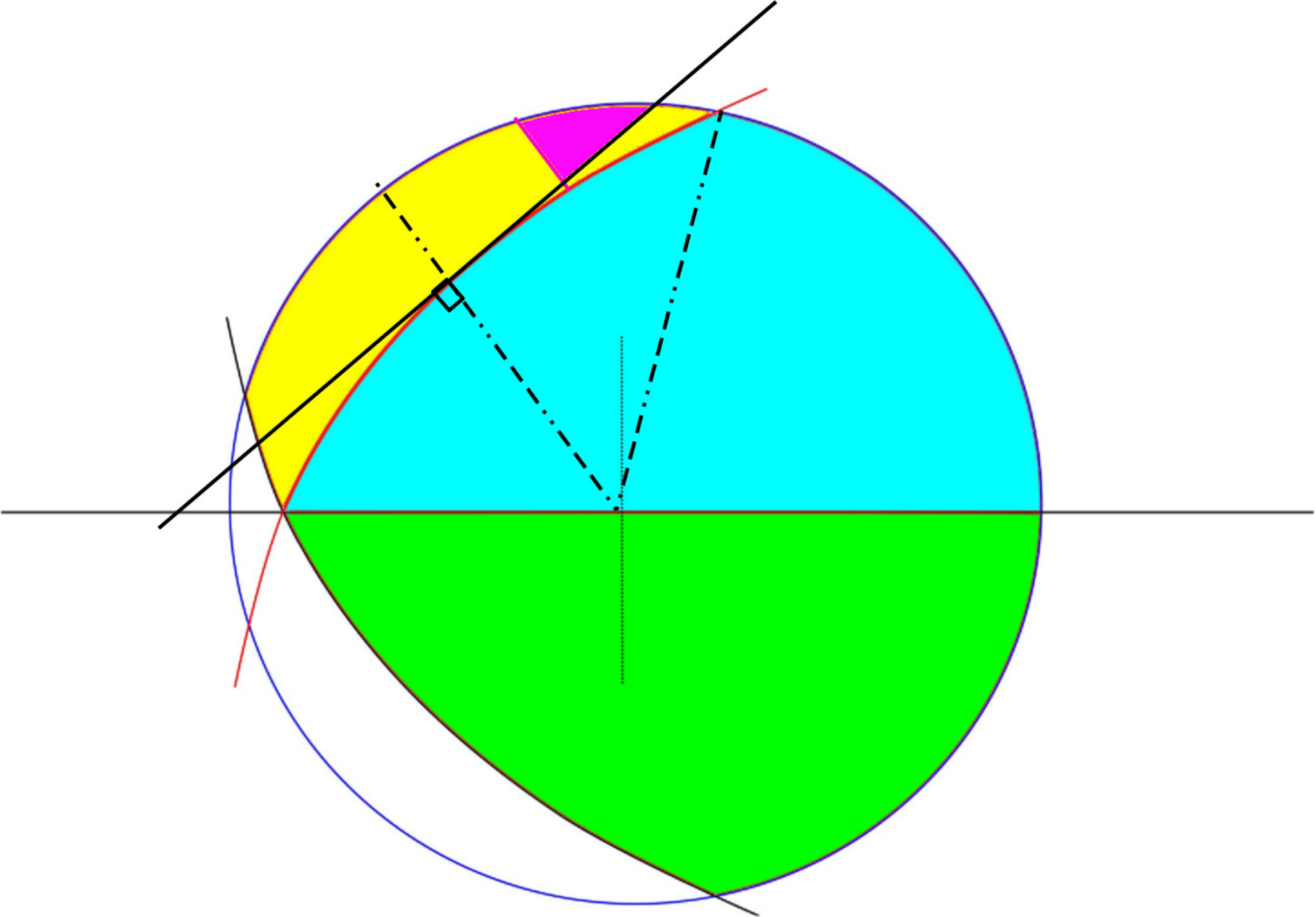}
 \put(-90,48){$\mathbf{x}$}
 \put(-25,48){$\ell(\mathbf{x})$}
 \put(-153,60){A}\put(-122,75){B}\put(-105,88){C}\put(-92,113){D}\put(-112,110){E}
 \put(-118,10){$\partial\Omega $}
 }
 \subfigure{\includegraphics[width=0.35\textwidth]{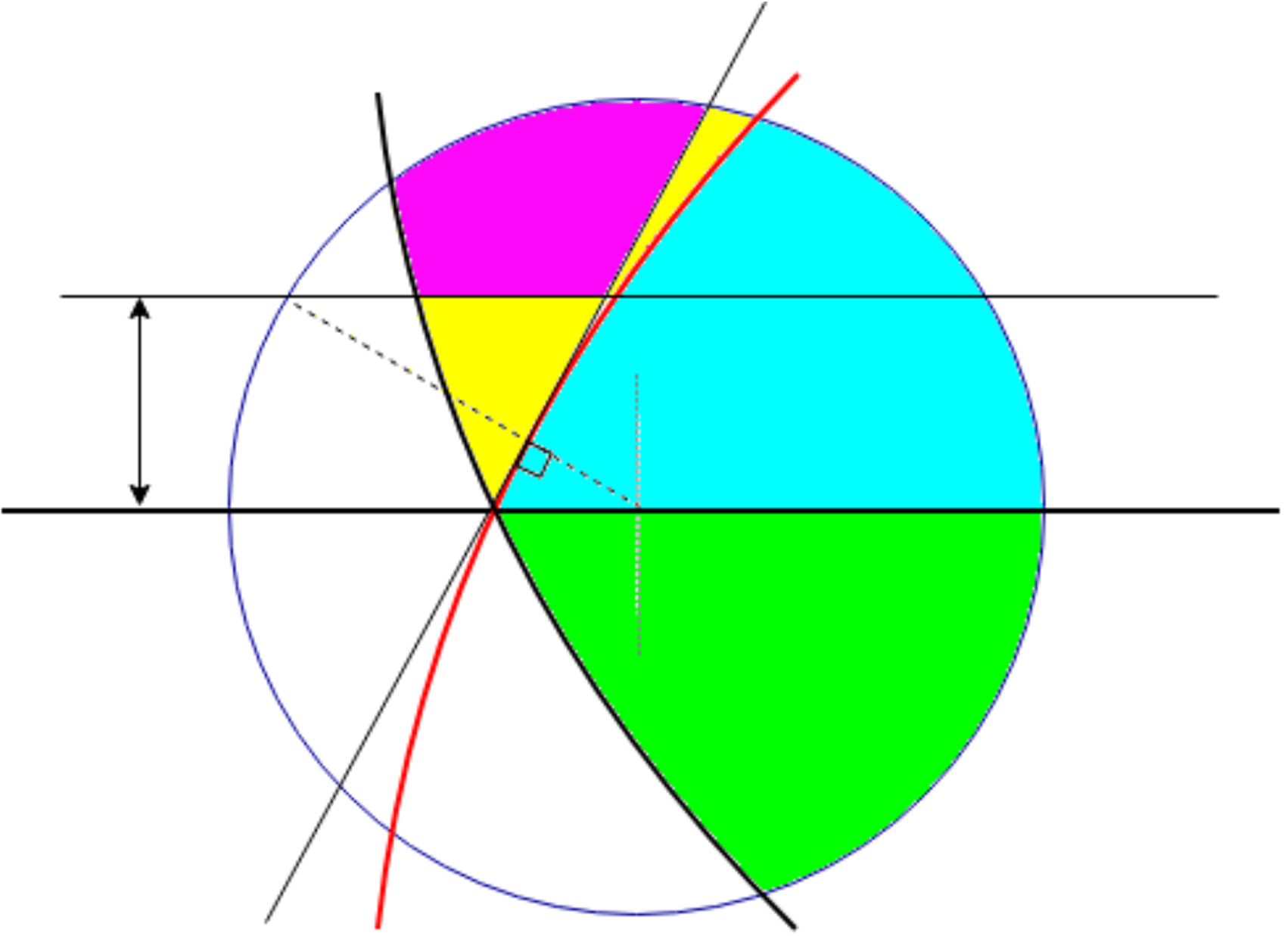}
 \put(-85,48){$\mathbf{x}$}
 \put(-25,48){$\ell(\mathbf{x})$}
 \put(-128,10){A}\put(-97,65){B}\put(-95,85){C}\put(-82,112){D}
 \put(-88,5){$\partial\Omega $}\put(-160,65){$\dfrac{\delta}{2}$}}
 \caption{Notation for estimating the bound of $\int_{\Omega \cap B(\mathbf{x},\delta)}(\phi(\mathbf{y})-\phi(\mathbf{x})) d\mathbf{y}$ when $\mathbf{x}\in \Omega_{N\delta}$, where the green and cyan regions denote $B^-(\mathbf{x},\delta)\cap \Omega$ and $[B^-(\mathbf{x},\delta)\cap \Omega]^*$, respectively. The union of yellow and purple regions represent
 $[\Omega \cap B^+(\mathbf{x},\delta) ]\setminus[B^-(\mathbf{x},\delta)\cap \Omega]^*$. Left: notation when $\text{dist}(\mathbf{x},\partial\Omega)> \delta/2$, where the purple region is chosen such that $\mathbf{p}(\mathbf{y}_B)\bot xB$, $\mathbf{n}(\mathbf{y}_B) \parallel CE$ and $|BC|=|BD|$. Right: notation when $\text{dist}(\mathbf{x},\partial\Omega)\leq \delta/2$, where the purple region is chosen such that $\mathbf{p}(\mathbf{y}_B)\bot xB$ and the distance from $C$ to $l(\mathbf{x})$ is $\delta/2$.}
\label{fig:ref1}
\end{figure}
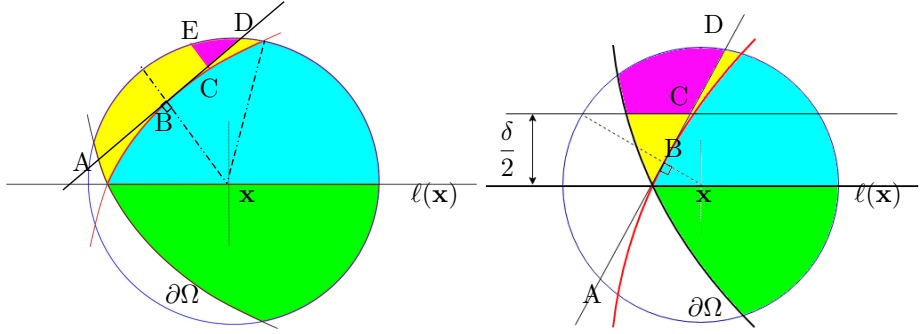


 {\em Step 3: bound of \eqref{positive2}.} For $\mathbf{x}\in\Omega_{N\delta}$, 
we will show that
{\begin{equation*}
\int_{\Omega \cap B(\mathbf{x},\delta)}(\phi(\mathbf{y})-\phi(\mathbf{x})) d\mathbf{y}\geq
 \left\{\begin{array}{cl}
  C \delta^3, & \text{ When }s_x\leq \delta/2,\\
  C \delta^4+C (\delta -s_x)^{3/2}\delta^{3/2}, & \text{ When }s_x> \delta/2.\\
 \end{array}\right.
\end{equation*}
}
 Let $[B^-(\mathbf{x},\delta)\cap \Omega]^*$ be the reflection of $B^-(\mathbf{x},\delta)\cap \Omega$
 across $\ell(\mathbf{x})$, as shown in Figure \ref{fig:5}. %
 Note the crucial geometric condition ensures that $[B^-(\mathbf{x},\delta)\cap \Omega]^*\subseteq B^+(\mathbf{x},\delta)\cap \Omega$. Since
\begin{equation}
 (\phi(\mathbf{y})-\phi(\mathbf{x})) +(\phi(\mathbf{y}^*)-\phi(\mathbf{x})) = 2\text{dist}(\mathbf{y},\ell(\mathbf{x}))^2,
 \label{def}
 \end{equation}
we have
 \begin{align*}
 \int_{\Omega \cap B(\mathbf{x},\delta)}&(\phi(\mathbf{y})-\phi(\mathbf{x})) d\mathbf{y}  
 =\int_{[\Omega \cap B^+(\mathbf{x},\delta) ]\setminus[B^-(\mathbf{x},\delta)\cap \Omega]^* }(\phi(\mathbf{y})-\phi(\mathbf{x})) d\mathbf{y}\\
 &+\int_{[B^-(\mathbf{x},\delta)\cap \Omega]^* }(\phi(\mathbf{y})-\phi(\mathbf{x})) d\mathbf{y}
 +\int_{[B^-(\mathbf{x},\delta)\cap \Omega]^* }(\phi(\mathbf{y}^*)-\phi(\mathbf{x})) d\mathbf{y}\\
 &\ge \int_{[B^-(\mathbf{x},\delta)\cap \Omega]^* }  [(\phi(\mathbf{y})-\phi(\mathbf{x})) +(\phi(\mathbf{y}^*)-\phi(\mathbf{x}))]  d\mathbf{y}\geq 0.
 \end{align*}
Since $\mathbf{x}\in \Omega$, one has $|B^-(\mathbf{x},\delta)\cap \Omega|>0$ unless 
$\mathbf{x}\in  \partial\Omega$ and $(\partial\Omega\cap B^-(\mathbf{x},\delta)) \subset(\ell(\mathbf{x})\cap B^-(\mathbf{x},\delta))$. 
Therefore, using \eqref{def}, when $s_x> \delta/2$ and $\delta\leq D/5\leq (5\sup_{\mathbf{z}\in\partial\Omega} |\kappa(\mathbf{z})|)^{-1}$, a direct computation gives
\begin{align*}
&\int_{[B^-(\mathbf{x},\delta)\cap \Omega]^* }(\phi(\mathbf{y})-\phi(\mathbf{x})) d\mathbf{y}
 +\int_{[B^-(\mathbf{x},\delta)\cap \Omega]^* }(\phi(\mathbf{y}^*)-\phi(\mathbf{x})) d\mathbf{y} \\
\ge & \int _0^{\delta/2-\delta/12} 2\rho^2 \sqrt{\delta^2-\rho^2} d\rho\ge \int _0^{5\delta/12} 2\rho^2 \sqrt{\delta^2-(\delta/2)^2} d\rho
 \ge \dfrac{125\sqrt{3}\delta^4}{5184}.
\end{align*}
On the other hand, when $s_x> \delta/2$ we calculate the integral on the purple region (denoted as $F$) shown in the left plot of Figure \ref{fig:ref1}. %
With the geometric assumption, we have %
$\phi(\mathbf{y}_B)-\phi(\mathbf{x})\geq 2\text{dist}(\mathbf{y}_B, l(\mathbf{x}))\geq 2s_x\sin(\tilde{\alpha})>\sin(\tilde{\alpha})\delta$, where $\mathbf{y}_B$ denotes the coordinate of point $B$. %
Since $\int_{[\Omega \cap B^+(\mathbf{x},\delta) ]\setminus[B^-(\mathbf{x},\delta)\cap \Omega]^* }(\phi(\mathbf{y})-\phi(\mathbf{x})) d\mathbf{y}\geq%
\int_{F}(\phi(\mathbf{y})-\phi(\mathbf{x})) d\mathbf{y}$ and $|CD|=\sqrt{\delta^2-s_x^2}/2$, $|CE|\geq (\delta-s_x)/2$ when $\delta\ll 1/\sup_{{\mathbf{z}}\in\partial\Omega}|\kappa({\mathbf{z}})|$, for $\mathbf{y}\in F$ we have
\[\phi(\mathbf{y})-\phi(\mathbf{x})\geq \phi(\mathbf{y}_B)-\phi(\mathbf{x}) \geq \sin(\tilde{\alpha})\delta,\]
and
\begin{align*}
\text{area}(F)> &\text{area}(\triangle_{CDE})=\frac{1}{8}\sqrt{|\delta-s_x|^{3}(\delta+s_x)}
\geq \frac{1}{8}\delta^{1/2}(\delta-s_x)^{3/2}.
\end{align*}
We then have
\begin{align*}
 \int_{[\Omega \cap B^+(\mathbf{x},\delta) ]\setminus[B^-(\mathbf{x},\delta)\cap \Omega]^* }&(\phi(\mathbf{y})-\phi(\mathbf{x})) d\mathbf{y}
\ge C(\delta-s_x)^{3/2}\delta^{3/2}.
\end{align*}
Similarly, for $s_x\leq \delta/2$ we have $G\subset[\Omega \cap B^+(\mathbf{x},\delta) ]\setminus[B^-(\mathbf{x},\delta)\cap \Omega]^* $ where $G$ is the purple set denoted in the right plot of Figure \ref{fig:ref1}.  For $\mathbf{y}\in G$ we have $\phi(\mathbf{y})-\phi(\mathbf{x})\geq \delta$ and 
$\text{area}(G)\geq \min\left\{\dfrac{(\sqrt{3}-1)^2\delta^2\tan(\tilde{\alpha})}{8},\dfrac{\sqrt{3}}{8}\delta^2\right\}=C\delta^2$. Therefore
\begin{align*}
 \int_{[\Omega \cap B^+(\mathbf{x},\delta) ]\setminus[B^-(\mathbf{x},\delta)\cap \Omega]^* }&(\phi(\mathbf{y})-\phi(\mathbf{x})) d\mathbf{y}
\ge C\delta^3.
\end{align*}
i.e. the contribution of a region that lies completely above $\ell(\mathbf{x})$
is of order $O(\delta^3)$, provided that it has positive area. 

Thus \eqref{positive3} and \eqref{positive2} are bounded. Combining with \eqref{positive1}, and recalling $J_\delta=A\delta^{-4}$ on its support,
we get\begin{align*}
-L_{\delta}\phi(\mathbf{x})&=2 \int_{\Omega \cup \partial\Omega_{D\delta}}
J_\delta(|\mathbf{x}-\mathbf{y}|)(\phi(\mathbf{y})-\phi(\mathbf{x})) d\mathbf{y}\ge  C
 \end{align*}
for all $\mathbf{x}\in\Omega\backslash\Omega_{N\delta}$, and
\begin{align*}
-L_{N\delta}\phi(\mathbf{x})&=2 \int_{\Omega }
J_\delta(|\mathbf{x}-\mathbf{y}|)(\phi(\mathbf{y})-\phi(\mathbf{x})) d\mathbf{y}%
+2M_\delta (\mathbf{x})\int_{-\delta}^{\delta} H(|l|) [u_0(\mathbf{x}_l)-u_0(\mathbf{x})] d \mathbf{x}_l\\
&\ge  C[\delta-s_x]^{3/2}\delta^{-5/2} +C_1>0
 \end{align*}
 for all $\mathbf{x}\in\Omega_{N\delta}$. 

Note that Lemma \ref{lem:T} and the above estimates on function $\phi$ are still insufficient to ensure second order $L^\infty(\Omega)$ convergence to the local limit, since Lemma \ref{lem:T} gives $T_\delta=O(\delta)$ on $\Omega_{N\delta}$, while the estimates for $\phi$ gives only 
\[-L_{N\delta}\phi \ge C[\delta-s_x]^{3/2}\delta^{-5/2}+C_1,\]
and it is unclear if $\dfrac{T_\delta}{-L_{N\delta}\phi}$ can be uniformly bounded from above by $C\delta^2$ as $\mathbf{x}$ approaches the inner boundary of $\Omega_{N\delta}$. The next Lemma aims to provide an estimate for $T_\delta$.


\begin{lem}
The term $T_\delta$ decays to $O(\delta^2)$ as $\mathbf{x}$ approaches the inner boundary of $\Omega_{N\delta}$, with the following bound:
{\[|T_\delta| \le C [\delta - s_x]^{3/2}\delta^{-1/2}+O(\delta^2).\]}
\end{lem}

\begin{proof}
By Lemma \ref{lem:T} and the facts $\int_{B(\mathbf{x},\delta)} J_{\delta}(|\mathbf{x}-\mathbf{y}|) ((\mathbf{x}-\mathbf{y})\cdot \mathbf{n}(\overline{\mathbf{x}}))^3d\mathbf{y}=0$, $\int_{B(\mathbf{x},\delta)} J_{\delta}(|\mathbf{x}-\mathbf{y}|) ((\mathbf{x}-\mathbf{y})\cdot \mathbf{n}(\overline{\mathbf{x}}))d\mathbf{y}=0$ and $\int_{B(\mathbf{x},\delta)} J_{\delta}(|\mathbf{x}-\mathbf{y}|) ((\mathbf{x}-\mathbf{y})\cdot \mathbf{n}(\overline{\mathbf{x}}))|(\mathbf{x}-\mathbf{y})\cdot \mathbf{p}(\overline{\mathbf{x}})|^2d\mathbf{y}=0$, we have
\begin{align*}
T_\delta 
 =&2\int_{E_\delta} J_{\delta}(|\mathbf{x}-\mathbf{y}|)\dfrac{\partial u_0({\mathbf{x}})}{\partial \mathbf{p}}%
 ((\mathbf{x}-\mathbf{y})\cdot \mathbf{p}(\overline{\mathbf{x}}))d\mathbf{y}\\
&-[u_0(\mathbf{x})]_{nnn}\left(\dfrac{1}{3}\int_{\partial\Omega_{N\delta}} J_{\delta}(|\mathbf{x}-\mathbf{y}|)%
 ((\mathbf{x}-\mathbf{y})\cdot \mathbf{n}(\overline{\mathbf{x}}))^3d\mathbf{y}%
 \right.\\
 &\left.-\int_{\partial\Omega_{N\delta}} J_{\delta}(|\mathbf{x}-\mathbf{y}|)((\mathbf{x}-{\mathbf{y}})\cdot \mathbf{n}(\overline{\mathbf{x}}))d\mathbf{y}|\overline{\mathbf{x}}-\mathbf{x}|^2\right)%
\\
&-\int_{\partial\Omega_{N\delta}} J_{\delta}(|\mathbf{x}-\mathbf{y}|)[u_0(\mathbf{x})]_{npp}%
 ((\mathbf{x}-\mathbf{y})\cdot \mathbf{n}(\overline{\mathbf{x}}))|(\mathbf{x}-\mathbf{y})\cdot \mathbf{p}(\overline{\mathbf{x}})|^2d\mathbf{y}\\
 & + \kappa(\overline{\mathbf{x}})[u_0(\mathbf{x})]_{nn}((\mathbf{x}-\overline{\mathbf{x}})\cdot \mathbf{n}(\overline{\mathbf{x}}))
 \int_{\partial\Omega_{N\delta}} J_{\delta}(|\mathbf{x}-\mathbf{y}|)\left(|(\mathbf{y}-\mathbf{x})\cdot\mathbf{p}(\overline{\mathbf{x}})|^2\right.\\
 &\left.-|(\mathbf{y}-\overline{\mathbf{x}})\cdot\mathbf{n}(\overline{\mathbf{x}})|^2+|(\mathbf{x}-\overline{\mathbf{x}})\cdot\mathbf{n}(\overline{\mathbf{x}})|^2\right)%
d\mathbf{y}+\textit{O}(\delta^2).
 \end{align*}
We firstly provide the bounds for the first term. Note that $|B({\mathbf{x}}, \delta) \cap \partial \Omega|\leq 2\sqrt{\delta^2-s_x^2}+C \kappa(\overline{\mathbf{x}})(\delta^2-s_x^2)$, therefore $|(\mathbf{x}-\mathbf{y})\cdot \mathbf{p}(\overline{\mathbf{x}})|\leq C \sqrt{\delta^2-s_x^2}$ for $\mathbf{y}\in E_\delta$. Moreover, as shown in Appendix \ref{App:T}, for the area of $E_\delta$ we have $|E_\delta|\leq C(\delta^2-s_x^2)^2+O(\delta^5)$. Then
\[\left|\int_{E_\delta} J_{\delta}(|\mathbf{x}-\mathbf{y}|)\dfrac{\partial u_0({\mathbf{x}})}{\partial \mathbf{p}}%
 ((\mathbf{x}-\mathbf{y})\cdot \mathbf{p}(\overline{\mathbf{x}}))d\mathbf{y}\right|\leq C(\delta^2-s_x^2)^{5/2}\delta^{-4}\leq C (\delta-s_x)^{3/2}\delta^{-1/2}.\]
For the rest of terms in $T_\delta$, note that the integrands
\begin{align*}
J_{\delta}(|\mathbf{x}-\mathbf{y}|)%
 ((\mathbf{x}-\mathbf{y})\cdot \mathbf{n}(\overline{\mathbf{x}}))^3 \le C\delta^{-1},\\
J_{\delta}(|\mathbf{x}-\mathbf{y}|)((\mathbf{x}-{\mathbf{y}})\cdot \mathbf{n}(\overline{\mathbf{x}})) |\overline{\mathbf{x}}-\mathbf{x}|^2
\le C\delta^{-1},\\
J_{\delta}(|\mathbf{x}-\mathbf{y}|)((\mathbf{x}-\mathbf{y})\cdot \mathbf{n}(\overline{\mathbf{x}}))|(\mathbf{x}-\mathbf{y})\cdot \mathbf{p}(\overline{\mathbf{x}})|^2
 \le C\delta^{-1},\\
J_{\delta}(|\mathbf{x}-\mathbf{y}|)|(\mathbf{y}-\mathbf{x})\cdot\mathbf{p}(\overline{\mathbf{x}})|^2 ((\mathbf{x}-\overline{\mathbf{x}})\cdot \mathbf{n}(\overline{\mathbf{x}}))\le C\delta^{-1},\\
J_{\delta}(|\mathbf{x}-\mathbf{y}|)|(\mathbf{y}-\overline{\mathbf{x}})\cdot\mathbf{n}(\overline{\mathbf{x}})|^2 ((\mathbf{x}-\overline{\mathbf{x}})\cdot \mathbf{n}(\overline{\mathbf{x}}))\le C\delta^{-1},\\
J_{\delta}(|\mathbf{x}-\mathbf{y}|)|(\mathbf{x}-\overline{\mathbf{x}})\cdot\mathbf{n}(\overline{\mathbf{x}})|^2 ((\mathbf{x}-\overline{\mathbf{x}})\cdot \mathbf{n}(\overline{\mathbf{x}}))\le C\delta^{-1}
\end{align*}
for some constant $C$. Thus it suffices to estimate the area of the domain of integration $\partial\Omega_{N\delta} \cap B(\mathbf{x},\delta)$. Since 
$|A_\delta|\leq C\delta^3$, it suffices the compute the area of $D_\delta$. 
Since $D_\delta$ is contained in the rectangle with side lengths $2\sqrt{\delta^2- s_x^2}$ and $ \delta-s_x$, direct computation then gives
\[|D_\delta| \le 2(\delta-s_x) \sqrt{\delta^2-s_x^2}\leq 2\sqrt{2}[\delta-s_x]^{3/2} \delta^{1/2}.\]
We then have $|\partial\Omega_{N\delta} \cap B(\mathbf{x},\delta)|\leq C[\delta-s_x]^{3/2} \delta^{1/2}+C\delta^3$ which together with the bounds of the integrands finishes the proof.
\end{proof}

With the above lemmas we obtain the main theorem of this section.
\begin{thm}
 Suppose $f\in C(\overline{\Omega})$, ${u}_\delta$ solves the nonlocal problem \eqref{eqn:nonlocalmix} and $u_0$ is the solution to the corresponding local problem %
 \eqref{eqn:localmix}, then for sufficiently small $\delta$ there exists a constant $C$ independent of $\delta$ such that
 \begin{equation}\label{eqn:s4bound}
  \sup_{\mathbf{x}\in\Omega}|u_\delta(\mathbf{x})-u_0(\mathbf{x})|\leq C\delta^2.
 \end{equation}
\end{thm}

{
\begin{proof}
With the barrier function $\phi$ defined as in \eqref{eqn:phi}, from the above lemmas and bounds we have
\begin{align*}
\dfrac{|T_\delta|}{-L_{\delta}\phi} \le C \delta^2, &\text{ for }\mathbf{x}\in \Omega\backslash \Omega_{N\delta},\\
\dfrac{|T_\delta|}{-L_{N\delta}\phi} \le \frac{ C_1 [\delta - s_x]^{3/2}\delta^{-1/2}
+C_2\delta^2 }{C_3[\delta - s_x]^{3/2} \delta^{-5/2}+C_4}, &\text{ for }\mathbf{x}\in \Omega_{N\delta}.
\end{align*}
Therefore, with Lemma \ref{thm:phi}, the proof of \eqref{eqn:s4bound} will be finished once we can show that $\dfrac{|T_\delta|}{-L_{N\delta}\phi}\le C \delta^2$ for $\mathbf{x}\in \Omega_{N\delta}$. Let
\[f(r):= \frac{ C_1 [\delta - r]^{3/2}\delta^{-1/2} +C_2\delta^2 }{C_3[\delta -r]^{3/2} \delta^{-5/2}+C_4},
\qquad r:= s_x \in [0,\delta].\]
Then 
\begin{align*}
f'(r)
&=\left(\dfrac{3}{2}\right)\frac{(C_2C_3-C_1C_4)[\delta - r]^{1/2}\delta^{-1/2}}{[C_3[\delta -r]^{3/2} \delta^{-5/2}+C_4]^2},
\end{align*}
and thus $f$ is monotone (either increasing or decreasing, depending on the sign of $C_2C_3-C_1C_4$).
Since
\[f(0) =\frac{ C_1 \delta  +C_2\delta^2 }{C_3 \delta^{-1}+C_4} = \frac{ C_1 \delta^2  +C_2\delta^3 }{C_3 +C_4\delta}
\le O(\delta^2),
\qquad f(\delta) = \frac{  C_2\delta^2 }{C_4},  \]
the monotonicity of $f$ ensures that $f \le O(\delta^2)$ for all $r\in [0,\delta]$, hence
we get
\[\sup_{\mathbf{x}\in \Omega_{N\delta}} \frac{|T_\delta|}{-L_{N\delta}\phi} 
\le O(\delta^2).\]
\end{proof}

}

\section{Meshfree Quadrature Rule and Numerical Solver}\label{sec:meshfree}

In this section we develop a discretization method based upon a meshfree quadrature rule for compactly supported nonlocal integro-differential equations (IDEs) with radial kernels. %
This approach is based upon the generalized moving least squares (GMLS) approximation framework \cite{Trask2018paper}, and falls within the scope of the well-established GMLS approximation theory.

We discretize the domain $\Omega$ and $\partial\Omega_{D\delta}$ by a collection of points $ \chi_{h} = \{\mathbf{x}_i\}_{\{i=1,2,\cdots,N_p\}} \subset \Omega \cup \partial \Omega_{D\delta}$, where the fill distance
\begin{equation}
	h := \sup\limits_{\mathbf{x}_i \in \chi_{h}} \min\limits_{1\leq j \leq N_p, j\neq i} |\mathbf{x}_i-\mathbf{x}_j|
\end{equation}
is a length scale characterizing the resolution of the point cloud, and $N_p$ denotes the total number of points. We define the separation distance
\begin{equation}
q_\chi = \frac12 \underset{i \neq j} \min |x_i - x_j|
\end{equation}
and assume that the point set is quasi-uniform, namely that there exists positive $c_{qu}$ satisfying
\begin{equation}
q_\chi \leq h \leq c_{qu} q_\chi.
\end{equation}

In a neighborhood of each point $\mathbf{x}_i\in \chi_h$, we reconstruct a polynomial approximation  $s_{u,\chi_{h},i}(\mathbf{x})$ to the nonlocal solution $u_\delta(\mathbf{x})$ in $B(\mathbf{x}_i,\delta)$. Specifically, we define $s_{u,\chi_{h},i}$ as the solution to the optimization problem 
\begin{equation}
s_{u,\chi_{h},i}(\mathbf{x}) = \min_{p \in \pi_m(\mathbb{R}^2)} \Bigg \{ \sum_{j=1}^{N_p}[u(\mathbf{x}_j)-p(\mathbf{x}_j)]^2w(\mathbf{x}_i,\mathbf{x}_j)\Bigg\},
\end{equation}
where $\pi_m(\mathbb{R}^2)$ are the $m$-th order polynomials in $\mathbb{R}^2$, and $w(\mathbf{x},\mathbf{y})$ is a translation-invariant positive weight function with compact support $\delta$. For concreteness we take in this work
\begin{displaymath}
w(\mathbf{x},\mathbf{y}) = \Phi_{\delta}(\mathbf{x}-\mathbf{y})=\left\{\begin{array}{cl}
(1-\frac{|\mathbf{x}-\mathbf{y}|}{\delta})^4, & \text{ when }|\mathbf{x}-\mathbf{y}|\leq \delta,\\
0& \text{ when }|\mathbf{x}-\mathbf{y}|> \delta.
\end{array}
\right.
\end{displaymath}
For a quasi-uniform pointset and sufficiently large $\delta$ the optimization problem possesses a unique solution \cite{wendland2004scattered}. We then use this polynomial reconstruction to approximate the nonlocal operator as follows.

For each point $\mathbf{x}_i$, denote the set of indices for points in $B(\mathbf{x}_i,\delta)$ as
\begin{equation}
I(\mathbf{x}_i) \equiv I(\mathbf{x}_i,\delta,\chi_{h}) := \{j\in\{1,\cdots,N_p\}:|\mathbf{x}_i-\mathbf{x}_j|<\delta\},
\end{equation}
and $\#I(\mathbf{x}_i)$ represents the number of indices in $I(\mathbf{x}_i)$. Define as a basis for $\pi_m(\mathbb{R}^2)$ the set $p_1(\mathbf{x}),p_2(\mathbf{x}),\cdots,p_Q(\mathbf{x})$, then the optimization problem has the following analytic solution.
\begin{equation}
s_{u,\chi_{h},i}(\mathbf{x}) = \tilde{u}DP(P^{\mathrm{T}}DP)^{-1}R(\mathbf{x}),
\end{equation}
where 
\begin{align*}
\tilde{u} &:= (u(\mathbf{x}_j):j\in I(\mathbf{x}_i))^{\mathrm{T}} \in \mathbb{R}^{\#I(\mathbf{x}_i)},\\
P &:= (p_k(\mathbf{x}_j))_{j\in I(\mathbf{x}_i),1\leq k \leq Q} \in \mathbb{R}^{\#I(\mathbf{x}_i) \times Q},\\
D &= diag(\Phi_\delta(\mathbf{x}_i-\mathbf{x}_j):j \in I(\mathbf{x}_i)) \in \mathbb{R}^{\#I(\mathbf{x}_i) \times 
\#I(\mathbf{x}_i)},\\
R(\mathbf{x}) &= (p_1(\mathbf{x}),\cdots,p_Q(\mathbf{x}))^{\mathrm{T}}\in \mathbb{R}^Q.
\end{align*}
This process exactly recovers $u \in \pi_m(\mathbb{R}^2)$. In the GMLS framework, the reconstruction may be used to approximate a linear bounded target functional $\varpi$ as
\begin{equation}
\varpi(u) \approx \varpi_h(u) := \varpi(s_{u,\chi_{h},i}) = \tilde{u}DP(P^{\mathrm{T}}DP)^{-1}\varpi(R(\mathbf{x})),
\end{equation}
where $\varpi(R)$ denotes the application of the target functional component-wise to each element of the polynomial basis. Classic examples of $\varpi$ include the point evaluation functional to develop meshfree approximants, point evaluations of derivatives of functions to develop meshfree collocation schemes, and integrals of functions over compact sets. In this work, we select $\varpi$ as the nonlocal operator in \eqref{eqn:nonlocaldelta} and \eqref{eqn:nonlocaleqn}, and thus obtain a meshfree estimator of the non-local operator that is exact when applied to $\pi_m(\mathbb{R}^2)$. To do so will require the computation of \eqref{eqn:nonlocaldelta} and \eqref{eqn:nonlocaleqn} applied to each member of the polynomial space.

In this paper we take $m=2$ and choose the quadratic basis functions as follows 
\begin{align*}
p_1(\mathbf{x}) &= 1, \quad
p_2(\mathbf{x}) = (\mathbf{x}-\mathbf{x}_i) \cdot \mathbf{e}_1,\quad
p_3(\mathbf{x}) = (\mathbf{x}-\mathbf{x}_i) \cdot \mathbf{e}_2, \quad
p_4(\mathbf{x}) = [(\mathbf{x}-\mathbf{x}_i) \cdot \mathbf{e}_1]^2,\\
p_5(\mathbf{x}) &= [(\mathbf{x}-\mathbf{x}_i) \cdot \mathbf{e}_2]^2,\quad
p_6(\mathbf{x}) = [(\mathbf{x}-\mathbf{x}_i) \cdot \mathbf{e}_1][(\mathbf{x}-\mathbf{x}_i) \cdot \mathbf{e_2}],
\end{align*}
where $\mathbf{e}_1 := \mathbf{n}(\overline{\mathbf{x}}_i)$, $\mathbf{e}_2: = \mathbf{p}(\overline{\mathbf{x}}_i)$ for $\mathbf{x}_i \in \Omega_{N\delta}$ and $\mathbf{e}_1 := (1,0)$, $\mathbf{e}_2 := (0,1)$ when $\mathbf{x}_i \in \Omega/\Omega_{N\delta}$. For $\mathbf{x}_i\in \Omega/\Omega_{N\delta}$, one may obtain the following formula for $\varpi_h$ in light of \eqref{eqn:nonlocaldelta}.
\begin{equation}\label{eqn:formuladis1}
-2\tilde{u}DP(P^{\mathrm{T}}DP)^{-1}\int_{B(\mathbf{x}_i,\delta)} J_{\delta}(|\mathbf{y}-\mathbf{x_i}|)(R(\mathbf{y})-R(\mathbf{x}_i))d\mathbf{y} 
= f(\mathbf{x}_i).
\end{equation}
Similarly, for $\mathbf{x}_i \in \Omega_{N\delta}$, we apply the Neumann boundary treatment and obtain the following formula for $\varpi_h$ in light of  \eqref{eqn:formula3}.
\begin{align}
 \nonumber&-2\tilde{u}DP(P^{\mathrm{T}}DP)^{-1}\int_{B(\mathbf{x}_i,\delta)\cap\Omega} J_{\delta}(|\mathbf{y}-\mathbf{x}_i|)(R(\mathbf{y})-R(\mathbf{x}_i))d\mathbf{y}\\
 \nonumber&-2 \tilde{u}DP(P^{\mathrm{T}}DP)^{-1}M_\delta(\mathbf{x}_i)\int_{-\delta}^{\delta} H_{\delta}(|l|) (R(\mathbf{x}_{l})-R(\mathbf{x}_i)) d\mathbf{x}_{l}\\
 \nonumber=&f(\mathbf{x}_i)+\left(2\int_{B(\mathbf{x}_i,\delta)\backslash\Omega} J_{\delta}(|\mathbf{y}-\mathbf{x}_i|)(\mathbf{y}-\mathbf{x}_i)\cdot\mathbf{n}(\overline{\mathbf{x}}_i)%
 d\mathbf{y}-M_\delta(\mathbf{x}_i)
 \kappa(\overline{\mathbf{x}}_i)\right)g(\overline{\mathbf{x}_i})\\
 &-\int_{B(\mathbf{x}_i,\delta)\backslash\Omega} J_{\delta}(|\mathbf{y}-\mathbf{x}_i|)%
 \left[|(\mathbf{y}-\overline{\mathbf{x}}_i)\cdot\mathbf{n}(\overline{\mathbf{x}}_i)|^2-|(\mathbf{x}_i-\overline{\mathbf{x}}_i)\cdot\mathbf{n}%
 (\overline{\mathbf{x}}_i)|^2\right]d\mathbf{y}f({{\mathbf{x}}_i}).\label{eqn:formuladis}
\end{align}
For $\mathbf{x} \in \partial \Omega_{D\delta}$, we apply the 
Dirichlet boundary condition and therefore $u_\delta(\mathbf{x})$ is given. We can then solve for $\tilde{u}$ with \eqref{eqn:formuladis1} and \eqref{eqn:formuladis}.

Numerically, the problem now reduces to how to integrate quadratic polynomials 
over $B(\mathbf{x}_i,\delta)\cap\Omega$ and $B(\mathbf{x}_i,\delta)\backslash\Omega$ properly. On simple geometries the integral in \eqref{eqn:formuladis1} and \eqref{eqn:formuladis} can be calculated analytically, while for more generalized cases where the boundary curve is more complicated, an analytic quadrature is intractable. We note that when $\delta$ is sufficiently small, $B(\mathbf{x}_i,\delta)\cap\Omega$ and $B(\mathbf{x}_i,\delta)\backslash\Omega$ can be written as the regions between two curves, and one can then evaluate the integral via numerical integration, for instance, with high-order Gaussian quadrature rules. 

\section{Numerical Results}\label{sec:test}

In this section we present the asymptotic convergence of the proposed boundary treatment by considering the nonlocal diffusion problem on three types of representative domains: a square domain in section \ref{sec:test1} which represents the case with $0$ curvature on $\partial\Omega_N$; a circular domain in section \ref{sec:test2}, which illustrates a case with constant curvature on $\partial\Omega$; and an elliptical domain in section \ref{sec:test3}, with varying curvatures along the domain boundary. Here we note that the square domain case does not satisfy the $C^3$ regularity requirement and it is therefore outside the scope of the model problem analysis presented earlier. Hence the results in section \ref{sec:test1} also demonstrate how robust the convergence rate results are when relaxing the  $C^3$ assumption on domain regularity. In this paper we focus on the type (3) convergence, i.e., the convergence of numerical solutions to the local solution as $h,\delta$ goes to $0$ simultaneously, by fixing $h/\delta=C$ and taking $h\rightarrow 0$.

\subsection{Test 1: curvature $\kappa(\mathbf{x})=0$}\label{sec:test1}

In this numerical example, we demonstrate a case where the Neumann boundary is a line segment. Specifically, we take the computational domain as $\Omega = [0,1] \times [0,1]$, with $\partial \Omega_N=\{(1,y):y\in[0,1]\}$ and $\partial \Omega_{D}=\partial\Omega\backslash \partial \Omega_N$. %
The local limit of the nonlocal problem has a smooth analytical solution 
$u_0(x,y) = \sin(\pi x)\cos(\pi y)$, together with $f(x,y) = 2\pi^2\sin(\pi x)\cos(\pi y)$ and $\frac{\partial u}{\partial \mathbf{n}}|_{x=1}=g(y)=-\pi \cos(\pi y)$. We apply the analytical local solution as a Dirichlet boundary condition over $\partial \Omega_{D\delta}$ by letting $u_{\delta} = u_0$, and impose the Neumann-type constraint \eqref{eqn:formula3} over the region $\Omega_{N\delta}=[1-\delta,1]\times[0,1]$. With uniform discretization of mesh size $h = \{1/16,1/32,1/64,1/128,/256\}$ and fixed ratio $\delta / h = 
\{4.0,3.5\}$, we demonstrate the difference between the numerical results and $u_0$ in the $L^{\infty}$-norm and $L^2$-norm in Table \ref{tab:1}. It may be seen that as $h,\delta\rightarrow 0$, the numerical solution from the proposed nonlocal Neumann-type constraint problem converges to the local analytical solution $u_0$ as $O(\delta^2)=O(h^2)$, which therefore verifies the analysis in section \ref{sec:linfty} and demonstrates the asymptotic compatibility of the numerical solver.

\begin{table}\renewcommand{\arraystretch}{0.8}
\centering
\begin{tabular}{|c|c|c|c|c|c|c|c|c|}
\hline
\multirow{2}{*}{h} & \multicolumn{4}{|c|}{$\delta /h = 4$} & \multicolumn{4}{|c|}{$\delta /h = 3.5$}\\
\cline{2-9}
&$||u_\delta-u_0||_{\infty}$ & order & $||u_\delta-u_0||_2$ & order&$||u_\delta-u_0||_{\infty}$ & order & $||u_\delta-u_0||_2$ & order\\
\hline
	$2^{-3}$ & $1.45 \times 10^{-1}$ & --  & $3.06 \times 10^{-2}$
    & -- &$9.04 \times 10^{-2}$ & -- & $3.01 \times 10^{-2}$ & -- \\
	$2^{-4}$ & $2.34 \times 10^{-2}$ & 2.63 & $5.80 \times 10^{-3}$ 
    & 2.39&$1.37 \times 10^{-2}$ & 2.72 & $5.40 \times 10^{-3}$ &
     2.48 \\
	$2^{-5}$ & $4.25 \times 10^{-3}$ & 2.38 & $1.30 \times 10^{-3}$ 
    & 2.16&$2.50 \times 10^{-3}$ & 2.45 & $1.10 \times 10^{-3}$ &
     2.30\\
	$2^{-6}$ & $1.00 \times 10^{-3}$ & 2.17 & $3.02 \times 10^{-4}$ 
    & 2.11&$5.65 \times 10^{-4}$ & 2.15 & $2.68 \times 10^{-4}$ &
     2.04\\
	$2^{-7}$ & $2.48 \times 10^{-4}$ & 2.01 & $7.38 \times 10^{-5}$ 
    & 2.03&$1.34 \times 10^{-4}$ &2.07 & $6.53 \times 10^{-5}$ &
     2.04 \\
\hline
\end{tabular}
\caption{Test 1: Convergence to the local solution for the $\kappa(\mathbf{x})=0$ case.}
\label{tab:1}
\end{table}

\subsection{Test 2: constant curvature $\kappa(\mathbf{x})$}\label{sec:test2}

We now consider as domain the unit circle $\Omega = \{(x,y)|x^2+y^2\leq 1 \}$, $\partial \Omega_N=\partial\Omega$ and with the value $u_\delta(0,-1)=u_0(0,-1)$ given to make the problem well-posed. Similarly as in test 1, we consider a smooth local solution $u_0(x,y) = \sin(\pi x)\cos(\pi y)$, with $f(x,y) = 2\pi^2\sin(\pi x)\cos(\pi y)$ and 
$\frac{\partial u}{\partial \mathbf{n}}|_{(x,y)\in \partial \Omega_N}=g(x,y)=\pi x\cos(\pi x)\cos(\pi y) - \pi y \sin(\pi x) \sin(\pi y)$,  
with uniform discretization of mesh-size $h = \{1/8,1/16,1/32,1/64,1/128\}$ and $\delta/h=\{4.0,3.5\}$. 
The $L^\infty$-norm 
and $L^2$-norm convergence results are presented in Table \ref{tab:2}. It can be observed that the convergence rate is $O(\delta^2)=O(h^2)$ as $\delta$ approaching $0$, consistent with the analysis in section \ref{sec:linfty}.

\begin{table}\renewcommand{\arraystretch}{0.8}
\centering
\begin{tabular}{|c|c|c|c|c|c|c|c|c|}
\hline
\multirow{2}{*}{h} & \multicolumn{4}{|c|}{$\delta /h = 4$} & \multicolumn{4}{|c|}{$\delta /h = 3.5$} \\
\cline{2-9}
&$||u_\delta-u_0||_{\infty}$ & order & $||u_\delta-u_0||_2$ & order&$||u_\delta-u_0||_{\infty}$ & order & $||u_\delta-u_0||_2$ & order\\
\hline
	$2^{-3}$ & $3.74 \times 10^{-1}$ & --  & $2.13 \times 10^{-1}$
    & -- &$2.98 \times 10^{-1}$ & -- & $1.77 \times 10^{-1}$ & -- \\
	$2^{-4}$ & $1.10 \times 10^{-1}$ & 1.77 & $6.88 \times 10^{-2}$ 
    & 1.63 &$8.21 \times 10^{-2}$ & 1.86 & $5.17 \times 10^{-2}$ &
     1.78\\
	$2^{-5}$ & $2.68 \times 10^{-2}$ & 2.04 & $1.68 \times 10^{-2}$ 
    & 2.03 &$1.98 \times 10^{-2}$ & 2.05 & $1.24 \times 10^{-2}$ &
     2.06 \\
	$2^{-6}$ & $6.30 \times 10^{-3}$ & 2.09 & $3.90 \times 10^{-3}$ 
    & 2.11 &$4.70 \times 10^{-3}$ & 2.07 & $2.90 \times 10^{-3}$ &
     2.10\\
	$2^{-7}$ & $1.50 \times 10^{-4}$ & 2.07 & $9.37 \times 10^{-4}$ 
    & 2.06 &$1.10 \times 10^{-3}$ &2.10 & $6.91 \times 10^{-4}$ &
     2.07 \\
\hline
\end{tabular}
\caption{Test 2: Convergence to the local solution for the $\kappa(\mathbf{x})=const$ case.}
\label{tab:2}
\end{table}

\subsection{Test 3: non-constant curvature $\kappa(\mathbf{x})$}\label{sec:test3}

In our previous two tests, the problem domains have either zero curvature or a constant curvature on the Neumann boundary. In this section we further consider a more generalized domain with a non-constant curvature on its boundary. 
We consider the ellipse $\Omega = \{(x,y)|x^2/4 + y^2 \leq 1 \}$ with $\partial \Omega_N=\partial\Omega\}$. $u_\delta(0,-1)=u_0(0,-1)$ is given to guarantee the compatibility condition. Here we note that when $\delta<1/2$, the orthogonal projection $\overline{\mathbf{x}}$ is well-defined for any $\mathbf{x}\in \Omega_{N\delta}$. We again consider a smooth local solution $u_0(x,y) = \sin(\pi x)\cos(\pi y)$ with $f(x,y) = 2\pi^2\sin(\pi x)\cos(\pi y)$, and we 
demonstrate the convergence of the numerical solution to the local solution with mesh-size $h = \{1/8,1/16,1/32,1/64,1/128\}$ and 
$\delta/h = \{4.0,3.5\}$. As shown in Table \ref{tab:3}, second order convergence is achieved which therefore verifies the estimates in section \ref{sec:linfty} and illustrates the asymptotic compatibility for a domain with nonuniform boundary curvature.

\begin{table}\renewcommand{\arraystretch}{0.8}
\centering
\begin{tabular}{|c|c|c|c|c|c|c|c|c|}
\hline
\multirow{2}{*}{h} & \multicolumn{4}{|c|}{$\delta /h = 4$}& \multicolumn{4}{|c|}{$\delta /h = 3.5$} \\
\cline{2-9}
&$||u_\delta-u_0||_{\infty}$ & order & $||u_\delta-u_0||_2$ & order&$||u_\delta-u_0||_{\infty}$ & order & $||u_\delta-u_0||_2$ & order\\
\hline
	$2^{-3}$ & $2.13 \times 10^{-1}$ & --  & $1.18 \times 10^{-1}$
    & -- &$1.73 \times 10^{-1}$ & -- & $9.60 \times 10^{-2}$ & -- \\
	$2^{-4}$ & $6.00 \times 10^{-2}$ & 1.83 & $3.32 \times 10^{-2}$ 
    & 1.78&$4.59 \times 10^{-2}$ & 1.91 & $2.53 \times 10^{-2}$ &
     1.92\\
	$2^{-5}$ & $1.43 \times 10^{-2}$ & 2.07 & $7.90 \times 10^{-3}$ 
    & 2.07&$1.08 \times 10^{-2}$ & 2.09 & $6.03 \times 10^{-3}$ &
     2.07 \\
	$2^{-6}$ & $3.40 \times 10^{-3}$ & 2.07 & $1.90 \times 10^{-3}$ 
    & 2.06&$2.60 \times 10^{-3}$ & 2.05 & $1.40 \times 10^{-3}$ &
    2.10\\
	$2^{-7}$ & $8.22 \times 10^{-4}$ & 2.05 & $4.49 \times 10^{-4}$ 
    & 2.08&$6.25 \times 10^{-4}$ &2.06 & $3.41 \times 10^{-4}$ &
     2.04 \\
\hline
\end{tabular}
\caption{Test 3: Convergence to the local solution for the non-constant $\kappa(\mathbf{x})$ case.}
\label{tab:3}
\end{table}

Moreover, we note that in the cases with constant curvature boundary, the Neumann-type constraint problem gives the analytical solution $u_\delta=u_0$ for the patch test problem with a linear solution $u_0(x,y)=x+y$. Therefore, in the previous two tests, the numerical solver passes the linear patch test with machine precision. In the elliptical domain with non-constant curvature, we further investigate the linear patch test problem, and the numerical results are illustrated in Table \ref{tab:4}. It can be observed that although the numerical solution is no longer within machine precision accuracy, the numerical solution converges to the analytical solution with an $O(h^2)$ rate as $h\rightarrow 0$.

\begin{table}\renewcommand{\arraystretch}{0.8}
\centering
\begin{tabular}{|c|c|c|c|c|c|c|c|c|}
\hline
\multirow{2}{*}{h} & \multicolumn{4}{|c|}{$\delta /h = 4$}&\multicolumn{4}{|c|}{$\delta /h = 3.5$}\\
\cline{2-9}
&$||u_\delta-u_0||_{\infty}$ & order & $||u_\delta-u_0||_2$ & order&$||u_\delta-u_0||_{\infty}$ & order & $||u_\delta-u_0||_2$ & order\\
\hline
	$2^{-3}$ & $1.71 \times 10^{-1}$ & --  & $7.87 \times 10^{-2}$
    & -- & $1.14 \times 10^{-1}$ & -- & $5.94 \times 10^{-2}$ & -- \\
	$2^{-4}$ & $2.89 \times 10^{-2}$ & 2.57 & $1.55 \times 10^{-2}$ 
    & 2.34&$2.16 \times 10^{-2}$ & 2.39 & $1.16 \times 10^{-2}$ &
     2.36\\
	$2^{-5}$ & $6.01 \times 10^{-3}$ & 2.27 & $3.20 \times 10^{-3}$ 
    & 2.28&$4.50 \times 10^{-3}$ & 2.26 & $2.40 \times 10^{-3}$ &
     2.27\\
	$2^{-6}$ & $1.20 \times 10^{-3}$ & 2.32 & $6.04 \times 10^{-4}$ 
    & 2.40&$8.35 \times 10^{-4}$ & 2.43 & $4.11 \times 10^{-4}$ &
     2.55\\
	$2^{-7}$ & $1.26 \times 10^{-4}$ & 3.25 & $4.69 \times 10^{-5}$ 
    & 3.69&$1.39\times 10^{-4}$ &2.58 & $6.20 \times 10^{-5}$ &
     2.72\\
\hline
\end{tabular}
\caption{Test 3: Linear patch test for convergence to the local solution for the non-constant $\kappa(\mathbf{x})$ case.}
\label{tab:4}
\end{table}

\section{Extension: domain with corners}\label{sec:corner}

In many popular nonlocal problem applications, it is common that the Neumann-type boundary contains corners. For example, on a peridynamic problem with damage, once a crack initiates and bifurcates, new zigzag boundary forms and the Neumann-type boundary condition must be applied on these new boundaries. To investigate how well the new Neumann-type constraint formulation extrapolates to the setting of Lipschitz domains, in this section we further extend the proposed formulation to boundaries with corners. We also numerically show the performance as well as asymptotic compatibility properties on a sample test problem with Neumann-type boundary on two sides of a square domain. Specifically, in section \ref{sec:cornerformula} we derive the formulation near a corner by approximating $-2\int_{\partial\Omega_{N\delta}} J_{\delta}(\mathbf{x}-\mathbf{y})(u(\mathbf{y})-u(\mathbf{x}))d\mathbf{y}$. Then in section \ref{sec:cornertest} we adopt a similar problem domain as in test 1 of section \ref{sec:test1} but with Neumann-type boundary conditions applied on two sides of the boundary including their intersecting corner, and demonstrate the convergence of the nonlocal solution to the corresponding local limit as $h,\delta\rightarrow 0$.

\subsection{Flux Condition and Numerical Setting}\label{sec:cornerformula}

 \begin{figure}\label{fig:corner}
 \centering
 \includegraphics[scale=0.6]{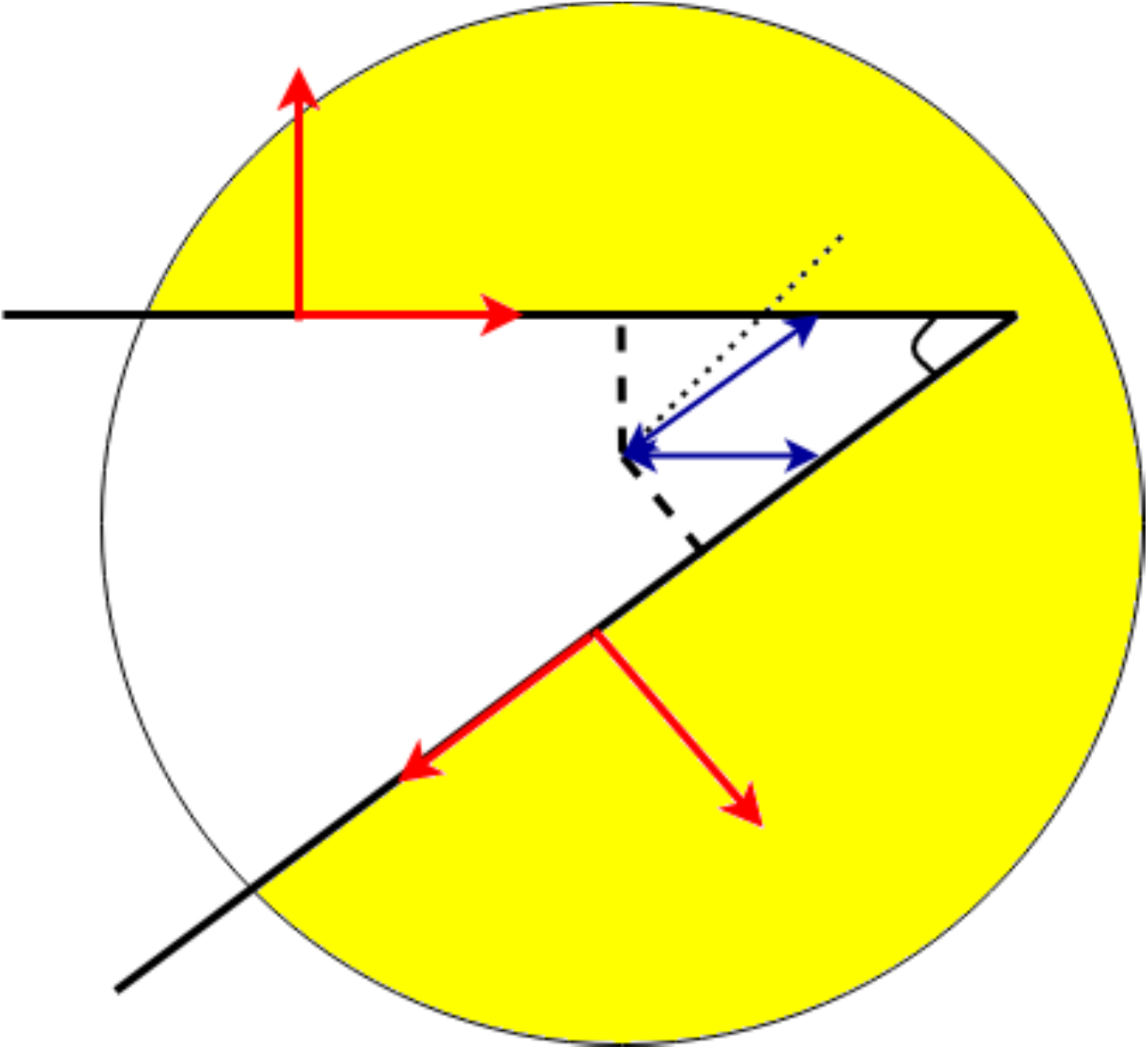}
 \put(-123,128){$\mathbf{x}$}\put(-118,163){$\overline{\mathbf{x}}_1$}\put(-97,93){$\overline{\mathbf{x}}_2$}
 \put(-100,38){$\mathbf{n}_2=(\sin\theta,-\cos\theta)$}\put(-200,215){$\mathbf{n}_1=(0,1)$}
 \put(-190,50){$\mathbf{p}_2=(-\cos\theta,-\sin\theta)$}\put(-170,143){$\mathbf{p}_1=(1,0)$}
 \put(-30,145){$\mathbf{c}=\partial\Omega_{N1}\bigcap \partial\Omega_{N2}$}
 \put(-240,162){$\partial\Omega_{N1}$}\put(-220,5){$\partial\Omega_{N2}$}
 \put(-60,140){$\theta$}\put(-63,180){$\mathbf{y}$}
 \put(-85,140){{\color{blue}$\delta_2$}}\put(-93,118){{\color{blue}$\delta_1$}}
 \caption{Geometric assumptions and notation for the corner case. Here the yellow region denotes $B(\mathbf{x},\delta)\cap \partial \Omega_{N\delta}$.}
\end{figure}
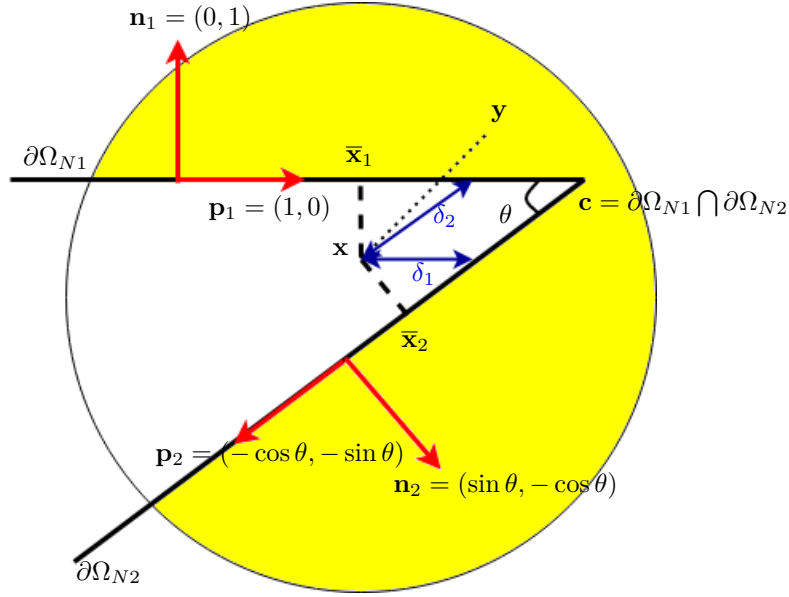

In this section we extend the numerical algorithm to a domains with corners. For simplicity, here we assume that there are two boundaries with Neumann-type boundary conditions:
\begin{align}
 \dfrac{\partial u}{\partial \mathbf{n}_1}=g_1,\quad&\text{ on }\partial\Omega_{N1},\\
 \dfrac{\partial u}{\partial \mathbf{n}_2}=g_2,\quad&\text{ on }\partial\Omega_{N2},
\end{align}
and the two boundaries intersect at $\mathbf{c}=\partial\Omega_{N1}\bigcap \partial\Omega_{N2}$. 
For any point $\mathbf{x}$ satisfying $|\mathbf{x}-\mathbf{c}|<\delta$, we project $\mathbf{x}$ onto the two boundaries respectively, i.e., %
$\mathbf{x}=\overline{\mathbf{x}}_1-s_{x1}\mathbf{n}_1(\overline{\mathbf{x}}_1)=\overline{\mathbf{x}}_2-s_{x2}\mathbf{n}_2(\overline{\mathbf{x}}_2)$. In this section, we assume that both $\partial\Omega_{N1}$ and $\partial\Omega_{N2}$ are straight lines near the corner $\mathbf{c}$, although the formulation can be further extended to more general cases. Denote $\theta$ as the angle between $\partial \Omega_{N_1}$ and $\partial \Omega_{N_2}$, without loss of generality we further denote $\mathbf{n}_1 = (0,1)$ and $\mathbf{n}_2 = (\sin \theta,-\cos \theta)$. Correspondingly, we have $\mathbf{p}_1 = (1,0)$ and $\mathbf{p}_2 = (-\cos \theta,-\sin \theta)$. We illustrate geometric assumptions and notation in Figure \ref{fig:corner}. For each point $\mathbf{x}= (x_1,x_2)$, with Taylor expansion we have the following approximation for $u(\mathbf{y})-u(\mathbf{x})$ with $\mathbf{y} = (y_1,y_2)\in B(\mathbf{x},\delta)\cap \partial \Omega_{N\delta}$:
\begin{align*}
 \nonumber u(\mathbf{y})-u(\mathbf{x})=&d_1\dfrac{\partial u(\mathbf{x})}{\partial \mathbf{n}_1}+d_2\dfrac{\partial u(\mathbf{x})}{\partial \mathbf{n}_2}+\dfrac{1}{2}d_1^2[u(\mathbf{x})]_{n_1n_1}+\dfrac{1}{2}d_2^2[u(\mathbf{x})]_{n_2n_2}
 \nonumber+d_1d_2[u(\mathbf{x})]_{n_1n_2}+O(\delta^3)\\
 =&d_1g_1(\overline{\mathbf{x}}_1)%
 +d_2g_2(\overline{\mathbf{x}}_2)+\left(\dfrac{1}{2}d_1^2-%
 (\overline{\mathbf{x}}_1-\mathbf{x})\cdot \mathbf{n}_1 d_1\right)\left(-f(\mathbf{x})-[u(\mathbf{x})]_{p_1p_1}\right)\\
 \nonumber&+\left(\dfrac{1}{2}d_2^2-%
 (\overline{\mathbf{x}}_2-\mathbf{x})\cdot \mathbf{n}_2 d_2\right)\left(-f(\mathbf{x})-[u(\mathbf{x})]_{p_2p_2}\right)\\
 \nonumber&+\dfrac{1}{2\sin \theta}d_1d_2%
 \left(\dfrac{\partial g_1(\overline{\mathbf{x}}_1)}{\partial \mathbf{p}_1}-\dfrac{\partial g_2(\overline{\mathbf{x}}_2)}{\partial \mathbf{p}_2}+f(\mathbf{x})\sin \theta \cos \theta \right)+O(\delta^3),
 \end{align*}
where
\begin{align*}
d_1 &:= \frac{\cos \theta}{\sin \theta}(y_1-x_1) + (y_2-x_2),\quad d_2 := \frac{1}{\sin \theta}(y_1-x_1).
\end{align*}
Moreover, we have
\begin{align*}
[u(\mathbf{x})]_{p_1p_1}+[u(\mathbf{x})]_{p_2p_2}=-f(\mathbf{x})+\cot\theta \dfrac{\partial g_1(\overline{\mathbf{x}}_1)}{\partial \mathbf{p}_1}-\cot\theta \dfrac{\partial g_1(\overline{\mathbf{x}}_2)}{\partial \mathbf{p}_2}+O (\delta).
\end{align*}
Let 
\[D_1 = 2\int_{\partial\Omega_{N\delta}}J_\delta(|\mathbf{x}-\mathbf{y}|)\left[\dfrac{1}{2}d_1^2-(\overline{\mathbf{x}}_1-\mathbf{x})\cdot \mathbf{n}_1 d_1\right]d \mathbf{y},\]
\[D_2 = 2\int_{\partial\Omega_{N\delta}}J_\delta(|\mathbf{x}-\mathbf{y}|)\left[\dfrac{1}{2}d_2^2-(\overline{\mathbf{x}}_2-\mathbf{x})\cdot \mathbf{n}_2 d_2\right] d\mathbf{y},\] 
substituting the above approximations into the nonlocal formulation and neglecting the higher order terms give the algorithm. For $D_1>D_2$, we take $\delta_1$ as the arc length from $\mathbf{x}$ to $\partial\Omega_{N}$ following the contour parallel to $\partial\Omega_{N1}$ and use $2\int_{-\delta_1}^{\delta_1} H_{\delta_1}(|l|)(u(\mathbf{x}_{l1})-u(\mathbf{x}))d\mathbf{x}_{l1}$ to denote the integral on this contour which approximates $[u(\mathbf{x})]_{p_1p_1}$:
\begin{align}
 \nonumber&-2\int_{\Omega} J_{\delta}(|\mathbf{x}-\mathbf{y}|)(u(\mathbf{y})-u(\mathbf{x}))d\mathbf{y}+4(D_1-D_2)\int_{-\delta_1}^{\delta_1} H_{\delta_1}(|l|)(u(\mathbf{x}_{l1})-u(\mathbf{x}))d\mathbf{x}_{l1}\\
 \nonumber=&f(\mathbf{x})-D_1 f(\mathbf{x})-D_2\cot \theta \left(\dfrac{\partial g_1(\overline{\mathbf{x}}_1)}{\partial \mathbf{p}_1}-\dfrac{\partial g_2(\overline{\mathbf{x}}_2)}{\partial \mathbf{p}_2}\right)+2\int_{\partial\Omega_{N\delta}} J_{\delta}(|\mathbf{x}-\mathbf{y}|)
 \bigg( d_1g_1(\overline{\mathbf{x}}_1)%
 \\
&+d_2g_2(\overline{\mathbf{x}}_2)+\dfrac{1}{2\sin\theta}d_1d_2%
 \left(\dfrac{\partial g_1(\overline{\mathbf{x}}_1)}{\partial \mathbf{p}_1}-\dfrac{\partial g_2(\overline{\mathbf{x}}_2)}{\partial \mathbf{p}_2}+f(\mathbf{x})\sin \theta \cos \theta \right) \bigg)%
 d\mathbf{y}.\label{eqn:cornerf1}
\end{align}
Else, we similarly take $\delta_2$ as the arc length from $\mathbf{x}$ to $\partial\Omega_{N}$ following the contour parallel to $\partial\Omega_{N2}$ and use $2\int_{-\delta_2}^{\delta_2} H_{\delta_2}(|l|)(u(\mathbf{x}_{l2})-u(\mathbf{x}))d\mathbf{x}_{l2}$ to denote the integral on this contour which approximates $[u(\mathbf{x})]_{p_2p_2}$:
\begin{align}
 \nonumber&-2\int_{\Omega} J_{\delta}(|\mathbf{x}-\mathbf{y}|)(u(\mathbf{y})-u(\mathbf{x}))d\mathbf{y}+4(D_2-D_1)\int_{-\delta_2}^{\delta_2} H_{\delta_2}(|l|)(u(\mathbf{x}_{l2})-u(\mathbf{x}))d\mathbf{x}_{l2}\\
 \nonumber=&f(\mathbf{x})-D_2f(\mathbf{x})-D_1\cot \theta \left(\dfrac{\partial g_1(\overline{\mathbf{x}}_1)}{\partial \mathbf{p}_1}-\dfrac{\partial g_2(\overline{\mathbf{x}}_2)}{\partial \mathbf{p}_2}\right)+2\int_{\partial\Omega_{N\delta}} J_{\delta}(|\mathbf{x}-\mathbf{y}|)
 \bigg( d_1g_1(\overline{\mathbf{x}}_1)%
 \\
&+d_2g_2(\overline{\mathbf{x}}_2)+\dfrac{1}{2\sin\theta}d_1d_2%
 \left(\dfrac{\partial g_1(\overline{\mathbf{x}}_1)}{\partial \mathbf{p}_1}-\dfrac{\partial g_2(\overline{\mathbf{x}}_2)}{\partial \mathbf{p}_2}+f(\mathbf{x})\sin \theta \cos \theta \right) \bigg)%
 d\mathbf{y}.\label{eqn:cornerf2}
\end{align}
Here we note that we lose coercivity in this formulation. However, numerical experiments in Section \ref{sec:cornertest} suggest that the method remains robust in practice.

\subsection{Numerical Results}\label{sec:cornertest}

In this section we investigate the numerical performance of formulation \eqref{eqn:cornerf1}-\eqref{eqn:cornerf2} on a square domain 
$\Omega = [0,1] \times [0,1]$ with Neumann-type boundary conditions applied on $\partial\Omega_{N1}=\{(1,y):y\in[0,1]\}$ and $\partial\Omega_{N2}=\{(x,1):x\in[0,1]\}$. Note that the Neumann-type boundary contains a corner $\mathbf{c}=(1,1)$ where the numerical algorithms \eqref{eqn:cornerf1}-\eqref{eqn:cornerf2} are employed. We set the analytical local solution as $u_0(x,y) = x^2y^2$, which then yields $f(x,y) = -2(x^2+y^2)$, $\frac{\partial u}{\partial \mathbf{n}}|_{x=1} = g_1(y) = 2y^2$ and $\frac{\partial u}{\partial \mathbf{n}}|_{y=1} = g_2(x) = 2x^2$. The Dirichlet-type condition $u=u_0$ is provided in a layer $\partial \Omega_{D\delta} = \{(x,y)|(x,y) \in [-\delta,1] \times [-\delta,1]/\Omega\}$. With mesh-sizes $h = \{1/16,1/32,1/64,1/128,1/256\}$ and a fixed ratio $\delta/h = \{4.0,3.5\}$, the numerical results are shown in Table \ref{tab:5}, illustrating an $O(\delta^2)=O(h^2)$ convergence rate to the local limit.

\begin{table}\renewcommand{\arraystretch}{0.8}
\centering
\begin{tabular}{|c|c|c|c|c|c|c|c|c|}
\hline
\multirow{2}{*}{h} & \multicolumn{4}{|c|}{$\delta /h = 4$}&\multicolumn{4}{|c|}{$\delta /h = 3.5$}\\
\cline{2-9}
&$||u_\delta-u_0||_{\infty}$ & order & $||u_\delta-u_0||_2$ & order&$||u_\delta-u_0||_{\infty}$ & order & $||u_\delta-u_0||_2$ & order\\
\hline
	$2^{-3}$ & $7.43 \times 10^{-2}$ & --  & $1.91 \times 10^{-2}$
    & -- &$5.45 \times 10^{-2}$ & -- & $1.46 \times 10^{-2}$ & --\\
	$2^{-4}$ & $1.52 \times 10^{-2}$ & 2.29& $4.01 \times 10^{-3}$ 
    & 2.26 &$1.13 \times 10^{-2}$ & 2.27 & $3.10 \times 10^{-3}$ &
     2.24\\
	$2^{-5}$ & $3.30 \times 10^{-3}$ & 2.20 & $9.12 \times 10^{-4}$ 
    & 2.13 &$2.40 \times 10^{-3}$ & 2.24 & $6.97 \times 10^{-4}$ &
     2.15\\
	$2^{-6}$ & $7.42 \times 10^{-3}$ & 2.15 & $2.17 \times 10^{-4}$ 
    & 2.06 &$5.60 \times 10^{-4}$ & 2.11 & $1.66 \times 10^{-4}$ &
     2.07\\
	$2^{-7}$ & $1.74 \times 10^{-4}$ & 2.09 & $5.32 \times 10^{-5}$ 
    & 2.03 &$1.31 \times 10^{-4}$ &2.09 & $4.04 \times 10^{-5}$ &
     2.03 \\
\hline
\end{tabular}
\caption{Convergence to the local solution in the case with corner.}
\label{tab:5}
\end{table}

\section{Conclusion and Future Work}\label{sec:conclusion}

In this paper we have introduced a new nonlocal Neumann-type constraint for the $2D$ nonlocal diffusion problem which is an analogue to the local flux boundary condition and for the first time achieved the optimal second-order convergence rate $O(\delta^2)$ to the local limit in the $L^{\infty}(\Omega)$ norm. The formulation is applied on a collar layer inside the domain and therefore requires no mesh or extrapolation outside the problem domain, which enables the possibility of applying the physical boundary conditions on a sharp interface. We have shown that when the problem domain is bounded, convex, connected and possesses sufficient regularity, the proposed nonlocal Neumann-type constraint with the nonlocal diffusion equation is well-posed. The nonlocal solution $u_\delta$ converges to the solution $u_0$ from the corresponding local problem in the $L^2(\Omega)$ norm as the horizon size $\delta\rightarrow 0$. Moreover, when the solution is continuous in $\overline{\Omega}$ and the Neumann type boundary is convex, we have further proved the second-order convergence of $u_\delta$ in the $L^{\infty}(\Omega)$ norm. Numerically, we have developed an asymptotically compatible particle method based on a meshfree quadrature rule for the Neumann-type constraint problem. Numerical examples on domains with representative geometries and boundary curvatures were investigated, and the optimal convergence rate $O(\delta^2)$ in the $L^{\infty}(\Omega)$ norm was observed in all instances, verifying the asymptotic compatibility of both the Neumann boundary treatment and discretization. Finally, we have demonstrated that the regularity assumption may be relaxed in practice and the formulation can be extended to domain with corners, which greatly improves the applicability of the proposed formulation for more complicated scenarios. Although the formulation does not preserve formal coercivity near the corner, numerical experiments indicate that the formulation is robust in practice and achieves the optimal convergence rate to the local limit.

We note that the formulation described in this paper actually provides an approach for applying the Neumann-type boundary condition on general compactly supported nonlocal integro-differential equations (IDEs) with radial kernels. As a natural extension, we are working on a nonlocal trace theorem which will immediately extend the current analysis results in the $L^2$ norm to problems with inhomogeneous boundary conditions, and we are also developing a sharp traction boundary condition for peridynamics which is consistent with the classical elasticity theory.

\section{Appendix}

\subsection{Proof of Lemma \ref{thm:l2conv}}\label{App:l2conv}

In this section we aim to provide the detailed proof for Lemma \ref{thm:l2conv}. Since $B_\delta(\tilde{u}_\delta,v)=(f,v)_{L^2(\Omega)}$ for any $v\in S_{\delta}$, with Lemma \ref{lemma:poincare}-\ref{thm:boundcoer} we have %
\begin{displaymath}
 ||\tilde{u}_\delta||^2_{S_{\delta}}\leq C B_\delta(\tilde{u}_\delta,\tilde{u}_\delta)=C(f,\tilde{u}_\delta)_{L^2(\Omega)}%
 \leq C||f||_{L^2(\Omega)}||\tilde{u}_\delta||_{L^2(\Omega)}\leq C ||f||_{L^2(\Omega)}||\tilde{u}_\delta||_{S_{\delta}}
\end{displaymath}
which yields the uniform boundedness of $\{\tilde{u}_\delta\}$. With Lemma \ref{thm:precompact}, we have the convergence of %
a subsequence of $\{\tilde{u}_\delta\}$ in $L^2(\Omega)$. Here we use the same $\tilde{u}_\delta$ to denote the convergent subsequence, then %
$\tilde{u}_\delta\rightarrow u_*\in S_0$. To proof the lemma, it suffices to show that $u_*=u_0$ or
\begin{equation}\label{eqn:B0equaf}
 B_0(u_*,v):=(\nabla u_*,\nabla v)=(f,v)_{L^2(\Omega)},\quad \forall v\in C^{\infty}(\Omega).
\end{equation}
Taking a standard mollifier $\phi_\epsilon$ satisfying $\int_{B(0,\epsilon)}\phi_\epsilon (\mathbf{x})d\mathbf{x}=1$ and letting $\tilde{u}_{\delta,\epsilon}=\int_{B(0,\epsilon)}\tilde{u}_{\delta}(\mathbf{x}-\mathbf{y})\phi_\epsilon(\mathbf{y})d\mathbf{y}$, we define %
$\Omega_\epsilon=\{\mathbf{x}\in\Omega:\text{dist}(\mathbf{x},\partial \Omega)<\epsilon\}$ and $\Omega^\epsilon=\{\mathbf{x}\in\Omega:\text{dist}(\mathbf{x},\partial \Omega)\geq\epsilon\}$. Assuming that $\epsilon>\delta$, for $v\in C^{\infty}(\Omega)$ %
we denote
\begin{align}
 \nonumber B_{\delta}^\epsilon(u,v)=&\int_{\Omega^\epsilon}\int_{\Omega^\epsilon}J_{\delta}(|\mathbf{x}-\mathbf{y}|)(u(\mathbf{y})-u(\mathbf{x}))%
 (v(\mathbf{y})-v(\mathbf{x}))d\mathbf{y}d\mathbf{x},\\
 \nonumber B_{0}^\epsilon(u,v)=&\int_{\Omega^\epsilon}\nabla u\cdot \nabla v d\mathbf{x}.
\end{align}
Since 
\begin{align}
 \nonumber&\nonumber B_{\delta}^\epsilon(\tilde{u}_{\delta,\epsilon},v)=%
 \int_{\Omega^\epsilon}\int_{\Omega^\epsilon}J_{\delta}(|\mathbf{x}-\mathbf{y}|)(\tilde{u}_{\delta,\epsilon}(\mathbf{y})-\tilde{u}_{\delta,\epsilon}(\mathbf{x}))%
 (v(\mathbf{y})-v(\mathbf{x}))d\mathbf{y}d\mathbf{x}\\
 \nonumber=&\int_{\Omega^\epsilon}\int_{\Omega^\epsilon}J_{\delta}(|\mathbf{x}-\mathbf{y}|)\left(\int_{B(0,\epsilon)}\phi_\epsilon (\mathbf{z})%
 \tilde{u}_{\delta}(\mathbf{y}-\mathbf{z})d\mathbf{z}-\int_{B(0,\epsilon)}\phi_\epsilon (\mathbf{z})\tilde{u}_{\delta}(\mathbf{x}-\mathbf{z})d\mathbf{z}\right)%
 (v(\mathbf{y})-v(\mathbf{x}))d\mathbf{y}d\mathbf{x}\\
 \nonumber=&\int_{B(0,\epsilon)}\phi_\epsilon (\mathbf{z})\left(\int_{\Omega^\epsilon}\int_{\Omega^\epsilon}J_{\delta}(|\mathbf{x}-\mathbf{y}|)%
 (\tilde{u}_{\delta}(\mathbf{y}-\mathbf{z})-\tilde{u}_{\delta}(\mathbf{x}-\mathbf{z}))%
 (v(\mathbf{y})-v(\mathbf{x}))d\mathbf{y}d\mathbf{x}\right)d\mathbf{z}\\
 \nonumber=&\int_{B(0,\epsilon)}\phi_\epsilon (\mathbf{z})B_{\delta}^\epsilon(\tilde{u}_{\delta}(\mathbf{x}-\mathbf{z}),v(\mathbf{x}))d\mathbf{z},
\end{align}
to show \eqref{eqn:B0equaf} it suffices to prove that when $\delta\rightarrow0$ first then $\epsilon\rightarrow0$, we have
\begin{equation}\label{eqn:Bde1}
 B_{\delta}^\epsilon(\tilde{u}_{\delta,\epsilon},v)\rightarrow %
B_{0}(u_*,v),
\end{equation}
and
\begin{equation}\label{eqn:Bde2}
 \int_{B(0,\epsilon)}\phi_\epsilon (\mathbf{z})B_{\delta}^\epsilon(\tilde{u}_{\delta}(\mathbf{x}-\mathbf{z}),v(\mathbf{x}))d\mathbf{z}\rightarrow %
(f,v)_{L^2(\Omega)}.
\end{equation}
To show \eqref{eqn:Bde1} we first fix $\epsilon$ and let $\delta\rightarrow 0$. Since $\epsilon>\delta$, $\Omega^\epsilon\bigcap\Omega_\delta=\Phi$ and %
$\Omega_\delta\subset\Omega_\epsilon$. %
Then
\begin{align}
 \nonumber B_{\delta}^\epsilon(\tilde{u}_{\delta,\epsilon},v)=&\int_{\Omega^\epsilon}\int_{\Omega}J_{\delta}(|\mathbf{x}-\mathbf{y}|)(\tilde{u}_{\delta,\epsilon}(\mathbf{y})-\tilde{u}_{\delta,\epsilon}(\mathbf{x}))%
 (v(\mathbf{y})-v(\mathbf{x}))d\mathbf{y}d\mathbf{x}\\
 \nonumber&-\int_{\Omega^\epsilon}\int_{\Omega_\epsilon}J_{\delta}(|\mathbf{x}-\mathbf{y}|)(\tilde{u}_{\delta,\epsilon}(\mathbf{y})-\tilde{u}_{\delta,\epsilon}(\mathbf{x}))%
 (v(\mathbf{y})-v(\mathbf{x}))d\mathbf{y}d\mathbf{x}.
\end{align}
Since $\tilde{u}_{\delta,\epsilon}\rightarrow {u}_{*,\epsilon}$ as $\delta\rightarrow 0$, with \cite[Proposition~3.4]{tian2014asymptotically} and %
the Dominated Convergence Theorem,
\begin{align}
&\lim_{\epsilon\rightarrow 0}\lim_{\delta\rightarrow 0}\int_{\Omega^\epsilon}\int_{\Omega}J_{\delta}(|\mathbf{x}-\mathbf{y}|)
 (\tilde{u}_{\delta,\epsilon}(\mathbf{y})-\tilde{u}_{\delta,\epsilon}(\mathbf{x}))%
 (v(\mathbf{y})-v(\mathbf{x}))d\mathbf{y}d\mathbf{x}=\lim_{\epsilon\rightarrow 0}B_0^\epsilon(u_{*,\epsilon},v)=B_0(u_{*},v).
\end{align}
On the other hand, for the second term, with the uniform boundedness
\begin{align*}
&\lim_{\epsilon\rightarrow 0}\lim_{\delta\rightarrow 0}\left|\int_{\Omega^\epsilon}\int_{\Omega_\epsilon}J_{\delta}(|\mathbf{x}-\mathbf{y}|)%
 (\tilde{u}_{\delta,\epsilon}(\mathbf{y})-\tilde{u}_{\delta,\epsilon}(\mathbf{x}))%
 (v(\mathbf{y})-v(\mathbf{x}))d\mathbf{y}d\mathbf{x}\right|\leq C \lim_{\epsilon\rightarrow 0}\text{area}(\Omega_{\epsilon})=0.
\end{align*}
Hence \eqref{eqn:Bde1} has been proved. For \eqref{eqn:Bde2} it suffices to show that
\begin{align}
 \nonumber&\lim_{\epsilon\rightarrow 0}\lim_{\delta\rightarrow 0}%
 |B_{\delta}^\epsilon(\tilde{u}_{\delta}(\mathbf{x}-\mathbf{z}),v(\mathbf{x}))-(f(\mathbf{x}),v(\mathbf{x}))_{L^2(\Omega)}|\\
 =&\lim_{\epsilon\rightarrow 0}\lim_{\delta\rightarrow 0}%
 |B_{\delta}^\epsilon(\tilde{u}_{\delta}(\mathbf{x}-\mathbf{z}),v(\mathbf{x}))-B_{\delta}(\tilde{u}_{\delta}(\mathbf{x}),v(\mathbf{x}))|=0.
\end{align}
Denote $\Omega^{\mathbf{z}\epsilon}=\{\mathbf{x}\in\Omega:\mathbf{x}-\mathbf{z}\in\Omega^\epsilon\}$, we have
{{\begin{align}
 \nonumber&|B_{\delta}^\epsilon(\tilde{u}_{\delta}(\mathbf{x}-\mathbf{z}),v(\mathbf{x}))-B_{\delta}(\tilde{u}_{\delta}(\mathbf{x}),v(\mathbf{x}))|\\
 \nonumber=&\left|\int_{\Omega^\epsilon}\int_{\Omega^\epsilon}J_{\delta}(|\mathbf{x}-\mathbf{y}|)%
 (\tilde{u}_{\delta}(\mathbf{y}-\mathbf{z})-\tilde{u}_{\delta}(\mathbf{x}-\mathbf{z}))%
 (v(\mathbf{y})-v(\mathbf{x}))d\mathbf{y}d\mathbf{x}\right.\\
 \nonumber &-\int_{\Omega}\int_{\Omega}J_{\delta}(|\mathbf{x}-\mathbf{y}|)%
 (\tilde{u}_{\delta}(\mathbf{y})-\tilde{u}_{\delta}(\mathbf{x}))%
 (v(\mathbf{y})-v(\mathbf{x}))d\mathbf{y}d\mathbf{x}\\
 \nonumber &-\int_{\Omega_\delta} M_\delta(\mathbf{x}) \int_{-\delta}^{\delta} H_{\delta}(|l|) [\tilde{u}^{\delta}(\mathbf{x}_{l})-\tilde{u}^{\delta}(\mathbf{x})]%
 [v(\mathbf{x}_{l})-v(\mathbf{x})] d\mathbf{x}_{l}d\mathbf{x}\\
 \nonumber &\left.-\int_{\Omega_\delta} \int_{-\delta}^{\delta}\left [M_\delta(\mathbf{x}_{l})\dfrac{|\mathbf{r}'(\mathbf{x})|}{|\mathbf{r}'(\mathbf{x}_l)|}-M_\delta(\mathbf{x})\right] H_{\delta}(|l|) [\tilde{u}^{\delta}(\mathbf{x}_{l})-\tilde{u}^{\delta}(\mathbf{x})] d\mathbf{x}_l v(\mathbf{x}) d\mathbf{x}\right|\\
 \nonumber\leq&\left|\int_{\Omega^{\mathbf{z}\epsilon}}\int_{\Omega^{\mathbf{z}\epsilon}}J_{\delta}(|\mathbf{x}-\mathbf{y}|)%
 (\tilde{u}_{\delta}(\mathbf{y})-\tilde{u}_{\delta}(\mathbf{x}))%
 (v(\mathbf{y}+\mathbf{z})-v(\mathbf{x}+\mathbf{z})-v(\mathbf{y})+v(\mathbf{x}))d\mathbf{y}d\mathbf{x}\right|\\
 \nonumber&+\left|\int_{\Omega^{\mathbf{z}\epsilon}}\int_{\Omega^{\mathbf{z}\epsilon}}J_{\delta}(|\mathbf{x}-\mathbf{y}|)%
 (\tilde{u}_{\delta}(\mathbf{y})-\tilde{u}_{\delta}(\mathbf{x}))%
 (v(\mathbf{y})-v(\mathbf{x}))d\mathbf{y}d\mathbf{x}\right.\\
 \nonumber&\left.-\int_{\Omega}\int_{\Omega}J_{\delta}(|\mathbf{x}-\mathbf{y}|)%
 (\tilde{u}_{\delta}(\mathbf{y})-\tilde{u}_{\delta}(\mathbf{x}))%
 (v(\mathbf{y})-v(\mathbf{x}))d\mathbf{y}d\mathbf{x}\right|\\
 \nonumber&+\left|\int_{\Omega_\delta} M_\delta(\mathbf{x}) \int_{-\delta}^{\delta} H_{\delta}(|l|) [\tilde{u}^{\delta}(\mathbf{x}_{l})-\tilde{u}^{\delta}(\mathbf{x})]%
 [v(\mathbf{x}_{l})-v(\mathbf{x})] d\mathbf{x}_{l}d\mathbf{x}\right.\\
\nonumber&\left. +\int_{\Omega_\delta} \int_{-\delta}^{\delta} \left[M_\delta(\mathbf{x}_{l})\dfrac{|\mathbf{r}'(\mathbf{x})|}{|\mathbf{r}'(\mathbf{x}_l)|}-M_\delta(\mathbf{x})\right] H_{\delta}(|l|) [\tilde{u}^{\delta}(\mathbf{x}_{l})-\tilde{u}^{\delta}(\mathbf{x})] d\mathbf{x}_l v(\mathbf{x}) d\mathbf{x}\right|\\
 \nonumber:=&I+II+III.
\end{align}}}
For the first term we have
\begin{align}
 \nonumber I
 \nonumber\leq&||\tilde{u}_{\delta}||_{S_\delta}||v(\mathbf{x}+\mathbf{z})-v(\mathbf{x})||_{S_\delta}%
 \leq ||\tilde{u}_{\delta}||_{S_\delta}||v(\mathbf{x}+\mathbf{z})-v(\mathbf{x})||_{S_0}
\end{align}
which goes to $0$ as $\epsilon\rightarrow0$ since $|\mathbf{z}|\leq\epsilon$ and $v\in C^{\infty}(\Omega)$. For the second term, since $|\mathbf{z}|<\epsilon$, $\Omega\backslash\Omega^{\mathbf{z}\epsilon}\subset\Omega_{2\epsilon}$. Therefore
{\begin{align}
 \nonumber II
 \nonumber\leq&\left|\int_{\Omega\backslash\Omega^{\mathbf{z}\epsilon}}\int_{\Omega}J_{\delta}(|\mathbf{x}-\mathbf{y}|)%
 (\tilde{u}_{\delta}(\mathbf{y})-\tilde{u}_{\delta}(\mathbf{x}))%
 (v(\mathbf{y})-v(\mathbf{x}))d\mathbf{y}d\mathbf{x}\right|\\
 \nonumber & +\left|\int_{\Omega\backslash\Omega^{\mathbf{z}\epsilon}}\int_{\Omega^{\mathbf{z}\epsilon}}J_{\delta}(|\mathbf{x}-\mathbf{y}|)%
 (\tilde{u}_{\delta}(\mathbf{y})-\tilde{u}_{\delta}(\mathbf{x}))%
 (v(\mathbf{y})-v(\mathbf{x}))d\mathbf{y}d\mathbf{x}\right|\\
 \nonumber\leq&2||\tilde{u}_{\delta}||_{S_\delta}\left(\int_{\Omega\backslash\Omega^{\mathbf{z}\epsilon}}\int_{\Omega}J_{\delta}(|\mathbf{x}-\mathbf{y}|)%
 (v(\mathbf{y})-v(\mathbf{x}))^2d\mathbf{y}d\mathbf{x}\right)^{1/2}\\
 \nonumber\leq&2||\tilde{u}_{\delta}||_{S_\delta}\left(\int_{\Omega_{2\epsilon}}\int_{\Omega}J_{\delta}(|\mathbf{x}-\mathbf{y}|)%
 (v(\mathbf{y})-v(\mathbf{x}))^2d\mathbf{y}d\mathbf{x}\right)^{1/2}.
\end{align}}
Since $v\in C^\infty(\Omega)$, we have $\lim_{\epsilon\rightarrow0}\lim_{\delta\rightarrow 0}\int_{\Omega_{2\epsilon}}\int_{\Omega}J_{\delta}(\mathbf{x}-\mathbf{y})%
 (v(\mathbf{y})-v(\mathbf{x}))^2d\mathbf{y}d\mathbf{x}=0$ and therefore $\lim_{\epsilon\rightarrow0}\lim_{\delta\rightarrow 0}II=0$. {For the third term {we first consider the curvature$\equiv 0$ case. When $\delta$ is sufficiently small, since $M_\delta(\mathbf{x})\leq 3\pi\sup_{r\leq1} J(|r|)$ we have}
{{\begin{align}
 \nonumber III\leq&\left|\int_{\Omega_\delta} M_\delta(\mathbf{x}) \int_{-\delta}^{\delta} H_{\delta}(|l|) [\tilde{u}_{\delta}(\mathbf{x}_{l})-\tilde{u}_{\delta}(\mathbf{x})]%
 [v(\mathbf{x}_{l})-v(\mathbf{x})] d\mathbf{x}_{l}d\mathbf{x}\right|\\
 \leq&C||\tilde{u}^{\delta}||_{S_\delta}\left(\int_{\Omega_\delta} M_\delta(\mathbf{x}) \int_{-\delta}^{\delta} H_{\delta}(|l|) l^2 d\mathbf{x}_{l}%
 \sup_{\mathbf{z}\in\Omega_\delta}\left|\dfrac{\partial v(\mathbf{z})}{\partial \mathbf{p}}\right|^2 d\mathbf{x}\right)^{1/2}
 \leq C||\tilde{u}^{\delta}||_{S_\delta}\sup_{\mathbf{z}\in\Omega_\delta}\left|\dfrac{\partial v(\mathbf{z})}{\partial \mathbf{p}}\right|%
 \left(\text{area}(\Omega_\delta)\right)^{1/2}.\label{eqn:term4}
\end{align}}}
Since $v\in C^\infty(\Omega)$, $\sup_{\mathbf{z}\in\Omega_\delta}\left|\dfrac{\partial v(\mathbf{z})}{\partial \mathbf{p}}\right|\leq\infty$. Since $\Omega$ is bounded, %
$\lim_{\delta\rightarrow 0}\text{area}(\Omega_\delta)=0$. %
Hence $\lim_{\epsilon\rightarrow 0}\lim_{\delta\rightarrow 0}III=0$. To prove the case of nonzero curvature, when $\delta$ is sufficiently small \eqref{eqn:term3} and \eqref{eqn:term4} yield
{\small\begin{align*}
  III &\le \left|\int_{\Omega_\delta} M_\delta(\mathbf{x}) \int_{-\delta}^{\delta} H_{\delta}(|l|) [\tilde{u}_{\delta}(\mathbf{x}_{l})-\tilde{u}_{\delta}(\mathbf{x})]%
 [v(\mathbf{x}_{l})-v(\mathbf{x})] d\mathbf{x}_{l}d\mathbf{x}\right|\\
 &+\left|\int_{\Omega_\delta} \int_{-\delta}^{\delta} \left[M_\delta(\mathbf{x}_{l})\dfrac{|\mathbf{r}'(\mathbf{x})|}{|\mathbf{r}'(\mathbf{x})_l|}-M_\delta(\mathbf{x})\right] H_{\delta}(|l|) [\tilde{u}_{\delta}(\mathbf{x}_{l})-\tilde{u}_{\delta}(\mathbf{x})]%
 d\mathbf{x}_l v(\mathbf{x}) d\mathbf{x}\right|\\
&\le C ||\tilde{u}^{\delta}||_{S_\delta}\sqrt{\int_{\Omega_\delta} v^2(\mathbf{x})
d\mathbf{x}+\sup_{\mathbf{z}\in \Omega_\delta}\left|\dfrac{\partial v(\mathbf{z})}{\partial \mathbf{p}}\right|^2%
 \text{area}(\Omega_\delta)}\overset{\delta\to 0}\to 0.
 \end{align*}}}
Due to $v\in C^\infty(\Omega)$, as $\delta\to0$ $\int_{\Omega_\delta} v^2(\mathbf{x})
d\mathbf{x} \to 0$. Moreover, $\lim_{\delta\rightarrow 0}\text{area}(\Omega_\delta)=0$. Therefore $\lim_{\epsilon\rightarrow 0}\lim_{\delta\rightarrow 0}III=0$ and %
we have then finished the proof.

\subsection{Proof of Lemma \ref{lem:T}}\label{App:T}

In this section we aim to provide the detailed derivation for Lemma \ref{lem:T}. For $\mathbf{x}\in \Omega\backslash\Omega_{N\delta}$,
\begin{align}
\nonumber T_\delta =&L_0u_0-L_\delta u_0=-\triangle u_0+2\int_{\Omega\cup\partial\Omega_{D\delta}} %
 J_{\delta}(|\mathbf{x}-\mathbf{y}|)(u_0(\mathbf{y})-u_0(\mathbf{x}))d\mathbf{y}\\
 \nonumber=&-\triangle u_0+\int_{\Omega\cup\partial\Omega_{D\delta}} %
 J_{\delta}(|\mathbf{x}-\mathbf{y}|)[u_0(\mathbf{x})]_{pp}((\mathbf{x}-\mathbf{y})\cdot \mathbf{p}(\overline{\mathbf{x}}))^2d\mathbf{y}\\
 \nonumber&+\int_{\Omega\cup\partial\Omega_{D\delta}} %
 J_{\delta}(|\mathbf{x}-\mathbf{y}|)[u_0(\mathbf{x})]_{nn}((\mathbf{x}-\mathbf{y})\cdot \mathbf{n}(\overline{\mathbf{x}}))^2d\mathbf{y}%
 +\textit{O}(\delta^2)=\textit{O}(\delta^2).
\end{align}
 
For $\mathbf{x}\in \Omega_{N\delta}$, we will first estimate %
$\int_{-\delta}^{\delta} H(|l|) [u_0(\mathbf{x}_l)-u_0(\mathbf{x})] d \mathbf{x}_l$. With Taylor's expansion we have
{\begin{align*}
 u_0(\mathbf{x}_{l})
=&u_0(\mathbf{x}) + \frac{\partial u_0(\mathbf{x})}{\partial \mathbf{n}} ( \mathbf{x}_{l } -\mathbf{x} ) \cdot \mathbf{n}(\mathbf{\overline{x}})%
+\frac{\partial u_0(\mathbf{x})}{\partial \mathbf{p}} ( \mathbf{x}_{l } -\mathbf{x} ) \cdot \mathbf{p}(\mathbf{\overline{x}})\\
&+\frac12 [u_0(\mathbf{x})]_{nn} |( \mathbf{x}_{l } -\mathbf{x} ) \cdot \mathbf{n}(\mathbf{\overline{x}})|^2+\frac12 [u_0(\mathbf{x})]_{pp} |( \mathbf{x}_{l } -\mathbf{x} ) \cdot \mathbf{p}(\mathbf{\overline{x}})|^2\\
&+[u_0(\mathbf{x})]_{pn} (( \mathbf{x}_{l } -\mathbf{x} ) \cdot \mathbf{n}(\mathbf{\overline{x}}))%
(( \mathbf{x}_{l } -\mathbf{x} ) \cdot \mathbf{p}(\mathbf{\overline{x}}))+\frac{1}{6} [u_0(\mathbf{x})]_{nnn} (( \mathbf{x}_{l } -\mathbf{x} ) \cdot \mathbf{n}(\mathbf{\overline{x}}))^3\\
&
+\frac{1}{6} [u_0(\mathbf{x})]_{ppp} (( \mathbf{x}_{l } -\mathbf{x} ) \cdot \mathbf{p}(\mathbf{\overline{x}}))^3%
+\frac12 [u_0(\mathbf{x})]_{pnn} |( \mathbf{x}_{l } -\mathbf{x} ) \cdot \mathbf{n}(\mathbf{\overline{x}})|^2%
(( \mathbf{x}_{l } -\mathbf{x} ) \cdot \mathbf{p}(\mathbf{\overline{x}}))\\
&+\frac12 [u_0(\mathbf{x})]_{ppn} (( \mathbf{x}_{l } -\mathbf{x} ) \cdot \mathbf{n}(\mathbf{\overline{x}}))%
|( \mathbf{x}_{l } -\mathbf{x} ) \cdot \mathbf{p}(\mathbf{\overline{x}})|^2+O(l^4).
\end{align*}}
Assuming the boundary $\partial\Omega$ is $C^3$ regular, we can approximate 
$\partial\Omega\cap B(\mathbf{x},\delta)$ with the osculating circle  $\mathcal{C}(\overline{\mathbf{x}})$. 
When $\partial\Omega$ does not coincide with $\mathcal{C}(\overline{\mathbf{x}})$, we denote $\mathbf{Px}_l$ as the point with distance $l$ to $\mathbf{x}$ along $\mathcal{C}(\overline{\mathbf{x}})$ following the $\mathbf{p}$ direction. For point $\mathbf{x}$, take the Cartesian coordinate system as shown in the right plot of Figure \ref{graph} and let $(c_1(s),c_2(s))$ be the curve of boundary $\partial\Omega$ which is parameterized by the arclength $s$. Then we have $\mathbf{x}_l=(c_1(l),c_2(l))^T$, and
\begin{displaymath}
\mathbf{x}_l=\mathbf{x}+\left(\begin{array}{c}
l\\
0\\
\end{array}\right)+\left(\begin{array}{c}
0\\
\frac{\kappa(\mathbf{x})l^2}{2}\\
\end{array}\right)+\left(\begin{array}{c}
c_1'''(0)\frac{l^3}{6} \\
c_2'''(0)\frac{l^3}{6}\\
\end{array}\right)+O(l^4),
\end{displaymath}
while $\mathbf{Px}_l=\mathbf{x}+\left(\dfrac{1}{\kappa(\mathbf{x})}\sin(l\kappa(\mathbf{x})),\dfrac{1}{\kappa(\mathbf{x})}(1-\cos(l\kappa(\mathbf{x})))\right)^T$. Therefore
\begin{displaymath}
\mathbf{x}_l-\mathbf{Px}_l=\left(\begin{array}{c}
\frac{c_1'''(0)+\kappa^2(\mathbf{x})}{6}l^3\\
\frac{c_2'''(0)}{6}l^3\\
\end{array}\right)+O(l^4).
\end{displaymath}
With $E_\delta$ to denote the region in $A_\delta$ which is asymmetric with respect to the $y$ axis in the right plot of Figure \ref{graph}, we then have the area of $E_\delta$ as $|E_\delta|\leq C(\delta^2-s_x^2)^2+O(\delta^5)$. 
Moreover, adopting the coordinates as shown in the right plot of Figure \ref{graph}, we have
$( \mathbf{x}_{l } -\mathbf{x} ) \cdot \mathbf{n}(\mathbf{\overline{x}}) = -\frac\kappa2 l^2-\frac{c_1'''(0)}{6}l^3+O(l^4)$, %
$( \mathbf{x}_{l } -\mathbf{x} ) \cdot \mathbf{p}(\mathbf{\overline{x}}) = l+\frac{c_2'''(0)}{6}l^3+O(l^4)$. {Therefore
{\begin{align*}
 u_0(\mathbf{x}_{l})-u_0(\mathbf{x})
=& -\frac{\partial u_0(\mathbf{x})}{\partial \mathbf{n}} \left(\frac\kappa2 l^2+{\frac{c_1'''(0)}{6}l^3}\right)%
+\frac{\partial u_0(\mathbf{x})}{\partial \mathbf{p}} \left(l+{\frac{c_2'''(0)}{6}l^3}\right)%
+\frac{l^2}{2} [u_0(\mathbf{x})]_{pp}\\
&-\dfrac{\kappa l^3}{2}[u_0(\mathbf{x})]_{pn}+\frac{l^3}{6} [u_0(\mathbf{x})]_{ppp}+O(l^4),
\end{align*}}
{\begin{align*}
 u_0(\mathbf{x}_{-l})-u_0(\mathbf{x})
=& -\frac{\partial u_0(\mathbf{x})}{\partial \mathbf{n}} \left(\frac\kappa2 l^2-{\frac{c_1'''(0)}{6}l^3}\right)%
+\frac{\partial u_0(\mathbf{x})}{\partial \mathbf{p}} \left(-l-{\frac{c_2'''(0)}{6}l^3}\right)+\frac{l^2}{2} [u_0(\mathbf{x})]_{pp}\\
&+\dfrac{\kappa l^3}{2}[u_0(\mathbf{x})]_{pn}-\frac{l^3}{6} [u_0(\mathbf{x})]_{ppp}+O(l^4),
\end{align*}}
which yield
{\begin{align*}
 &u_0(\mathbf{x}_{l})+u_0(\mathbf{x}_{-l})-2u_0(\mathbf{x})
= \kappa l^2\frac{\partial u_0(\mathbf{x})}{\partial \mathbf{n}}+l^2[u_0(\mathbf{x})]_{pp}+O(l^4)\\
=&\kappa l^2\frac{\partial u_0(\overline{\mathbf{x}})}{\partial \mathbf{n}}+{\kappa l^2[u_0(\mathbf{x})]_{nn} ((\mathbf{x} - \overline{\mathbf{x}}) \cdot \mathbf{n}(\mathbf{\overline{x}}))}%
+l^2[u_0(\mathbf{x})]_{pp}+O(l^4)\\
=&\kappa l^2[u_0(\mathbf{x})]_{nn} ((\mathbf{x} - \overline{\mathbf{x}}) \cdot \mathbf{n}(\mathbf{\overline{x}}))+l^2[u_0(\mathbf{x})]_{pp}+O(l^4).
\end{align*}}
Therefore
{\begin{align}
\nonumber &\int_{-\delta}^{\delta} H(|l|) [u_0(\mathbf{x}_l)-u_0(\mathbf{x})] d \mathbf{x}_l\\
\nonumber=&\int_{0}^{\delta} H(|l|) \left[\kappa l^2[u_0(\mathbf{x})]_{nn} ((\mathbf{x} - \overline{\mathbf{x}}) \cdot \mathbf{n}(\mathbf{\overline{x}}))+l^2[u_0(\mathbf{x})]_{pp}
+O(l^4)\right] d \mathbf{x}_l\\
=&\dfrac{\kappa}{2} [u_0(\mathbf{x})]_{nn} %
((\mathbf{x} - \overline{\mathbf{x}}) \cdot \mathbf{n}(\mathbf{\overline{x}}))+\dfrac{1}{2}[u_0(\mathbf{x})]_{pp}+O(\delta^2),\label{eqn:H}
\end{align}}
and
\begin{align*}
 &2M_\delta(\mathbf{x}) \int_{-\delta}^{\delta} H(|l|) [u_0(\mathbf{x}_l)-u_0(\mathbf{x})] d \mathbf{x}_l\\
 =& M_\delta(\mathbf{x})[u_0(\mathbf{x})]_{pp}+\kappa M_\delta(\mathbf{x}) [u_0(\mathbf{x})]_{nn} ((\mathbf{x} - \overline{\mathbf{x}}) \cdot \mathbf{n}(\mathbf{\overline{x}}))+ O(\delta^2)\\
  =&\kappa M_\delta(\mathbf{x}) [u_0(\mathbf{x})]_{nn} ((\mathbf{x} - \overline{\mathbf{x}}) \cdot \mathbf{n}(\mathbf{\overline{x}}))+[u_0(\mathbf{x})]_{pp}\int_{\partial\Omega_{N\delta}} J_{\delta}(|\mathbf{x}-\mathbf{y}|)|(\mathbf{y}-\mathbf{x})\cdot\mathbf{p}(\overline{\mathbf{x}})|^2d\mathbf{y}\\
 &-[u_0(\mathbf{x})]_{pp}\int_{\partial\Omega_{N\delta}} J_{\delta}(|\mathbf{x}-\mathbf{y}|)\left(|(\mathbf{y}-\overline{\mathbf{x}})\cdot\mathbf{n}(\overline{\mathbf{x}})|^2%
 -|(\mathbf{x}-\overline{\mathbf{x}})\cdot\mathbf{n}(\overline{\mathbf{x}})|^2\right)d\mathbf{y} + O(\delta^2).
\end{align*}
With the above properties one has the following approximation via Taylor expansion:
{\begin{align}
 \nonumber&2\int_{\Omega} %
 J_{\delta}(|\mathbf{x}-\mathbf{y}|)(u_0(\mathbf{y})-u_0(\mathbf{x}))d\mathbf{y}\\
 \nonumber =&-2\int_{\Omega} J_{\delta}(|\mathbf{x}-\mathbf{y}|)[u_0(\mathbf{x})]_{nn}%
 ((\mathbf{x}-\mathbf{y})\cdot (\overline{\mathbf{x}}-\mathbf{x}))d\mathbf{y}+2\int_{E_\delta} J_{\delta}(|\mathbf{x}-\mathbf{y}|)\dfrac{\partial u_0({\mathbf{x}})}{\partial \mathbf{p}}%
 ((\mathbf{x}-\mathbf{y})\cdot \mathbf{p}(\overline{\mathbf{x}}))d\mathbf{y}\\
 \nonumber&+\int_{\Omega} J_{\delta}(|\mathbf{x}-\mathbf{y}|)[u_0(\mathbf{x})]_{nnn}%
 ((\mathbf{x}-\mathbf{y})\cdot \mathbf{n}(\overline{\mathbf{x}})) (-|\overline{\mathbf{x}}-\mathbf{x}|^2+\dfrac{1}{3}|(\mathbf{x}-\mathbf{y})\cdot \mathbf{n}(\overline{\mathbf{x}})|^2)d\mathbf{y}\\
 \nonumber&+\int_{\Omega} J_{\delta}(|\mathbf{x}-\mathbf{y}|)[u_0(\mathbf{x})]_{nn}%
 |(\mathbf{x}-\mathbf{y})\cdot \mathbf{n}(\overline{\mathbf{x}})|^2d\mathbf{y}+\int_{\Omega} J_{\delta}(|\mathbf{x}-\mathbf{y}|)[u_0(\mathbf{x})]_{pp}%
 |(\mathbf{x}-\mathbf{y})\cdot \mathbf{p}(\overline{\mathbf{x}})|^2d\mathbf{y}\\
 &+\int_{\Omega} J_{\delta}(|\mathbf{x}-\mathbf{y}|)[u_0(\mathbf{x})]_{npp}%
 ((\mathbf{x}-\mathbf{y})\cdot \mathbf{n}(\overline{\mathbf{x}}))|(\mathbf{x}-\mathbf{y})\cdot \mathbf{p}(\overline{\mathbf{x}})|^2d\mathbf{y}+O(\delta^2),\label{eqn:taylor}
\end{align}}
}
and the estimate for $T_\delta$ with $\mathbf{x}\in \Omega_{N\delta}$:
{{\begin{align}
\nonumber T_\delta =&(L_0u_0-L_{N\delta} u_0)+(f_\delta-f)\\
 \nonumber =&-\triangle u_0(\mathbf{x})+2\int_{\Omega} %
 J_{\delta}(|\mathbf{x}-\mathbf{y}|)(u_0(\mathbf{y})-u_0(\mathbf{x}))d\mathbf{y}+2M_\delta(\mathbf{x})
 \int_{-\delta}^{\delta} H(|l|) [u_0(\mathbf{x}_l)-u_0(\mathbf{x})] d \mathbf{x}_l \\
 \nonumber&-\int_{\partial\Omega_{N\delta}} J_{\delta}(|\mathbf{x}-\mathbf{y}|)%
 (|(\mathbf{y}-\overline{\mathbf{x}})\cdot\mathbf{n}(\overline{\mathbf{x}})|^2-|(\mathbf{x}-\overline{\mathbf{x}})\cdot\mathbf{n}%
 (\overline{\mathbf{x}})|^2)(-\triangle u_0(\mathbf{x}))d\mathbf{y}\\
 \nonumber =&2\int_{E_\delta} J_{\delta}(|\mathbf{x}-\mathbf{y}|)\dfrac{\partial u_0({\mathbf{x}})}{\partial \mathbf{p}}%
 ((\mathbf{x}-\mathbf{y})\cdot \mathbf{p}(\overline{\mathbf{x}}))d\mathbf{y}\\
 \nonumber&+\int_{\Omega} J_{\delta}(|\mathbf{x}-\mathbf{y}|)[u_0(\mathbf{x})]_{nnn}%
 ((\mathbf{x}-\mathbf{y})\cdot \mathbf{n}(\overline{\mathbf{x}})) (-|\overline{\mathbf{x}}-\mathbf{x}|^2+\dfrac{1}{3}|(\mathbf{x}-\mathbf{y})\cdot \mathbf{n}(\overline{\mathbf{x}})|^2)d\mathbf{y}\\
 \nonumber&+\int_{\Omega} J_{\delta}(|\mathbf{x}-\mathbf{y}|)[u_0(\mathbf{x})]_{npp}%
 ((\mathbf{x}-\mathbf{y})\cdot \mathbf{n}(\overline{\mathbf{x}}))|(\mathbf{x}-\mathbf{y})\cdot \mathbf{p}(\overline{\mathbf{x}})|^2d\mathbf{y}\\
&+{\kappa M_\delta(\mathbf{x}) [u_0(\mathbf{x})]_{nn} ((\mathbf{x} - \overline{\mathbf{x}}) \cdot \mathbf{n}(\mathbf{\overline{x}}))}+O(\delta^2).\label{eqn:taylorT}
\end{align}}}


We have then finished the proof.

\section{Acknowledgments}

Sandia National Laboratories is a multimission laboratory managed and operated by National Technology and Engineering Solutions of Sandia, LLC., a wholly owned subsidiary of Honeywell %
International, Inc., for the U.S. 555 Department of Energys National Nuclear Security Administration under contract DE-NA-0003525. H. You and Y. Yu would like to acknowledge support %
from the National Science Foundation under awards DMS 1753031. Y. Yu is also partially supported by the Lehigh faculty research grant. X.Y. Lu acknowledges the partial support of %
Lakehead University internal grants
10-50-16422410 and 10-50-16422409, and NSERC Discovery Grant 10-50-16420120.
The authors want to express their appreciation of the critical suggestions from Dr. Xiaochuan Tian, Dr. Marta D'Elia and Dr. Michael Parks, which improved the clarity and quality of this work.

\bibliographystyle{plain}
\bibliography{yyu}

\end{document}